\documentclass[a4paper]{amsart}
\usepackage[T1]{fontenc}
\usepackage[latin1]{inputenc}
\usepackage{amssymb}
\usepackage{amsfonts}
\usepackage{amsxtra}
\usepackage[all]{xy}
\usepackage{enumerate}
\usepackage{mathrsfs}

\usepackage[lite]{amsrefs}

\newcommand*{\KK}{\mathrm{KK}}
\newcommand*{\RKK}{\mathrm{RKK}}
\newcommand*{\CRKK}{\mathcal{R}\mathrm{KK}}

\newcommand*{\UCT}{\mathrm{UCT}}
\newcommand*{\K}{\mathrm{K}}
\newcommand*{\Ktop}{\mathrm{K}^{\mathrm{top}}}

\DeclareMathOperator{\ev}{ev}
\DeclareMathOperator{\Obs}{Obs}
\DeclareMathOperator{\Res}{Res}
\DeclareMathOperator{\Ind}{Ind}

\DeclareMathOperator{\Ext}{Ext}

\DeclareMathOperator{\cone}{cone}
\DeclareMathOperator{\cyl}{cyl}
\DeclareMathOperator{\const}{const}
\DeclareMathOperator{\coker}{coker}

\DeclareMathOperator{\hoinjlim}{ho-\!\varinjlim}

\DeclareMathOperator{\MCAC}{LC}

\newcommand*{\Cat}{{\mathfrak{C}}}
\newcommand*{\Tri}{{\mathcal{T}}}
\newcommand*{\CC}{{\mathcal{CC}}}
\newcommand*{\CI}{{\mathcal{CI}}}
\newcommand*{\CIZ}{{\mathcal{CI}_0}}
\newcommand*{\Prfin}{{\mathcal{CIF}}}
\newcommand*{\Ab}{{\mathsf{Ab}}}

\newcommand*{\pt}{{\star}}

\newcommand*{\EG}{{\mathcal{E}G}}
\newcommand*{\Dirac}{{\mathsf{D}}}
\newcommand*{\ADir}{{\mathsf{P}}}
\newcommand*{\AN}{{\mathsf{N}}}
\newcommand*{\Proj}{{\mathcal{P}}}
\newcommand*{\Null}{{\mathcal{N}}}
\newcommand*{\RC}{\mathsf{R}}

\newcommand*{\ID}{{\mathrm{id}}}

\newcommand*{\Mult}{{\mathcal{M}}}
\newcommand*{\Cred}{C^*_{\mathrm{r}}}

\newcommand*{\C}{{\mathbb{C}}}
\newcommand*{\R}{{\mathbb{R}}}
\newcommand*{\Z}{{\mathbb{Z}}}
\newcommand*{\Ztwo}{{\mathbb{Z}/2}}
\newcommand*{\N}{{\mathbb{N}}}
\newcommand*{\Comp}{{\mathbb{K}}}

\newcommand*{\brd}{-\hspace{0pt}}
\newcommand*{\nbd}{\nobreakdash-\hspace{0pt}}

\newcommand*{\abs}[1]{\lvert#1\rvert}

\newcommand*{\gen}[1]{\langle#1\rangle}
\newcommand*{\GEN}{\mathcal{G}}

\newcommand*{\rcross}{\mathbin{\ltimes_{\mathrm{r}}}}
\newcommand*{\cross}{\mathbin{\ltimes}}
\newcommand*{\blank}{{\llcorner\!\!\lrcorner}}
\newcommand*{\defeq}{\mathrel{:=}}
\newcommand*{\congto}{\mathrel{\overset{\cong}{\to}}}
\newcommand*{\into}{\rightarrowtail}
\newcommand*{\prto}{\twoheadrightarrow}

\newcommand*{\Left}{\mathbb{L}}

\newcommand*{\Lotimes}{\mathbin{\otimes^{\mathbb{L}}}}
\newcommand*{\Lrcross}{\mathbin{\ltimes^{\mathbb{L}}_{\mathrm{r}}}}
\newcommand*{\Obscross}{\mathbin{\ltimes^{\mathrm{Obs}}_{\mathrm{r}}}}
\newcommand*{\Lcross}{\mathbin{\ltimes^{\mathbb{L}}}}

\theoremstyle{plain}
\newtheorem{theorem}{Theorem}[section]
\newtheorem{proposition}[theorem]{Proposition}
\newtheorem{lemma}[theorem]{Lemma}
\newtheorem{corollary}[theorem]{Corollary}

\theoremstyle{definition}
\newtheorem{definition}[theorem]{Definition}
\theoremstyle{remark}
\newtheorem{remark}[theorem]{Remark}

\begin{document}

\title[The Baum-Connes Conjecture via Localisation of Categories]{The
  Baum-Connes Conjecture via\\ Localisation of Categories}

\author{Ralf Meyer}
\address{Mathematisches Institut\\
         Westfälische Wilhelms-Universität Münster\\
         Einsteinstr.\ 62\\
         48149 Münster\\
         Germany
}
\email{rameyer@math.uni-muenster.de}

\author{Ryszard Nest}
\address{Københavns Universitets Institut for Matematiske Fag\\
         Universitetsparken 5\\ 2100 København\\ Denmark
}
\email{rnest@math.ku.dk}


\thanks{This research was supported by the EU-Network \emph{Quantum
    Spaces and Noncommutative Geometry} (Contract HPRN-CT-2002-00280)
  and the \emph{Deutsche Forschungsgemeinschaft} (SFB 478).}

\begin{abstract}
  We redefine the Baum-Connes assembly map using simplicial approximation in
  the equivariant Kasparov category.  This new interpretation is ideal for
  studying functorial properties and gives analogues of the Baum-Connes
  assembly map for other equivariant homology theories.  We extend many of the
  known techniques for proving the Baum-Connes conjecture to this more general
  setting.
\end{abstract}
\maketitle

\tableofcontents

\section{Introduction}
\label{sec:intro}

Let~$G$ be a second countable locally compact group.  Let~$A$ be a separable
$C^*$\nbd{}algebra with a strongly continuous action of~$G$ and let $G\rcross
A$ be the reduced crossed product, which is another separable
$C^*$\nbd{}algebra.  The aim of the Baum-Connes conjecture (with coefficients)
is to compute the $\K$\nbd{}theory of $G\rcross A$.  For the trivial action
of~$G$ on~$\C$ (or~$\R$), this specialises to $\K_*(\Cred(G))$, the
$\K$\nbd{}theory of the reduced $C^*$\nbd{}algebra of~$G$.  One defines a
certain graded Abelian group $\Ktop_*(G,A)$, called the \emph{topological
  $\K$\nbd{}theory of~$G$ with coefficients~$A$}, and a homomorphism
\begin{equation}  \label{eq:BC_assembly}
  \mu_A\colon \Ktop_*(G,A)\to\K_*(G\rcross A),
\end{equation}
which is called the \emph{Baum-Connes assembly map}.  The \emph{Baum-Connes
  conjecture for~$G$ with coefficients~$A$} asserts that this map is an
isomorphism.  It has important applications in topology and ring theory.  The
conjecture is known to hold in many cases, for instance, for amenable groups
(\cite{Higson-Kasparov:Amenable}).  A recent survey article on the Baum-Connes
conjecture is~\cite{Higson:Survey_Article}.

Despite its evident success, the usual definition of the Baum-Connes assembly
map has some important shortcomings.  At first sight $\Ktop_*(G,A)$ may seem
even harder to compute than $\K_*(G\rcross A)$.  Experience shows that this is
not the case.  Nevertheless, there are situations where $\Ktop_*(G,A)$ creates
more trouble than $\K_*(G\rcross A)$.  For instance, most of the work required
to prove the permanence properties of the Baum-Connes conjecture is needed to
extend evident properties of $\K_*(G\rcross A)$ to $\Ktop_*(G,A)$.  The
meaning of the Baum-Connes conjecture is rather mysterious: it is not \emph{a
priori} clear that $\Ktop_*(G,A)$ should have anything to do with
$\K_*(G\rcross A)$.  A related problem is that the Baum-Connes assembly map
only makes sense for $\K$\nbd{}theory and not for other interesting
equivariant homology theories.  For instance, in connection with the Chern
character it would be desirable to have a Baum-Connes assembly map for local
cyclic homology as well.

Our alternative description of the assembly map addresses these shortcomings.
It applies to any equivariant homology theory, that is, any functor defined on
the equivariant Kasparov category $\KK^G$.  For instance, we can also apply
$\K$\nbd{}homology and local cyclic homology to the crossed product.
Actually, this is nothing so new.  Gennadi Kasparov did this using his Dirac
dual Dirac method---for all groups to which his method applies (see
\cites{Kasparov:Novikov, Kasparov-Skandalis:Buildings}).  In his approach, the
topological side of the Baum-Connes conjecture appears as the
$\gamma$\nbd{}part of $\K_*(G\rcross A)$, and this $\gamma$\nbd{}part makes
sense for any functor defined on $\KK^G$.  Indeed, our approach is very close
to Kasparov's.  We show that one half of Kasparov's method, namely, the Dirac
morphism, exists in complete generality, and we observe that this suffices to
construct the assembly map.  From the technical point of view, this is the
main innovation in this article.

Our approach is very suitable to state and prove general functorial properties
of the assembly map.  The various known permanence results of the Baum-Connes
conjecture become rather transparent in our setup.  Such permanence results
have been investigated by several authors.  There is a series of papers by
Jérôme Chabert, Siegfried Echterhoff and Hervé Oyono-Oyono
(\cites{Chabert:BC_product, Chabert-Echterhoff:Permanence,
  Chabert-Echterhoff:Twisted, Chabert-Echterhoff-Oyono:Shapiro,
  Chabert-Echterhoff-Oyono:Going_down}).  Both authors of this article have
been quite familiar with their work, and it has greatly influenced this
article.  We also reprove a permanence result for unions of groups by Paul
Baum, Stephen Millington and Roger Plymen (\cite{Baum-Millington-Plymen}) and
a result relating the real and complex versions of the Baum-Connes conjecture
by Paul Baum and Max Karoubi (\cite{Baum-Karoubi}) and independently by Thomas
Schick (\cite{Schick:Real}).  In addition, we use results
of~\cite{Schick:Real} to prove that the existence of a $\gamma$\nbd{}element
for a group~$G$ for real and complex coefficients is equivalent.

A good blueprint for our approach towards the Baum-Connes conjecture is the
work of James Davis and Wolfgang Lück in~\cite{Davis-Lueck:Assembly}.  As
kindly pointed out by the referee, the approach of Paul Balmer and Michel
Matthey in \cites{Balmer-Matthey:Foundations, Balmer-Matthey:Cofibrant,
  Balmer-Matthey:Model_BC} is even closer.  However, these are only formal
analogies, as we shall explain below.

Davis and Lück only consider discrete groups and reinterpret the Baum-Connes
assembly map for $\K_*(G\rcross C_0(X))$ as follows.  A proper
$G$\nbd{}CW\brd{}complex~$\tilde{X}$ with a $G$\nbd{}equivariant continuous
map $\tilde{X}\to X$ is called a \emph{proper $G$\nbd{}CW\brd{}approximation}
for~$X$ if it has the following universal property: any map from a proper
$G$\nbd{}CW\brd{}complex to~$X$ factors through~$\tilde{X}$, and this
factorisation is unique up to equivariant homotopy.  Such approximations
always exist and are unique up to equivariant homotopy equivalence.  Given a
functor~$F$ on the category of $G$\nbd{}spaces, one defines its localisation
by $\Left F(X)\defeq F(\tilde{X})$ (up to isomorphism).  It comes equipped
with a map $\Left F(X)\to F(X)$.  For suitable~$F$, this is the Baum-Connes
assembly map.

We replace the homotopy category of $G$\nbd{}spaces by the
$G$\nbd{}equivariant Kasparov category $\KK^G$, whose objects are the
separable $G$-$C^*$\brd{}algebras and whose morphism spaces are the bivariant
groups $\KK^G_0(A,B)$ defined by Kasparov.  We need some extra structure, of
course, in order to do algebraic topology.  For our purposes, it is enough to
turn $\KK^G$ into a triangulated category (see \cites{Neeman:Triangulated,
Verdier:Thesis}).  The basic examples of triangulated categories are the
derived categories in homological algebra and the stable homotopy category in
algebraic topology.  They have enough structure to localise and to do
rudimentary homological algebra.  According to our knowledge, Andreas Thom's
thesis~\cite{Thom:Thesis} is the first work on $C^*$\nbd{}algebras where
triangulated categories are explicitly used.  Since this structure is crucial
for us and not well-known among operator algebraists, we discuss it in an
operator algebraic context in Section~\ref{sec:triangulated_categories}.  We
also devote an appendix to a detailed proof that $\KK^G$ is a triangulated
category.  This verification of axioms is not very illuminating.  The reason
for including it is that we could not find a good reference.

We call $A\in\KK^G$ \emph{compactly induced} if it is $\KK^G$\brd{}equivalent
to $\Ind_H^G A'$ for some compact subgroup $H\subseteq G$ and some
$H$-$C^*$\brd{}algebra~$A'$.  We let $\CI\subseteq\KK^G$ be the full
subcategory of compactly induced objects and $\gen{\CI}$ the localising
subcategory generated by it.  The objects of $\gen{\CI}$ are our substitute
for proper $G$\nbd{}CW\brd{}complexes.  The objects of $\CI$ behave like the
cells out of which proper $G$\nbd{}CW\brd{}complexes are built.  We define a
\emph{$\CI$\brd{}simplicial approximation} of $A\in\KK^G$ as a morphism
$\tilde{A}\to A$ in $\KK^G$ with $\tilde{A}\in\gen{\CI}$ such that
$\KK^G(P,\tilde{A})\cong\KK^G(P,A)$ for all $P\in\gen{\CI}$.  We show that
\emph{$\CI$\brd{}simplicial approximations} always exist, are unique,
functorial, and have good exactness properties.  Therefore, if
$F\colon\KK^G\to\Cat$ is any homological functor into an Abelian category,
then its \emph{localisation} $\Left F(A)\defeq F(\tilde{A})$ is again a
homological functor $\KK^G\to\Cat$.  It comes equipped with a natural
transformation $\Left F(A)\to F(A)$.  For the functor $F(A)\defeq
\K_*(G\rcross A)$, this map is naturally isomorphic to the Baum-Connes
assembly map.  In particular, $\Ktop_*(G,A)\cong\Left F(A)$.  Thus we have
redefined the Baum-Connes assembly map as a localisation.

Of course, we do not expect the map $\Left F(A)\to F(A)$ to be an isomorphism
for all functors~$F$.  For instance, consider the $\K$\nbd{}theory of the full
and reduced crossed products.  We will show that both functors have the same
localisation.  However, the full and reduced group $C^*$\nbd{}algebras may
have different $\K$\nbd{}theory.

A variant of Green's Imprimitivity Theorem (\cite{Green:Imprimitivity}) for
reduced crossed products says that $G\rcross \Ind_H^G A$ for a compact
subgroup $H\subseteq G$ is Morita-Rieffel equivalent to $H\cross A$.
Combining this with the Green-Julg Theorem (\cite{Julg:Green-Julg}), we get
$$
\K_*(G\rcross \Ind_H^G A)\cong \K_*(H\cross A) \cong \K^H_*(A).
$$
Hence $\K_*(G\rcross B)$ is comparatively easy to compute for $B\in\CI$.  For
an object of $\gen{\CI}$, we can, in principle, compute its $\K$\nbd{}theory
by decomposing it into building blocks from $\CI$.  In a forthcoming article,
we will discuss a spectral sequence that organises this computation.  As a
result, $\K_*(G\rcross \tilde{A})$ is quite tractable for
$\tilde{A}\in\gen{\CI}$.  The $\CI$\brd{}simplicial approximation replaces an
arbitrary coefficient algebra~$A$ by the best approximation to~$A$ in this
tractable subcategory in the hope that $\K_*(G\rcross\tilde{A})\cong
\Ktop_*(G,A)$ is then a good approximation to $\K_*(G\rcross A)$.

Above we have related the Baum-Connes assembly map to simplicial approximation
in homotopy theory.  Alternatively, we can use an analogy to homological
algebra.  In this picture, the category $\KK^G$ corresponds to the homotopy
category of chain complexes over an Abelian category.  The latter has chain
complexes as objects and homotopy classes of chain maps as morphisms.  To do
homological algebra, we also need exact chain complexes and
quasi-isomorphisms.  In our context, these have the following analogues.

A $G$\nbd{}$C^*$\brd{}algebra is called \emph{weakly contractible} if it is
$\KK^H$\nbd{}equivalent to~$0$ for all compact subgroups $H\subseteq G$.  We
let $\CC\subseteq\KK^G$ be the full subcategory of weakly contractible
objects.  This is a localising subcategory of $\KK^G$.  We call
$f\in\KK^G(A,B)$ a \emph{weak equivalence} if it is invertible in $\KK^H(A,B)$
for all compact subgroups ${H\subseteq G}$.  The weakly contractible objects
and the weak equivalences determine each other: a morphism is a weak
equivalence if and only if its ``mapping cone'' is weakly contractible,
and~$A$ is weakly contractible if and only if the zero map $0\to A$ is a weak
equivalence.

The subcategories $\CC$ and $\gen{\CI}$ are ``orthogonal complements'' in the
sense that $B\in\CC$ if and only if $\KK^G(A,B)=0$ for all $A\in\gen{\CI}$,
and $A\in\gen{\CI}$ if and only if $\KK^G(A,B)=0$ for all $B\in\CC$.  Hence
$f\in\KK^G(B,B')$ is a weak equivalence if and only if the induced map
$\KK^G(A,B)\to\KK^G(A,B')$ is an isomorphism for all $A\in\gen{\CI}$.
Therefore, a $\CI$\nbd{}simplicial approximation for~$A$ is the same as a weak
equivalence $f\in\KK^G(\tilde{A},A)$ with $\tilde{A}\in\gen{\CI}$.

We now return to our analogy with homological algebra.  The weakly
contractible objects play the role of the exact chain complexes and the weak
equivalences play the role of the quasi-isomorphisms.  Objects of $\gen{\CI}$
correspond to projective chain complexes as defined in~\cite{Keller:Handbook}.
Hence $\CI$\brd{}simplicial approximations correspond to projective
resolutions.  In homological algebra, we can compute the total left derived
functor of a functor~$F$ by applying~$F$ to a projective resolution.  Thus
$\Left F$ as defined above corresponds to the total left derived functor
of~$F$.  In particular, $\Ktop_*(G,A)$ appears as the total left derived
functor of $\K_*(G\rcross A)$.

Bernhard Keller's presentation of homological algebra
in~\cite{Keller:Handbook} is quite close to our constructions because it
relies very much on triangulated categories.  This is unusual because most
authors prefer to use the finer structure of Abelian categories.  However,
nothing in our setup corresponds to the underlying Abelian category.  Hence we
only get an analogue of the total derived functor, not of the satellite
functors that are usually called derived functors.  A more serious difference
is that there are almost no interesting exact functors in homological algebra.
In contrast, the Baum-Connes conjecture asserts that the functor
$\K_*(G\rcross A)$ agrees with its total derived functor, which is equivalent
to exactness in classical homological algebra.  Hence the analogy to
homological algebra is somewhat misleading.

Using weak equivalences, we can also formulate the Baum-Connes conjecture with
coefficients as a \emph{rigidity} statement.  The assembly map $\Left F(A)\to
F(A)$ is an isomorphism for all~$A$ if and only if~$F$ maps all weak
equivalences to isomorphisms.  If~$F$ satisfies some exactness property, this
is equivalent to $F(A)=0$ for all $A\in\CC$.  If $A\in\CC$, then~$A$ is
$\KK^G$\brd{}equivalent to a $G$-$C^*$\brd{}algebra that is
$H$\nbd{}equivariantly contractible for any compact subgroup $H\subseteq G$
(replace~$A$ by the universal algebra $q_s A$ defined in~\cite{Meyer:KKG}).
Thus the Baum-Connes conjecture with coefficients is equivalent to the
statement that $\K_*(G\rcross A)=0$ if~$A$ is $H$\nbd{}equivariantly
contractible for all compact subgroups $H\subseteq G$.  Another equivalent
formulation that we obtain in Section~\ref{sec:strong_BC} is the following.
The Baum-Connes conjecture with coefficients is equivalent to the statement
that $\K_*(G\rcross A)=0$ if $\K_*(H\cross A)=0$ for all compact subgroups
$H\subseteq G$.  Both reformulations of the Baum-Connes conjecture with
coefficients are as elementary as possible: they involve nothing more than
compact subgroups, $\K$\nbd{}theory and reduced crossed products.

The localisation of the homotopy category of chain complexes over an Abelian
category at the subcategory of exact chain complexes, is its derived category,
which is the category of primary interest in homological algebra.  In our
context, it corresponds to the localisation $\KK^G/\CC$.  We describe
$\KK^G/\CC$ in more classical terms, using the universal proper $G$\nbd{}space
$\EG$.  We identify the space of morphisms $A\to B$ in $\KK^G/\CC$ with the
group $\RKK^G(\EG;A,B)$ as defined by Kasparov (\cite{Kasparov:Novikov}).  The
canonical functor $\KK^G\to\KK^G/\CC$ is the obvious one,
$$
p_\EG^*\colon \KK^G(A,B)\to\RKK^G(\EG;A,B).
$$
As a consequence, if~$A$ is weakly contractible, then $p^*_Y(A)\cong0$ for any
proper $G$\nbd{}space~$Y$.  This means that the homogeneous spaces $G/H$ for
$H\subseteq G$ compact, which are implicitly used in the definition of weak
contractibility, already generate all proper $G$\nbd{}spaces. Another
consequence is that proper $G$\nbd{}$C^*$\brd{}algebras in the sense of
Kasparov belong to $\gen{\CI}$.  Conversely, for many groups any object of
$\gen{\CI}$ is $\KK^G$\brd{}equivalent to a proper $G$-$C^*$\brd{}algebra (see
the end of Section~\ref{sec:derived_proper}).

Let $\pt\in\KK^G$ be the real or complex numbers, depending on the category we
work with.  We have a tensor product operation in $\KK^G$, which is nicely
compatible with the subcategories $\CC$ and $\CI$.  Therefore, if
$\Dirac\in\KK^G(\ADir,\pt)$ is a $\CI$\brd{}simplicial approximation
for~$\pt$, then $\Dirac\otimes\ID_A\in\KK^G(\ADir\otimes A,A)$ is a
$\CI$\brd{}simplicial approximation for $A\in\KK^G$.  Thus we can describe the
localisation of a functor more explicitly as $\Left F(A)\defeq F(\ADir\otimes
A)$.  We call~$\Dirac$ a \emph{Dirac morphism} for~$G$.  Its existence is
equivalent to the representability of a certain functor.  Eventually, this is
deduced from a generalisation of Brown's Representability Theorem to
triangulated categories.

The following example of a Dirac morphism motivates our terminology.  Suppose
that $\EG$ is a smooth manifold.  Replacing it by a suspension of $T^*\EG$, we
achieve that $\EG$ has a $G$\nbd{}invariant spin structure and that $8 \mid
\dim \EG$.  Then the Dirac operator on $\EG$ defines an element of
$\KK^G_0(C_0(\EG),\pt)$; this is a Dirac morphism for~$G$.  We can also
describe it as the element $p!\in
\KK^G_0(C_0(\EG),\pt)$ associated to the constant map $p\colon \EG\to\pt$ by
wrong way functoriality.  Unfortunately, wrong way functoriality only works
for manifolds.  Extending it to non-Hausdorff manifolds as
in~\cite{Kasparov-Skandalis:Buildings}, one can construct explicit Dirac
morphisms also for groups acting properly and simplicially on finite
dimensional simplicial complexes.  However, it is unclear how to adapt this to
infinite dimensional situations.

Since we work in the Kasparov category, Bott periodicity is an integral part
of our setup.  The above example of a Dirac morphism shows that wrong way
functoriality and hence Bott periodicity indeed play significant roles.  This
distinguishes our approach from \cites{Balmer-Matthey:Foundations,
Balmer-Matthey:Cofibrant, Balmer-Matthey:Model_BC, Davis-Lueck:Assembly}.  The
bad news is that we cannot treat homology theories such as algebraic
$\K$-theory that do not satisfy periodicity.  The good news is that the Dirac
dual Dirac method, which is one of the main proof techniques in connection
with the Baum-Connes conjecture, is also part of our setup.  In examples, this
method usually arises as an equivariant version of Bott periodicity.

A \emph{dual Dirac morphism} is an element $\eta\in\KK^G(\pt,\ADir)$ that is a
left-inverse to the Dirac morphism $\Dirac\in\KK^G(\ADir,\pt)$, that is,
$\eta\Dirac=\ID_\ADir$.  Suppose that it exists.  Then $\gamma=\Dirac\eta$ is
an idempotent in $\KK^G(\pt,\pt)$.  By exterior product, we get idempotents
$\gamma_A\in\KK^G(A,A)$ for all $A\in\KK^G$.  We have $A\in\CC$ if and only if
$\gamma_A=0$, and $A\in\gen{\CI}$ if and only if $\gamma_A=1$.  The category
$\KK^G$ is equivalent to the direct product $\KK^G\cong\CC\times\gen{\CI}$.
Therefore, the assembly map is split injective for any covariant functor.  For
groups with the Haagerup property and, in particular, for amenable groups, a
dual Dirac morphism exists and we have $\gamma=1$.  This important theorem is
due to Nigel Higson and Gennadi Kasparov (\cite{Higson-Kasparov:Amenable}).
In this case, weak equivalences are already isomorphisms in $\KK^G$.  Hence
$\Left F=F$ for any functor~$F$.

When we compose two functors in homological algebra, it frequently happens
that $\Left (F'\circ F)\cong\Left F'\circ\Left F$.  This holds, for instance,
if~$F$ maps projectives to projectives.  We check that the restriction and
induction functors preserve the subcategories $\CC$ and $\gen{\CI}$.  The same
holds for the complexification functor from real to complex $\KK$\brd{}theory
and many others.  The ensuing identities of localised functors imply
permanence properties of the Baum-Connes conjecture.

Another useful idea that our new approach allows is the following.  Instead of
deriving the functor $A\mapsto \K_*(G\rcross A)$, we may also derive the
crossed product functor $A\mapsto G\rcross A$ itself.  Its localisation
$G\Lrcross A$ is a triangulated functor from $\KK^G$ to $\KK$.  It can be
described explicitly as $G\Lrcross A=G\rcross (\ADir\otimes A)$ if
$\Dirac\in\KK^G(\ADir,\pt)$ is a Dirac morphism.  The Baum-Connes conjecture
asks for $\Dirac_*\in\KK(G\Lrcross A,G\rcross A)$ to induce an isomorphism on
$\K$\nbd{}theory.  Instead, we can ask it to be a $\KK$\brd{}equivalence.
Then the Baum-Connes conjecture holds for $F(G\rcross A)$ for any split exact,
stable homotopy functor~$F$ on $C^*$\nbd{}algebras because such functors
descend to the category $\KK$.  For instance, this covers local cyclic
(co)homology and $\K$\nbd{}homology.

This stronger conjecture is known to be false in some cases where the
Baum-Connes conjecture holds.  Nevertheless, it holds in many examples.  For
groups with the Haagerup property, we have $\gamma=1$, so that $\Left F=F$ for
any functor, anyway.  If both $G\rcross A$ and $G\Lrcross A$ satisfy the
Universal Coefficient Theorem (UCT) in $\KK$, then an isomorphism on
$\K$\nbd{}theory is automatically a $\KK$\brd{}equivalence.  Since
$G\Lrcross\pt$ always satisfies the UCT, the strong Baum-Connes conjecture
with trivial coefficients holds if and only if the usual Baum-Connes
conjecture holds and $\Cred(G)$ satisfies the UCT.  This is known to be the
case for almost connected groups and linear algebraic groups over
$p$\nbd{}adic number fields, see
\cites{Chabert-Echterhoff-Oyono:Going_down,
Chabert-Echterhoff-Nest:Connes_Kasparov}.

This article is the first step in a programme to extend the Baum-Connes
conjecture to quantum group crossed products.  It does not seem a good idea to
extend the usual construction in the group case because it is not clear
whether the resulting analogue of $\Ktop_*(G,A)$ can be computed.  Even if we
had a good notion of a proper action of a quantum group, these actions would
certainly occur on very noncommutative spaces, so that we have to ``quantise''
the algebraic topology needed to compute $\Ktop_*(G,A)$.  The framework of
triangulated categories and localisation of functors is ideal for this
purpose.  In the group case, the homogeneous spaces $G/H$ for compact
subgroups $H\subseteq G$ generate all proper actions in the group case.  Thus
we expect that we can formulate the Baum-Connes conjecture for quantum groups
using quantum homogeneous spaces instead of proper actions.  However, we still
need some further algebraic structure: restriction and induction functors and
tensor products of coactions.  We plan to treat this additional structure and
to construct a Baum-Connes assembly map for quantum groups in a sequel to this
paper.  Here we only consider the classical case of group actions.

\subsection{Some general conventions}
\label{sec:conventions}

Let~$\Cat$ be a category.  We write $A\in\Cat$ to denote that~$A$ is an object
of~$\Cat$, and $\Cat(A,B)$ for the space of morphisms $A\to B$ in~$\Cat$.

It makes no difference whether we work with real, ``real'', or complex
$C^*$\nbd{}algebras.  Except for section~\ref{sec:real_complex}, we do not
distinguish between these cases in our notation.  Of course, standard
$C^*$\nbd{}algebras like $C_0(X)$ and $\Cred(G)$ have to be taken in the
appropriate category.  We denote the one-point space by~$\pt$ and also write
$\pt=C(\pt)$.  Thus~$\pt$ denotes the complex or real numbers depending on the
category we use.

Locally compact groups and spaces are tacitly assumed to be second countable,
and $C^*$\brd{}algebras are tacitly assumed to be separable.  Let~$G$ be a
locally compact group and let~$X$ be a locally compact $G$\nbd{}space.  A
\emph{$G$-$C^*$\brd{}algebra} is a $C^*$\nbd{}algebra equipped with a strongly
continuous action of~$G$ by automorphisms.  A \emph{$G\cross
  X$-$C^*$\brd{}algebra} is a $G$-$C^*$\brd{}algebra equipped with a
$G$\nbd{}equivariant essential $*$\nbd{}homomorphism from $C_0(X)$ to the
centre of its multiplier algebra.  Kasparov defines bivariant $\K$\nbd{}theory
groups $\CRKK^G_*(X;A,B)$ involving these data in
\cite{Kasparov:Novikov}*{Definition 2.19}.  The notation $\CRKK^G_*(X;A,B)$
should be distinguished from $\RKK^G_*(X;A,B)$.  The latter is defined for two
$G$-$C^*$-algebras $A$ and~$B$ by
\begin{equation}  \label{eq:RKKG}
  \RKK^G_*(X;A,B) \defeq \CRKK^G_*(X;C_0(X,A),C_0(X,B)).
\end{equation}

Since $\CRKK^G_*(X;A,B)$ agrees with the bivariant $\K$\nbd{}groups for the
groupoid $G\cross X$ as defined in~\cite{LeGall:KK_groupoid}, we denote it by
$\KK^{G\cross X}_*(A,B)$.  For several purposes, it is useful to generalise
from groups to groupoids.  However, we do not treat arbitrary groupoids
because it is not so clear what should correspond to the compact subgroups in
this case.  We work with transformation groups throughout because this
generalisation is not more difficult than the group case and useful for
several applications.

We write $\K_*(A)$ for the graded Abelian group $n\mapsto \K_n(A)$, $n\in\Z$,
and similarly for $\KK^{G\cross X}_*(A,B)$.  We usually omit the
subscript~$0$, that is, $\K(A)\defeq \K_0(A)$, etc..

The \emph{$G\cross X$\nbd{}equivariant Kasparov category} is the additive
category whose objects are the $G\cross X$-$C^*$\brd{}algebras and whose group
of morphisms $A\to B$ is $\KK^{G\cross X}_0(A,B)$; the composition is the
Kasparov product.  We denote this category by $\KK^{G\cross X}$.

The notion of equivalence for $G\cross X$-$C^*$-algebras that we encounter
most frequently is $\KK$-equivalence, that is, isomorphism in $\KK^{G\cross
  X}$, which we simply denote by ``$\cong$''.  Sometimes we may want to stress
that two $G\cross X$-$C^*$-algebras are more than just $\KK$\nbd{}equivalent.
We write $A\approx B$ if $A$ and~$B$ are isomorphic as $G\cross
X$-$C^*$-algebras and $A\sim_M B$ if $A$ and~$B$ are $G\cross X$-equivariantly
Morita-Rieffel equivalent.  Both relations imply $A\cong B$.

\section{Triangulated categories of operator algebras}
\label{sec:triangulated_categories}

In this section, we explain triangulated categories in the context of
equivariant Kasparov theory.  The purpose is to introduce operator algebraists
to the language of triangulated categories, which we are using throughout this
article.  We hope that it allows them to understand this article without
having to read the specialised literature on triangulated categories (like
\cites{Neeman:Triangulated, Verdier:Thesis}).  Thus we translate various known
results of non-commutative topology into the language of triangulated
categories.  In addition, we sketch how to prove basic facts about
localisation of triangulated categories in the special case where there are
enough projectives.  There is nothing essentially new in this section.  The
only small exception is the rather satisfactory treatment of inductive limits
of $C^*$\nbd{}algebras in Section~\ref{sec:homotopy_limits}.

The two motivating examples of triangulated categories are the stable homotopy
category from algebraic topology and the derived categories of Abelian
categories from homological algebra.  The definition of a triangulated
category formalises some important structure that is present in these
categories.  The additional structure consists of a \emph{translation
  automorphism} and a class of \emph{exact triangles} (often called
\emph{distinguished triangles}).  We first explain what these are for
$\KK^{G\cross X}$.

\subsection{Suspensions and mapping cones}
\label{sec:suspend_cone}

Let $\Sigma\colon \KK^{G\cross X}\to\KK^{G\cross X}$ be the suspension functor
$\Sigma A\defeq C_0(\R)\otimes A$.  This is supposed to be the translation
automorphism in our case.  However, it is only an equivalence and not an
isomorphism of categories.  This defect is repaired by the following trick.
We replace $\KK^{G\cross X}$ by the category $\widetilde{\KK}{}^{G\cross X}$
whose objects are pairs $(A,n)$ with $A\in\KK^{G\cross X}$, $n\in\Z$, with
morphisms
$$
\widetilde{\KK}{}^{G\cross X}\bigl((A,n),(B,m)\bigr) \defeq
\varinjlim_{p\in\N} \KK^{G\cross X}(\Sigma^{p+n} A,\Sigma^{p+m} B).
$$
Actually, since the maps $\KK^{G\cross X}(A,B)\to\KK^{G\cross X}(\Sigma
A,\Sigma B)$ are isomorphisms by Bott periodicity, we can omit the direct
limit over~$p$.  Morphisms in $\widetilde{\KK}{}^{G\cross X}$ are composed in
the obvious fashion.  We define the \emph{translation} or \emph{suspension}
automorphism on $\widetilde{\KK}{}^{G\cross X}$ by $\Sigma (A,n)\defeq
(A,n+1)$.  The evident functor $\KK^{G\cross X}\to\widetilde{\KK}{}^{G\cross
  X}$, $A\mapsto (A,0)$, identifies $\KK^{G\cross X}$ with a full subcategory
of $\widetilde{\KK}{}^{G\cross X}$.  Any object of $\widetilde{\KK}{}^{G\cross
  X}$ is isomorphic to one from this subcategory because Bott periodicity
yields $(A,n)\cong (\Sigma^{n\bmod 8} A,0)$ for all $n\in\Z$,
$A\in\KK^{G\cross X}$.  Thus the categories $\widetilde{\KK}{}^{G\cross X}$
and $\KK^{G\cross X}$ are equivalent.  It is not necessary to distinguish
between $\widetilde{\KK}{}^{G\cross X}$ and $\KK^{G\cross X}$ except in very
formal arguments and definitions.  Most of the time, we ignore the difference
between these two categories.

Let $f\colon A\to B$ be an equivariant $*$\nbd{}homomorphism.  Then
its \emph{mapping cone}
\begin{equation}  \label{eq:def_cone}
  \cone(f)\defeq
  \{(a,b)\in A\times C_0(\mathopen]0,1],B) \mid f(a)=b(1)\}
\end{equation}
is again a $G\cross X$-$C^*$\brd{}algebra and there are natural
equivariant $*$\nbd{}homomorphisms
\begin{equation}  \label{eq:standard_triangle}
  \Sigma B \overset{\iota}\longrightarrow
  \cone(f) \overset{\epsilon}\longrightarrow
  A \overset{f}\longrightarrow
  B.
\end{equation}
Such diagrams are called \emph{mapping cone triangles}.  A diagram $\Sigma
B'\to C'\to A'\to B'$ in $\widetilde{\KK}{}^{G\cross X}$ is called an
\emph{exact triangle} if it is isomorphic to a mapping cone triangle.  That
is, there is an equivariant $*$\nbd{}homomorphism $f\colon A\to B$ and a
commutative diagram
$$
\xymatrix{
  {\Sigma B} \ar[r] \ar[d]^{\Sigma\beta}_{\cong} &
  \cone(f) \ar[r] \ar[d]^{\gamma}_{\cong} &
  A \ar[r] \ar[d]^{\alpha}_{\cong} &
  B \ar[d]^{\beta}_{\cong} \\
  {\Sigma B'} \ar[r] &
  {C'} \ar[r] &
  {A'} \ar[r] &
  {B'}
}
$$
where $\alpha,\beta,\gamma$ are isomorphisms in $\widetilde{\KK}{}^{G\cross
  X}$ and $\Sigma\beta$ is the suspension of~$\beta$.

\begin{proposition}  \label{pro:KKG_triangulated}
  The category $\widetilde{\KK}{}^{G\cross X}$ with~$\Sigma^{-1}$ as
  translation functor and with the exact triangles as described above
  is a triangulated category.  It has countable direct sums: they are the
  usual $C^*$\nbd{}direct sums.
\end{proposition}

It is proven in the appendix that $\widetilde{\KK}{}^{G\cross X}$ is
triangulated.  It is shown in~\cite{Kasparov:Novikov} that $\KK^{G\cross
X}(\bigoplus A_n,B)\cong\prod \KK^{G\cross X}(A_n,B)$.  This means that the
usual $C^*$\nbd{}direct sum is a direct sum operation also in $\KK^{G\cross
X}$.

Notice that the translation functor is the inverse of the suspension~$\Sigma$.
The reason for this is as follows.  The axioms of a triangulated category are
modelled after the stable homotopy category, and the functor from spaces to
$C^*$\nbd{}algebras is contravariant.  Hence we ought to work with the
opposite category of $\KK^{G\cross X}$.  The opposite category of a
triangulated category becomes again triangulated if we use ``the same''
exact triangles and replace the translation functor by its inverse.
Since we want to work with $\KK^{G\cross X}$ and not its opposite and retain
the usual constructions from the stable homotopy category, we sometimes
deviate in our conventions from the usual ones for a triangulated category.
For instance, we always write exact triangles in the form $\Sigma B\to
C\to A\to B$.

One of the axioms of a triangulated category requires any $f\in\KK^{G\cross
  X}(A,B)$ to be part of an exact triangle $\Sigma B\to C\to
A\overset{f}\to B$.  We call this triangle a \emph{mapping cone triangle}
for~$f$ and~$C$ a \emph{mapping cone} for~$f$.  We can use the mapping cone
triangle~\eqref{eq:standard_triangle} if~$f$ is an equivariant
$*$\nbd{}homomorphism.  In general, we replace~$f$ by an equivariant
$*$\nbd{}isomorphism $f'\colon q_s A\to q_s B$ and then take a mapping cone
triangle of~$f'$ as in~\eqref{eq:standard_triangle}.  The universal
$C^*$\nbd{}algebra $q_s A$ is defined in~\cite{Meyer:KKG}.  It is important
that $q_s A$ is isomorphic to~$A$ in $\KK^{G\cross X}$.  We warn the reader
that the above construction only works for ungraded $C^*$\nbd{}algebras.  For
this reason, the Kasparov category of graded $C^*$\nbd{}algebras is \emph{not}
triangulated.  We can represent elements of $\KK(A,B)$ by equivariant, grading
preserving $*$\nbd{}homomorphisms $\chi A\to \Comp\otimes B$ as
in~\cite{Meyer:KKG}.  However, $\chi A$ is no longer $\KK$\nbd{}equivalent
to~$A$.

The mapping cone triangle has the weak functoriality property that for any
commutative diagram
$$
\xymatrix{
  {\Sigma B} \ar[r] \ar[d]^{\Sigma\beta} &
  C \ar[r] &
  A \ar[r] \ar[d]^{\alpha} &
  B \ar[d]^{\beta} \\
  {\Sigma B'} \ar[r] &
  {C'} \ar[r] &
  {A'} \ar[r] &
  {B'}
}
$$
whose rows are exact triangles there is a morphism $\gamma\colon C\to
C'$ making the diagram commute.  The triple $(\gamma,\alpha,\beta)$ is called
a \emph{morphism of triangles}.  We do not have a functor essentially
because~$\gamma$ is not unique.  At least, the axioms of a triangulated
category guarantee that~$\gamma$ is an isomorphism if $\alpha$ and~$\beta$
are.  Thus the mapping cone and the mapping cone triangle are unique up to a
non-canonical isomorphism.

The following facts are proven in~\cite{Neeman:Triangulated}.

\begin{lemma}  \label{lem:triangle_lemma}
  Let $\Sigma B\to C\to A\overset{f}\to B$ be an exact triangle.  Then
  $C=0$ if and only if~$f$ is an isomorphism.  That is, a morphism~$f$ is an
  isomorphism if and only if its mapping cone vanishes.

  If the map $\Sigma B\to C$ vanishes, then there is an isomorphism $A\cong
  C\oplus B$ such that the maps $C\to A\to B$ become the obvious ones.  That
  is, the triangle is isomorphic to a ``direct sum'' triangle.  Conversely,
  direct sum triangles are exact.
\end{lemma}

\subsection{Long exact sequences}
\label{sec:homological}

Let~$\Tri$ be a triangulated category, for instance, $\KK^{G\cross X}$, let
$\Ab$ be the category of Abelian groups (or any Abelian category).  We call a
covariant functor $F\colon \Tri\to\Ab$ \emph{homological} if $F(C)\to F(A)\to
F(B)$ is exact for any exact triangle $\Sigma B\to C\to A\to B$.  We
define $F_n(A)\defeq F(\Sigma^n A)$ for $n\in\Z$.  Similarly, we call a
contravariant functor $F\colon \Tri\to\Ab$ \emph{cohomological} if $F(B)\to
F(A)\to F(C)$ is exact for any exact triangle, and we define
$F^n(A)\defeq F(\Sigma^n A)$.  The functor $A\mapsto \Tri(A,B)$ is
cohomological for any fixed~$B$ and the functor $B\mapsto \Tri(A,B)$ is
homological for any fixed~$A$.  This follows from the axioms of a triangulated
category.  Since we can rotate exact triangles, we obtain a long exact
sequence (infinite in both directions)
$$
\dotso\to F_n(C) \to F_n(A)\to F_n(B) \to F_{n-1}(C) \to F_{n-1}(A) \to
F_{n-1}(B) \to\dotso
$$
if~$F$ is homological, and a dual long exact sequence for
cohomological~$F$.  The maps in this sequence are induced by the maps of the
exact triangle, of course.

\subsection{Extension triangles}
\label{sec:extension_triangles}

Since any exact triangle in $\KK^{G\cross X}$ is isomorphic to a
mapping cone triangle, Section~\ref{sec:homological} only yields long exact
sequences for mapping cone triangles.  As in~\cite{Cuntz-Skandalis:Puppe},
this suffices to get long exact sequences for suitable extensions.  Let
$K\overset{i}\into E\overset{p}\prto Q$ be an extension of $G\cross
X$-$C^*$\brd{}algebras.  There is a canonical equivariant
$*$\nbd{}homomorphism $K\to\cone(p)$ that makes the diagram
\begin{equation}  \label{eq:extension_cone}
\begin{gathered}
  \xymatrix{
    \Sigma Q \ar@{=}[d] \ar@{.>}[r] &
    K \ar[r]^-{i} \ar[d] &
    E \ar[r]^-{p} \ar@{=}[d] &
    Q \ar@{=}[d] \\
    \Sigma Q \ar[r]^-{\iota} &
    \cone(p) \ar[r]^-{\epsilon} &
    E \ar[r]^-{p} &
    Q
  }
\end{gathered}
\end{equation}
commute.  The bottom row is the mapping cone triangle, of course.  In the
non-equivariant case, there is a canonical isomorphism $\KK(\Sigma Q,K)\cong
\Ext(Q,K)$.  There also exist similar results in the equivariant case
(\cite{Thomsen:Equivariant_Ext}).  If the extension has a completely positive,
contractive, equivariant cross section, then it defines an element of
$\Ext^{G\cross X}(Q,K)\cong\KK^{G\cross X}(\Sigma Q,K)$.  This provides the
dotted arrow in~\eqref{eq:extension_cone}.  Furthermore, the vertical map
$K\to\cone(p)$ is invertible in $\KK^{G\cross X}$ in this case.  This can be
proven directly and then used to prove excision for the given extension.
Conversely, it follows from excision and the Puppe sequence using the Five
Lemma.

\begin{definition}  \label{def:admissible_ext}
  We call the extension \emph{admissible} if the map $K\to\cone(p)$
  in~\eqref{eq:extension_cone} is invertible in $\KK^{G\cross X}$.  Then there
  is a unique map $\Sigma Q\to K$ that makes~\eqref{eq:extension_cone}
  commute.  The triangle $\Sigma Q\to K\to E\to Q$ is called an
  \emph{extension triangle}.
\end{definition}

If an extension is admissible, then the vertical maps
in~\eqref{eq:extension_cone} form an isomorphism of triangles.  Hence
extension triangles are exact.  Not every extension is admissible.  As we
remarked above, extensions with an equivariant, contractive, completely
positive section are admissible.  If we replace $\KK^{G\cross X}$ by
$E^{G\cross X}$, then every extension becomes admissible.

Let $f\colon A\to B$ be an equivariant $*$\nbd{}homomorphism.  We claim that
the mapping cone triangle for~$f$ is the extension triangle for an appropriate
extension.  For this we need the \emph{mapping cylinder}
\begin{equation}  \label{eq:def_cyl}
  \cyl(f)\defeq
  \{(a,b)\in A\times C([0,1],B) \mid f(a)=b(1)\}.
\end{equation}
Given $b\in B$, let $\const b\in C([0,1],B)$ be the constant function with
value~$b$.  Define natural $*$\nbd{}homomorphisms
\begin{alignat*}{2}
  p_A&\colon \cyl(f)\to A, &\qquad (a,b)&\mapsto a,
  \\
  j_A&\colon A\to\cyl(f),  &\qquad a&\mapsto (a,\const f(a)),
  \\
  \tilde{f}&\colon \cyl(f)\to B, &\qquad (a,b)&\mapsto b(0).
\end{alignat*}
Then $p_Aj_A=\ID_A$, $\tilde{f} j_A=f$, and $j_Ap_A$ is homotopic to the
identity map in a natural way.  Thus $\cyl(f)$ is homotopy equivalent to~$A$
and this homotopy equivalence identifies the maps $\tilde{f}$ and~$f$.  We
have a natural $C^*$\nbd{}extension
\begin{equation}  \label{eq:cylinder_extension}
  0 \to \cone(f)\to \cyl(f) \overset{\tilde{f}}\to B \to 0.  
\end{equation}
Build the diagram~\eqref{eq:extension_cone} for this extension.  One checks
easily that the resulting map $\cone(f)\to\cone(\tilde{f})$ is a homotopy
equivalence, so that the extension~\eqref{eq:cylinder_extension} is
admissible.  The composition $\Sigma B\to \cone(\tilde{f})
\overset{\sim}\leftarrow \cone(f)$ is naturally homotopy equivalent to the
inclusion map $\Sigma B\to\cone(f)$.  Thus the extension triangle for the
admissible extension~\eqref{eq:cylinder_extension} is isomorphic to the
mapping cone triangle for~$f$.  It follows that any exact triangle in
$\KK^{G\cross X}$ is isomorphic to an extension triangle for some admissible
extension.

\subsection{Homotopy limits}
\label{sec:homotopy_limits}

Let $(A_n,\alpha_m^n)$ be a countable inductive system in $\KK^{G\cross X}$,
with structure maps $\alpha_m^n\colon A_m\to A_n$ for $m\le n$.  (Of course,
it suffices to give the maps $\alpha_m^{m+1}$.)  Roughly speaking, its
homotopy direct limit is the correct substitute for the inductive limit for
homological computations.  Homotopy direct limits play an important role in
the proof of the Brown Representability Theorem~\ref{the:Brown_rep}.  They
also occur in connection with the behaviour of the Baum-Connes conjecture for
unions of open subgroups in Section~\ref{sec:direct_unions}.

The \emph{homotopy direct limit} $\hoinjlim A_m$ is defined to fit in an
exact triangle
\begin{equation}  \label{eq:def_hoinjlim}
  \Sigma \hoinjlim A_m \to
  \bigoplus A_m \overset{\ID-S}\longrightarrow
  \bigoplus A_m \to \hoinjlim A_m.
\end{equation}
Here~$S$ is the shift map that maps the summand~$A_m$ to the summand~$A_{m+1}$
via $\alpha_m^{m+1}$.  Thus the homotopy direct limit is
$\Sigma^{-1}\cone(\ID-S)$; it is well-defined up to non-canonical isomorphism
and has the same weak kind of functoriality as mapping cones.  The
(de)suspensions are due to the passage to opposite categories that is implicit
in our conventions.  This also means that homotopy direct limits in $\KK$
behave like homotopy inverse limits of spaces.  The map $\bigoplus A_m\to
\hoinjlim A_m$ in~\eqref{eq:def_hoinjlim} is equivalent to maps
$\alpha_m^\infty\colon A_m\to \hoinjlim A_m$ with
$\alpha_n^\infty\circ\alpha_m^n=\alpha_m^\infty$ for $m\le n$.

To formulate the characteristic properties of the homotopy limit, we consider
(co)homological functors to the category of Abelian groups that are
\emph{compatible with direct sums}.  The latter means $F(\bigoplus A_m) \cong
\bigoplus F(A_m)$ in the covariant case and $F(\bigoplus A_m) \cong \prod
F(A_m)$ in the contravariant case.  The functor $B\mapsto \Tri(A,B)$ is not
always compatible with direct sums.  We call~$A$ \emph{compact} if it is.  The
functor $A\mapsto \Tri(A,B)$ is always compatible with direct sums: this is
just the universal property of direct sums.  Hence the following lemma applies
to $F(A)\defeq \Tri(A,B)$ for any~$B$.

\begin{lemma}[\cite{Neeman:Thomason}]  \label{lem:holim}
  If~$F$ is homological and compatible with direct sums, then the maps
  $\alpha_m^\infty\colon A_m\to\hoinjlim A_m$ yield an isomorphism $\varinjlim
  F_n(A_m) \congto F_n(\hoinjlim A_m)$.  If~$F$ is cohomological and
  compatible with direct sums, then there is a short exact sequence
  \begin{displaymath}
    0 \to \varprojlim\nolimits^1 F^{n-1}(A_m)
    \to F^n(\hoinjlim A_m)
    \to \varprojlim F^n(A_m) \to 0.
  \end{displaymath}
  The map $F^n(\hoinjlim A_m) \to \varprojlim F^n(A_m)$ is induced by
  $(\alpha_m^\infty)_{m\in\N}$.
\end{lemma}

\begin{proof}
  Consider the homological case first.  Apply the long exact homology sequence
  to~\eqref{eq:def_hoinjlim} and cut the result into short exact sequences of
  the form
  \begin{multline*}
    \coker \Bigl(\ID-S\colon \bigoplus F_n(A_m)\to \bigoplus F_n(A_m)\Bigr)
    \into F_n(\hoinjlim A_m)
    \\ \prto \ker \Bigl(\ID-S\colon \bigoplus F_{n-1}(A_m)
    \to \bigoplus F_{n-1}(A_m)\Bigr).
  \end{multline*}
  The kernel of $\ID-S$ vanishes and its cokernel is, by definition,
  $\varinjlim F_n(A_m)$.  Whence the assertion.  In the cohomological case, we
  get a short exact sequence
  \begin{multline*}
    \coker \Bigl(\ID-S\colon \prod F^{n-1}(A_m)\to \prod F^{n-1}(A_m)\Bigr)
    \into F^n(\hoinjlim A_m)
    \\ \prto \ker \Bigl(\ID-S\colon \prod F^n(A_m)\to \prod F^n(A_m)\Bigr).
  \end{multline*}
  By definition, the kernel is $\varprojlim F^n(A_m)$ and the cokernel is
  $\varprojlim\nolimits^1 F^{n-1}(A_m)$.
\end{proof}

We now specialise to the category $\KK^{G\cross X}$ and relate homotopy direct
limits to ordinary direct limits via mapping telescopes.  This is used in our
discussion of unions of groups in Section~\ref{sec:direct_unions}.  Any
inductive system in $\KK^{G\cross X}$ is isomorphic to the image of a direct
system of $G\cross X$-$C^*$\brd{}algebras.  That is, the maps~$\alpha_m^n$ are
equivariant $*$\nbd{}homomorphisms and satisfy
$\alpha_m^n\circ\alpha_l^m=\alpha_l^n$ as such.  To get this, replace
the~$A_m$ by the universal algebra $q_s(A_m)$ as in the appendix.

The following discussion follows the treatment of inductive limits
in~\cite{Schochet:Axiomatic}.  Let $(A_m,\alpha_m^n)$ be an inductive system
of $G\cross X$\brd{}$C^*$\brd{}algebras.  We let $A_\infty\defeq\varinjlim
A_m$ and denote the natural maps $A_m\to A_\infty$ by~$\alpha_m^\infty$.  The
\emph{mapping telescope} of the system is defined as the $G$-$C^*$-algebra
\begin{multline*}
  T(A_m,\alpha_m^n) \defeq
  \biggl\{ (f_m)\in \bigoplus_{m\in\N} C([m,m+1],A_m) \biggm|
  \\
  \text{$f_0(0)=0$ and $f_{m+1}(m+1)=\alpha_m^{m+1}\bigl(f_m(m+1)\bigr)$}
  \biggr\}.
\end{multline*}
In the special case where the homomorphisms $\alpha_m^n$ are injective,
$T(A_m,\alpha_m^n)$ is the space of all $f\in
C_0(\mathopen]0,\infty\mathclose[,A_\infty)$ with $f(t)\in A_m$ for $t\le
m+1$.  In particular, for the constant inductive system $(A_\infty,\ID)$ we
obtain just the suspension $\Sigma A_\infty$.  Since the mapping telescope
construction is functorial, there is a natural equivariant
$*$\nbd{}homomorphism $T(A_n,\alpha_m^n)\to \Sigma A_\infty$.

\begin{definition}  \label{def:admissible_ind}
  An inductive system $(A_m,\alpha_m^n)$ is called \emph{admissible} if the
  map $T(A_n,\alpha_m^n)\to \Sigma A_\infty$ is invertible in $\KK^{G\cross
  X}$.
\end{definition}

\begin{proposition}  \label{pro:holim}
  We have $\varinjlim {} (A_m,\alpha_m^n) \cong \hoinjlim {}(A_m,\alpha_m^n)$
  for admissible inductive systems.
\end{proposition}

\begin{proof}
  Evaluation at positive integers defines a natural, surjective, equivariant
  $*$\nbd{}homomorphism $\pi\colon T(A_m,\alpha_m^n)\to \bigoplus A_m$.  Its
  kernel is naturally isomorphic to $\bigoplus \Sigma A_m$.  Thus we obtain a
  natural extension
  $$
  0 \longrightarrow \bigoplus \Sigma A_m
  \overset{\iota}\longrightarrow T(A_m,\alpha_m^n)
  \overset{\pi}\longrightarrow \bigoplus A_m
  \longrightarrow 0.
  $$
  Build the diagram~\eqref{eq:extension_cone} for this extension.  The map
  $\bigoplus \Sigma A_n\to \cone(\pi)$ is a homotopy equivalence in a natural
  and hence equivariant fashion.  Hence the extension is admissible.
  Moreover, one easily identifies the map $\Sigma\bigl(\bigoplus
  A_m\bigr)\to\bigoplus \Sigma A_m$ with $S-\ID$, where~$S$ is the shift map
  defined above.  Rotating the extension triangle, we obtain an exact triangle
  $$
  \xymatrix@1@C=3.5em{
    T(A_m,\alpha_m^n) \ar[r]^-{-\pi} &
    \bigoplus A_m \ar[r]^-{\ID-S} &
    \bigoplus A_m \ar[r]^-{\Sigma^{-1}\iota} &
    \Sigma^{-1} T(A_m,\alpha_m^n).
  }
  $$
  This implies $\Sigma^{-1} T(A_m,\alpha_m^n) \cong \hoinjlim
  (A_m,\alpha_m^n)$ and hence the assertion.
\end{proof}

To obtain a concrete criterion for admissibility, we let
$\tilde{T}(A_m,\alpha_m^n)$ be the variant of $T(A_m,\alpha_m^n)$ where we
require $\lim_{t\to\infty} \alpha_m^\infty(f_m(t))$ to exist in~$A_\infty$
instead of $\lim f_m(t)=0$.  The algebra $\tilde{T}(A_m,\alpha_m^n)$ is
equivariantly contractible in a natural way.  The contracting homotopy is
obtained by making sense of the formula $H_s f(t)\defeq \alpha_{[st]}^{[t]}
f(st)$ for $0\le s\le 1$.  There is a natural commutative diagram
$$
\xymatrix{
  0 \ar[r] &
  T(A_m,\alpha_m^n) \ar[r]^{\subseteq} \ar[d] &
  \tilde{T}(A_m,\alpha_m^n) \ar[r]^-{\ev_\infty} \ar[d] &
  A_\infty \ar[r] \ar@{=}[d] &
  0 \\
  0 \ar[r] &
  T(A_\infty,\ID) \ar[r]^{\subseteq} &
  \tilde{T}(A_\infty,\ID) \ar[r]^-{\ev_\infty} &
  A_\infty \ar[r] &
  0
}
$$
whose rows are short exact sequences.  The bottom extension is evidently
admissible.  By definition, the vertical map on $T(\dotso)$ is invertible in
$\KK^{G\cross X}$ if and only if the inductive system is admissible.  The
other vertical maps are invertible in any case because
$\tilde{T}(\dotso)\cong0$ in $\KK^{G\cross X}$.  Therefore, if the inductive
system is admissible, then the top row is an admissible extension whose
extension triangle is isomorphic to the one for the bottom row.  Conversely,
if the top row is an admissible extension, then the vertical map on
$T(\dotso)$ is invertible in $\KK^{G\cross X}$ by the uniqueness of mapping
cones.  As a result, the inductive system is admissible if and only if the
extension in the top row above is admissible.  In $E^{G\cross X}$, all
inductive systems are admissible because all extensions are admissible.

\begin{lemma}  \label{lem:admissible_ind}
  An inductive system $(A_m,\alpha_m^n)$ is admissible if there exist
  equivariant completely positive contractions $\phi_m\colon A_\infty\to A_m$
  such that $\alpha_m^\infty\circ\phi_m\colon A_\infty\to A_\infty$ converges
  in the point norm topology towards the identity.
\end{lemma}

\begin{proof}
  By the above discussion, the inductive system is admissible if there is an
  equivariant, contractive, completely positive cross section $A_\infty\to
  \tilde{T}(A_m,\alpha_m^n)$.  It is not hard to see that such a cross section
  exists if and only if there are maps~$\phi_m$ as in the statement of the
  lemma.
\end{proof}

\subsection{Triangulated functors and subcategories}
\label{sec:functors_subcategories}

A \emph{triangulated subcategory} of a triangulated category~$\Tri$ is a full
subcategory $\Tri'\subseteq\Tri$ that is closed under suspensions and has the
exactness property that if $\Sigma B\to C\to A\to B$ is an exact triangle with
$A,B\in\Tri'$, then $C\in\Tri'$ as well.  In particular, $\Tri'$ is closed
under isomorphisms and finite direct sums.  A triangulated subcategory is
called \emph{thick} if all retracts (direct summands) of objects of~$\Tri'$
belong to~$\Tri'$.  A triangulated subcategory is indeed a triangulated
category in its own right.  Given any class of objects~$\GEN$, there is a
smallest (thick) triangulated subcategory containing~$\GEN$.  This is called
the \emph{(thick) triangulated subcategory generated by~$\GEN$}.  Since a full
subcategory is determined by its class of objects, we do not distinguish
between full subcategories and classes of objects.

Let~$\aleph$ be some infinite regular cardinal number.  In our applications we
only use the countable cardinal number~$\aleph_0$.  We suppose that direct
sums of cardinality~$\aleph$ exist in~$\Tri$.  A subcategory of~$\Tri$ is
called \emph{($\aleph$\nbd{})localising} if it is triangulated and closed
under direct sums of cardinality~$\aleph$.  We can define the \emph{localising
  subcategory generated by some class~$\GEN$ of objects} as above.  We denote
it by $\gen{\GEN}$ or $\gen{\GEN}^\aleph$.  Notice that a triangulated
subcategory that is closed under direct sums is also closed under homotopy
direct limits.  Localising subcategories are automatically thick
(see~\cite{Neeman:Triangulated}).

It is easy to see that an $\aleph_0$\brd{}localising subcategory of
$\KK^{G\cross X}$ amounts to a class~$\Null$ of $G\cross
X$-$C^*$\brd{}algebras with the following properties:
\begin{enumerate}[(1)]
\item if $A$ and~$B$ are $\KK^{G\cross X}$\brd{}equivalent and $A\in\Null$,
  then $B\in\Null$;

\item $\Null$ is closed under suspension;
  
\item if $f\colon A\to B$ is an equivariant $*$\nbd{}homomorphism with
  $A,B\in\Null$, then also $\cone(f)\in\Null$;
  
\item if $A_n\in\Null$ for all $n\in\N$, then also $\bigoplus_{n\in\N}
  A_n\in\Null$.

\end{enumerate}
We can replace (3) and~(4) by the equivalent conditions
\begin{enumerate}[(1')]
\setcounter{enumi}{2}
\item if $K\into E\prto Q$ is an admissible extension and two of $K$, $E$, $Q$
  belong to~$\Null$, so does the third;
  
\item if $(A_n,\alpha_m^n)$ is an admissible inductive system with
  $A_n\in\Null$ for all $n\in\N$, then $\varinjlim A_n\in\Null$ as well.

\end{enumerate}

Thus the localising subcategory generated by a class~$\GEN$ of $G\cross
X$\nbd{}$C^*$\brd{}algebras is the smallest class of $G\cross
X$\nbd{}$C^*$\brd{}algebras containing~$\GEN$ with the above four properties.
For example, the localising subcategory of $\KK$ generated by~$\pt$ is exactly
the \emph{bootstrap category} (see~\cite{Blackadar:Book}).  The proof uses
that extensions and inductive systems of nuclear $C^*$\nbd{}algebras are
automatically admissible.  Another example of a localising subcategory of
$\KK$ is the class of $C^*$\nbd{}algebras with vanishing $\K$\nbd{}theory.  We
discuss these two subcategories further in Section~\ref{sec:representability}
to give an easy application of the Brown Representability Theorem.

Let $\Tri$ and~$\Tri'$ be triangulated categories.  A functor $F\colon
\Tri\to\Tri'$ is called \emph{triangulated} if it is additive, intertwines the
translation automorphisms, and maps exact triangles to exact triangles.
Although the latter condition may look like an exactness condition, it is
almost empty.  Since any exact triangle in $\KK^{G\cross X}$ is isomorphic to
a mapping cone triangle, a functor is triangulated once it commutes with
suspensions and preserves mapping cone triangles.  For instance, the functor
$A\mapsto A\otimes_{\mathrm{min}} B$ has this property regardless of
whether~$B$ is exact.  Similarly, the full and reduced crossed product
functors $\KK^{G\cross X}\to\KK$ are triangulated.  An analogous situation
occurs in homological algebra: any additive functor between Abelian categories
gives rise to a triangulated functor between the homotopy categories of chain
complexes.  The exactness of the functor only becomes relevant for the derived
category.

Let $F\colon \Tri\to\Tri'$ be a triangulated functor.  Its \emph{kernel} is
the class $\ker F$ of all objects~$X$ of~$\Tri$ with $F(X)\cong 0$.  It is
easy to see that $\ker F$ is a thick triangulated subcategory of~$\Tri$.
If~$F$ commutes with direct sums of cardinality~$\aleph$, then $\ker F$ is
$\aleph$\nbd{}localising.

\subsection{Localisation of categories and functors}
\label{sec:localisation}

A basic (and not quite correct) result on triangulated categories asserts that
any thick triangulated subcategory $\Null\subseteq\Tri$ arises as the kernel
of a triangulated functor.  Even more, there exists a universal triangulated
functor $\Tri\to\Tri/\Null$ with kernel~$\Null$, called \emph{localisation
  functor}, such that any other functor whose kernel contains~$\Null$
factorises uniquely through $\Tri/\Null$ (see~\cite{Neeman:Triangulated}).
Its construction is quite involved and may fail to work in general because the
morphism spaces in $\Tri/\Null$ may turn out to be classes and not sets.

There are two basic examples of localisations, which have motivated the whole
theory of triangulated categories.  They come from homological algebra and
homotopy theory, respectively.  In homological algebra, the ambient
category~$\Tri$ is the homotopy category of chain complexes over an Abelian
category.  The subcategory $\Null\subseteq\Tri$ consists of the exact
complexes, that is, complexes with vanishing homology.  A chain map is called
a \emph{quasi\brd{}isomorphism} if it induces an isomorphism on homology.  The
localisation $\Tri/\Null$ is, by definition, the derived category of the
underlying Abelian category.  One of the motivations for developing the theory
of triangulated categories was to understand what additional structure of the
homotopy category of chain complexes is inherited by the derived category.

In homotopy theory there are several important instances of localisations.  We
only discuss one very elementary situation which provides a good analogy for
our treatment of the Baum-Connes assembly map.  Let~$\Tri$ be the stable
homotopy category of all topological spaces.  We call an object of~$\Tri$
\emph{weakly contractible} if its stable homotopy groups vanish.  A map is
called a \emph{weak homotopy equivalence} if it induces an isomorphism on
stable homotopy groups.  Let $\Null\subseteq\Tri$ be the subcategory of weakly
contractible objects.  In homotopy theory one often wants to disregard objects
of~$\Null$, that is, work in the localisation $\Tri/\Null$.

The concepts of a weak equivalence in homotopy theory and of a
quasi\brd{}isomorphism in homological algebra become equivalent once
formulated in terms of triangulated categories: we call a morphism
$f\in\Tri(A,B)$ an \emph{$\Null$\nbd{}weak equivalence} or an
\emph{$\Null$\nbd{}quasi\brd{}isomorphisms} if $\cone(f)\in\Null$.  Since
$N\in\Null$ if and only if $0\to N$ is an $\Null$\nbd{}weak equivalence, the
weak equivalences and~$\Null$ determine each other uniquely.

A morphism is a weak equivalence if and only if its image in the localisation
$\Tri/\Null$ is an isomorphism.  This implies several cancellation assertions
about weak equivalences.  For instance, if $f$ and~$g$ are composable and two
of the three morphisms $f$, $g$ and $f\circ g$ are weak equivalences, so is
the third.  The localisation has the universal property that any functor out
of~$\Tri$, triangulated or not, that maps $\Null$\nbd{}weak equivalences to
isomorphisms, factorises uniquely through the localisation.

In many examples of localisation, some more structure is present.  This is
formalised in the following definition.  In the simplicial approximation
example, let $\Proj\subseteq\Tri$ be the subcategory of all objects that have
the stable homotopy type of a CW\brd{}complex, that is, are isomorphic
in~$\Tri$ to a CW\brd{}complex.  The Whitehead Lemma implies that $f\colon
X\to Y$ is a weak homotopy equivalence if and only if $f_*\colon
\Tri(P,X)\to\Tri(P,Y)$ is an isomorphism for all $P\in\Proj$.  Similarly,
$N\in\Null$ if and only $\Tri(P,N)=0$ for all $P\in\Proj$.  Another important
fact is that any space~$S$ has a \emph{simplicial approximation}.  This is
just a weak equivalence $\tilde{X}\to X$ with $\tilde{X}\in\Proj$.  In
homological algebra, a similar situation arises if there are ``enough
projectives''.  Then one lets~$\Proj$ be the subcategory of projective chain
complexes (see~\cite{Keller:Handbook}).

\begin{definition}  \label{def:complementary}
  Let~$\Tri$ be a triangulated category and let $\Proj$ and~$\Null$ be
  thick triangulated subcategories of~$\Tri$.  We call the pair
  $(\Proj, \Null)$ \emph{complementary} if $\Tri(P,N)=0$ for all
  $P\in\Proj$, $N\in\Null$ and if for any $A\in\Tri$ there is an
  exact triangle $\Sigma N\to P\to A\to N$ with $P\in\Proj$,
  $N\in\Null$.
\end{definition}

We shall only need localisations in the situation of complementary
subcategories.  In this case, the construction of $\Tri/\Null$ is easier and
there is some important (and well-known) additional structure
(see~\cite{Neeman:Thomason}).  We prove some basic results because they are
important for our treatment of the Baum-Connes assembly map.

\begin{proposition}  \label{pro:complementary}
  Let~$\Tri$ be a triangulated category and let $(\Proj,\Null)$ be
  complementary thick triangulated subcategories of~$\Tri$.
  \begin{enumerate}[\ref{pro:complementary}.1.]
  \item \label{pro::PN_orthogonal} We have $N\in\Null$ if and only if
    $\Tri(P,N)=0$ for all $P\in\Proj$, and $P\in\Proj$ if and only if
    $\Tri(P,N)=0$ for all $N\in\Null$; thus $\Proj$ and~$\Null$ determine each
    other.
    
  \item \label{pro::PN_functors} The exact triangle $\Sigma N\to P\to
    A\to N$ with $P\in\Proj$ and $N\in\Null$ is uniquely determined up to
    isomorphism and depends functorially on~$A$.  In particular, its entries
    define functors $P\colon \Tri\to\Proj$ and $N\colon \Tri\to\Null$.
    
  \item \label{pro::PN_functors_triangulated} The functors
    $P,N\colon\Tri\to\Tri$ are triangulated.
    
  \item \label{pro::localisexists} The localisations $\Tri/\Null$ and
    $\Tri/\Proj$ exist.
    
  \item \label{pro::localise_sub} The compositions $\Proj\to\Tri\to\Tri/\Null$
    and $\Null\to\Tri\to\Tri/\Proj$ are equivalences of triangulated
    categories.
    
  \item The functors \label{pro::PN_descend} $P,N\colon\Tri\to\Tri$ descend to
    triangulated functors $P\colon \Tri/\Null\to\Proj$ and $N\colon
    \Tri/\Proj\to\Null$, respectively, that are inverse (up to isomorphism) to
    the functors in \ref{pro:complementary}.\ref{pro::localise_sub}.
    
  \item \label{pro::PN_adjoint} The functors $P\colon \Tri/\Null\to\Tri$ and
    $N\colon \Tri/\Proj\to\Tri$ are left and right adjoint to the localisation
    functors $\Tri\to\Tri/\Null$ and $\Tri\to\Tri/\Proj$, respectively; that
    is, we have natural isomorphisms
    $$
    \Tri(P(A),B) \cong \Tri/\Null(A,B), \qquad
    \Tri(A,N(B)) \cong \Tri/\Proj(A,B),
    $$
    for all $A,B\in\Tri$.

  \end{enumerate}
\end{proposition}

\begin{proof}
  We can exchange the roles of $\Proj$ and~$\Null$ by passing to opposite
  categories.  Hence it suffices to prove the various assertions about one of
  them.
  
  By hypothesis, $N\in\Null$ implies $\Tri(P,N)=0$ for all $P\in\Proj$.
  Conversely, suppose $\Tri(P,A)=0$ for all $P\in\Proj$.  Let $\Sigma N\to
  P\to A\to N$ be an exact triangle with $P\in\Proj$ and $N\in\Null$.
  The map $P\to A$ vanishes by hypothesis.  Lemma~\ref{lem:triangle_lemma}
  implies $N\cong A\oplus \Sigma^{-1} P$.  Since~$\Null$ is thick,
  $A\in\Null$.  This proves \ref{pro:complementary}.\ref{pro::PN_orthogonal}.
  
  Let $\Sigma N\to P\to A\to N$ and $\Sigma N'\to P'\to A'\to N'$ be
  exact triangles with $P,P'\in\Proj$ and $N,N'\in\Null$ and let
  $f\in\Tri(A,A')$.  Since $\Tri(P,N')=0$, the map $P'\to A'$ induces an
  isomorphism $\Tri(P,P')\cong\Tri(P,A')$.  Hence there is a unique and hence
  natural way to lift the composite map $P\to A\to A'$ to a map $P\to P'$.  By
  the axioms of a triangulated category, there exists a morphism of
  exact triangles from $\Sigma N\to P\to A\to N$ to $\Sigma N'\to
  P'\to A'\to N'$ that extends $f\colon A\to A'$ and its lifting $P(f)\colon
  P\to P'$.  An argument as above shows that there is a unique way to lift~$f$
  to a map $N\to N'$.  Thus the morphism of triangles that extends~$f$ is
  determined uniquely, so that the triangle $\Sigma N\to P\to A\to N$ depends
  functorially on~$A$.  This proves
  \ref{pro:complementary}.\ref{pro::PN_functors}.
  
  Next we show that~$P$ is a triangulated functor on~$\Tri$.  Let $\Sigma B\to
  C \to A\to B$ be an exact triangle.  Consider the solid arrows in the
  diagram in Figure~\ref{fig:PN_exact}.
  \begin{figure}
    $$
    \xymatrix{ \Sigma^2 N(B) \ar@{.>}[r] \ar[d] \ar@{}[dr]|{-}
      & \Sigma N'(C) \ar@{.>}[r] \ar@{.>}[d] & \Sigma N(A) \ar[r]
      \ar[d] &
      \Sigma N(B) \ar[d] \\
      \Sigma P(B) \ar@{.>}[r] \ar[d] & P'(C) \ar@{.>}[r] \ar@{.>}[d] &
      P(A) \ar@{.>}[r] \ar[d] &
      P(B) \ar[d] \\
      \Sigma B \ar[r] \ar[d] & C \ar[r] \ar@{.>}[d] & A \ar[r] \ar[d]
      &
      B \ar[d] \\
      \Sigma N(B) \ar@{.>}[r] & N'(C) \ar@{.>}[r] & N(A) \ar[r] &
      N(B) \\
    }
    $$
    \caption{Exactness of $P$ and~$N$}
    \label{fig:PN_exact}
  \end{figure}
  We can find objects $N'(C)$ and $P'(C)$ of~$\Tri$ and the dotted arrows
  in this diagram so that all rows and columns are exact and such that the
  diagram commutes except for the square marked with a~$-$, which
  anti-commutes (see \cite{Beilinson-Bernstein-Deligne}*{Proposition 1.1.11}).
  Since $\Proj$ and~$\Null$ are triangulated subcategories and the rows in
  this diagram are exact triangles, we get $P'(C)\in\Proj$ and
  $N'(C)\in\Null$.  Hence the column over~$C$ is as in the definition of the
  functors $P$ and~$N$.  Therefore, we can replace this column by the exact
  triangle $\Sigma N(C)\to P(C)\to C\to N(C)$.  Our proof of
  \ref{pro:complementary}.\ref{pro::PN_functors} shows that the rows must be
  obtained by applying the functors $P$ and~$N$ to the given exact triangle
  $\Sigma B\to C\to A\to B$.  Since the rows are exact triangles by
  construction, the functors $P$ and~$N$ preserve exact triangles.  They
  evidently commute with suspensions.  This proves
  \ref{pro:complementary}.\ref{pro::PN_functors_triangulated}.
  
  Next we construct a candidate~$\Tri'$ for the localisation $\Tri/\Null$.  We
  let~$\Tri'$ have the same objects as~$\Tri$ and morphisms $\Tri'(A,B)\defeq
  \Tri(P(A),P(B))$.  The identity map on objects and the map~$P$ on morphisms
  define a canonical functor $\Tri\to\Tri'$.  We define the suspension
  on~$\Tri'$ to be the same as for~$\Tri$.  A triangle in~$\Tri'$ is called
  exact if it is isomorphic to the image of an exact triangle in~$\Tri$.  We
  claim that~$\Tri'$ with this additional structure is a triangulated category
  and that the functor $\Tri\to\Tri'$ is the localisation functor at~$\Null$.
  
  The uniqueness of the exact triangle $\Sigma N(A)\to P(A)\to A\to
  N(A)$ yields that the natural map $P(A)\to A$ is an isomorphism for
  $A\in\Proj$.  Therefore, the map $\Tri(A,B)\to\Tri'(A,B)$ is an isomorphism
  for $A,B\in\Proj$.  That is, the restriction of the functor $\Tri\to\Tri'$
  to~$\Proj$ is fully faithful and identifies~$\Proj$ with a full subcategory
  of~$\Tri'$.  Moreover, since $P(A)\in\Proj$, the map $P^2(A)\to P(A)$ is an
  isomorphism.  This implies that the map $P(A)\to A$ is mapped to an
  isomorphism in~$\Tri'$.  Thus any object of~$\Tri'$ is isomorphic to one in
  the full subcategory~$\Proj$.  Therefore, the category~$\Tri'$ is equivalent
  to the subcategory~$\Proj$.  Using that~$P$ is a triangulated functor
  on~$\Tri$, one shows easily that both functors $\Proj\to\Tri'$ and
  $\Tri'\to\Proj$ map exact triangles to exact triangles.
  They commute with suspensions anyway.  Since they are equivalences of
  categories and since~$\Proj$ is a triangulated category, the
  category~$\Tri'$ is triangulated and the equivalence $\Proj\cong\Tri'$ is
  compatible with the triangulated category structure.
  
  We define the functor $P\colon \Tri'\to\Proj$ to be~$P$ on objects and the
  identity on morphisms.  This functor is clearly inverse to the above
  equivalence $\Proj\to\Tri'$ and has the property that the composition
  $\Tri\to\Tri'\overset{P}\to\Proj\subseteq\Tri$ agrees with $P\colon
  \Tri\to\Tri$.  Moreover, we have observed already above that
  $\Tri(P(A),B))\cong \Tri(P(A),P(B))$ for all $A,B\in\Tri$.  Hence all the
  remaining assertions follow once we show that~$\Tri'$ has the universal
  property of $\Tri/\Null$.  It is easy to see that~$\Null$ is equal to the
  kernel of $\Tri\to\Tri'$.  If $F\colon \Tri\to\Tri''$ is a triangulated
  functor with kernel~$\Null$, then the maps $P(A)\to A$ induce isomorphisms
  $F(P(A))\to F(A)$ by Lemma~\ref{lem:triangle_lemma}.  Therefore,
  $\Tri'\overset{P}\to\Proj\subseteq\Tri\overset{F}\to\Tri''$ is the required
  factorisation of~$F$ through~$\Tri'$.
\end{proof}

We call the map $P(A)\to A$ an \emph{$\Null$\nbd{}projective resolution} or a
\emph{$\Proj$\nbd{}simplicial approximation} of~$A$.  The first term comes
from homological algebra, the second one from homotopy theory.  We prefer the
terminology from homotopy theory because it gives a more accurate analogy for
the Baum-Connes assembly map.

Finally, we consider the localisation of functors.  Let $F\colon \Tri\to\Tri'$
be a covariant triangulated functor to another triangulated category~$\Tri'$.
Then its \emph{localisation} or \emph{left derived functor} $\Left F\colon
\Tri/\Null\to\Tri'$ is, in general, defined by a certain universal property.
In the case of a complementary pair of subcategories, it is given simply by
$\Left F\cong F\circ P$.  This makes sense for any functor~$F$, triangulated
or not.  If~$F$ is triangulated, then so is $\Left F$.  If~$F$ is
(co)homological, then so is $\Left F$.  Both assertions follow from
Proposition \ref{pro:complementary}.\ref{pro::PN_functors_triangulated}.  In
the following discussion, we assume~$F$ to be triangulated or homological.

The functor $\Left F$ descends to the category $\Tri/\Null$ and comes equipped
with a natural transformation $\Left F\to F$ which comes from the natural
transformation $P(A)\to A$.  The universal property that characterises~$\Left
F$ is the following.  If $F'\colon\Tri/\Null\to\Tri'$ is any functor together
with a natural transformation $F'\to F$, then this natural transformation
factorises uniquely through $\Left F$.  This factorisation is obtained as the
composition
$$
F'(A) \overset{\cong}\leftarrow F'(P(A))\to F(P(A)) \cong \Left F(A).
$$
Thus we may view $\Left F$ as the best approximation to~$F$ that factors
through $\Tri/\Null$.  In particular, we have $\Left F\cong F$ if and only if
$\Null\subseteq \ker F$ if and only if~$F$ maps $\Null$\nbd{}weak equivalences
to isomorphisms in~$\Tri'$.

Alternatively, we may view $\Left F(A)$ as the best approximation to $F(A)$
that uses only the restriction of~$F$ to~$\Proj$.  The simplicial
approximation $P(A)\to A$ has the universal property that any map $B\to A$
with $B\in\Proj$ factors uniquely through $P(A)$.  In this sense, $P(A)$ is
the best possible approximation to~$A$ inside~$\Proj$ and $F(P(A))$ is the
best guess we can make for $F(A)$ if we want the guess to be of the form
$F(B)$ for some $B\in\Proj$.

We can also use the functor $N\colon \Tri/\Proj\to\Tri$ to define an
\emph{obstruction functor} $\Obs F\defeq F\circ N$.  It comes equipped with a
natural transformation $F\to\Obs F$.
Proposition~\ref{pro:complementary}.\ref{pro::PN_functors} shows that if the
functor~$F$ is triangulated then $\Left F$, $F$ and $\Obs F$ are related by a
natural exact triangle
$$
\Sigma \Obs F(A)\to \Left F(A) \to F(A)\to \Obs F(A).
$$
Thus $\Obs F(A)$ measures the lack of invertibility of the map $\Left
F(A)\to F(A)$.  In particular, $\Obs F(A)=0$ if and only if $\Left F(A)\cong
F(A)$.  Similar remarks apply if~$F$ is homological.  In that case, the
functors $\Left F$, $F$ and $\Obs F$ are related by a long exact sequence.

\section{Preliminaries on compact subgroups and some functors}
\label{sec:preliminaries}

We first recall some structural results about compact subgroups in locally
compact groups.  Then we recall the well-known formal properties of tensor
product, restriction and induction functors.  We discuss them in some detail
because they are frequently used.  We apply the universal property of
$\KK$-theory to treat them.  This has the advantage that proofs do not require
the definition of $\KK$.

\subsection{Compact subgroups}
\label{sec:compact_subgroups}

Let~$G$ be a locally compact group.  Let $G_0\subseteq G$ be the connected
component of the identity element.  We call~$G$ \emph{almost totally
disconnected} if~$G_0$ is compact, and \emph{almost connected} if $G/G_0$ is
compact.  If~$G$ is almost totally disconnected, then~$G$ contains a compact
open subgroup (and vice versa) by~\cite{Hewitt-Ross}*{Theorem 7.5}.
Therefore, if~$G$ is arbitrary, then there exists an open almost connected
subgroup $U\subseteq G$: take the preimage of a compact open subgroup in
$G/G_0$.  Almost connected groups are very closely related to Lie groups (with
finitely many connected components) by~\cite{Montgomery-Zippin}: if~$U$ is
almost connected, then each neighbourhood of the identity element contains a
compact normal subgroup $N\subseteq U$ such that $U/N$ is a Lie group (the
smooth structure on $U/N$ is unique if it exists).

Let~$U$ be almost connected and let $K\subseteq U$ be maximal compact.  We
recall some structural results about $U/K$ from~\cite{Abels:Slices}.  Let
$\mathfrak{k}\subseteq\mathfrak{u}$ be the Lie algebras of $K$ and~$U$,
respectively, and let $\mathfrak{p} \defeq \mathfrak{u}/\mathfrak{k}$.  This
quotient is a finite dimensional $\R$\nbd{}vector space, on which~$K$ acts
linearly by conjugation.  There exists a $K$\nbd{}equivariant homeomorphism
$U/K\cong\mathfrak{p}$.  Thus $U/K$ as a $K$-space is homeomorphic to a linear
action of~$K$ on a real vector space.  This fact is crucial for our purposes.
Moreover, Abels shows in~\cite{Abels:Slices} that $U/K$ is a universal proper
$U$\nbd{}space.  This contains the assertion that any compact subgroup of~$U$
is subconjugate to~$K$ (because it fixes a point in $U/K$).  Especially, any
two maximal compact subgroups are conjugate.

We define some classes of special compact subgroups that we shall use later.
Let $H\subseteq G$ be a compact subgroup.  We call~$H$ \emph{strongly smooth}
if its normaliser $N_GH\subseteq G$ is open in~$G$ and $N_GH/H$ is a Lie
group.  We call~$H$ \emph{smooth} if it contains a strongly smooth subgroup
of~$G$.  Finally, we call~$H$ \emph{large} if it is a maximal compact subgroup
of some open almost connected subgroup of~$G$.  We let $\MCAC=\MCAC(G)$ be the
set of large compact subgroups of~$G$.

Of course, strongly smooth subgroups are smooth.  Large subgroups are also
smooth because if $L\subseteq U$ is maximal compact and $N\subseteq U$ is a
smooth, compact normal subgroup, then $NL$ is a compact subgroup as well by
normality.  Hence $N\subseteq L$ by maximality.

\begin{lemma}  \label{lem:MCAC_enough}
  Any compact subgroup of~$G$ is contained in a large compact subgroup.
  
  If $H\subseteq G$ is a large compact subgroup, then the open almost
  connected subgroup $U\subseteq G$ in which~$H$ is maximal is unique and
  denoted by~$U_H$.
  
  Suppose $H,L\in\MCAC$ satisfy $H\subseteq L$.  Then $H=U_H\cap L$, so
  that~$H$ is open in~$L$.  The natural map $U_H/H\to U_L/L$ is a
  homeomorphism.
  
  If $H\subseteq G$ is smooth, then the homogeneous space $G/H$ is a smooth
  manifold in a canonical way.
\end{lemma}

\begin{proof}
  We claim that any compact subgroup~$H$ of a totally disconnected group~$G$
  is contained in a compact open one.  Let $U\subseteq G$ be any compact open
  subgroup.  Then $H\cap U$ has finite index in~$H$.  Therefore, $U'\defeq
  \bigcap_{h\in H} hUh^{-1}$ is again a compact open subgroup of~$G$.  By
  construction, it is normalised by~$H$, so that $HU'$ is again a subgroup.
  It is compact and open and contains~$H$.  Since $HU'$ is almost connected,
  $H$ is contained in some maximal compact subgroup of $HU'$.  This yields the
  first assertion.
  
  Suppose~$H$ is maximal compact in the open almost connected subgroups $U$
  and~$V$ of~$G$.  We claim that $U=V$.  Since~$H$ is still maximal compact in
  $U\cap V$, we may assume that $U\subseteq V$.  Hence $V/H$ is a disjoint
  union of $U:V$ copies of $U/H$.  However, $V/H$ is homeomorphic to a vector
  space and therefore connected, forcing $U=V$.  Let $H$ and~$L$ be large
  compact subgroups of~$G$ that satisfy $H\subseteq L$.  Then $H\subseteq
  U_H\cap L\subseteq U_H$, so that $H=U_H\cap L$ by maximality.  Hence the
  natural map $U_H/H\to U_L/L$ is injective.  Its image is both open and
  closed and hence must be all of $U_L/L$ by connectedness.
  
  Let $H\subseteq G$ be smooth.  Then we can find an almost connected open
  subgroup $U\subseteq G$ that contains~$H$ and a subgroup $N\subseteq H$ that
  is normal in~$U$ such that $U/N$ is a Lie group.  Write $G/N$ as a disjoint
  union of copies of the Lie group $U/N$.  This reveals that $G/H$ is a
  disjoint union of copies of the homogeneous space $(U/N)/(H/N)$ and hence a
  smooth manifold in a canonical way.
\end{proof}

\subsection{Functors on Kasparov categories}
\label{sec:KK_functors}

The (minimal) $C^*$\nbd{}tensor product gives rise to bifunctors
\begin{alignat*}{2}
  \KK^{G\cross X} \times \KK^G &\to \KK^{G\cross X},
  &\qquad (A,B) &\mapsto A\otimes B, \\
  \KK^{G\cross X} \times \KK^{G\cross X} &\to \KK^{G\cross X},
  &\qquad (A,B) &\mapsto A\otimes_X B,
\end{alignat*}
see \cite{Kasparov:Novikov}*{Definition 2.12 and Proposition 2.21}.  We
briefly recall how $A\otimes_X B$ looks like.  If $A,B\in\KK^{G\cross X}$,
then $A\otimes B$ is a $G\cross (X\times X)$-$C^*$\brd{}algebra, and
$A\otimes_X B$ is defined as its ``restriction'' to the diagonal.  That is, we
divide out elements of the form $f\cdot a$ with $f\in C_0(X\times X)$,
$f(x,x)=0$ for all $x\in X$, $a\in A\otimes_X B$.  See also
\cite{Kasparov:Novikov}*{Definition 1.6}.

The full and reduced descent functors $\KK^{G\cross X}\to \KK$ are defined in
\cite{Kasparov:Novikov}*{page 170--173}.  On objects, they act by $A\mapsto
(G\cross X)\cross A$ and $A\mapsto (G\cross X)\rcross A$.  We remark that the
space~$X$ has no effect here, that is, $(G\cross X)\cross A=G\cross A$ and
$(G\cross X)\rcross A=G\rcross A$ (see also
\cite{Chabert-Echterhoff-Oyono:Shapiro}).

Let $H\subseteq G$ be a closed subgroup.  Then we have functors
\begin{align*}
  \Res_G^H &\colon \KK^{G\cross X}\to \KK^{H\cross X}, \\
  \Ind_H^G &\colon \KK^{H\cross X}\to \KK^{G\cross X},
\end{align*}
called \emph{restriction} and \emph{induction}, respectively.  The restriction
functor is a special case of the functoriality of $\KK^{G\cross X}$ in~$G$:
any group homomorphism $H\to G$ induces a functor $\KK^{G\cross
  X}\to\KK^{H\cross X}$ by \cite{Kasparov:Novikov}*{Definition 3.1}.  The
induction functor is introduced in \cite{Kasparov:Novikov}*{Section 3.6}, see
also~\cite{Chabert-Echterhoff:Permanence}.

Finally, a $G$\nbd{}equivariant continuous map $f\colon X\to Y$ induces
functors
\begin{align*}
  f_* &\colon \KK^{G\cross X}\to \KK^{G\cross Y}, \\
  f^* &\colon \KK^{H\cross Y}\to \KK^{G\cross X}.
\end{align*}
The functor~$f_*$ is just a forgetful functor: to view a $G\cross
X$-$C^*$-algebra~$A$ as a $G\cross Y$-$C^*$-algebra, compose $f^*\colon
C_0(Y)\to C_b(X)$ and the canonical extension of the structural homomorphism
$C_b(X)\to Z\Mult(A)$.  The functor~$f^*$ is defined on objects by
$f^*(A)\defeq C_0(X)\otimes_Y A$.  Clearly, $\ID_*\cong\ID$, $\ID^*\cong\ID$
and $f_*g_*\cong (fg)_*$, $g^*f^*\cong (fg)^*$ if $f$ and~$g$ are composable.

Nowadays, we can treat these functors much more easily than
in~\cite{Kasparov:Novikov} using the universal property of Kasparov theory.
In the non-equivariant case, Nigel Higson has shown that the functor from
$C^*$-algebras to $\KK$ is the universal split exact stable homotopy functor,
that is, any functor from $C^*$-algebras to some category with these
properties factors uniquely through $\KK$.  This result has been extended by
Klaus Thomsen to the equivariant case and also works $G\cross
X$\nbd{}equivariantly by~\cite{Meyer:KKG}.  The above functors on
$\KK$\nbd{}categories all come from functors~$F$ between categories of
$C^*$-algebras, which are much easier to describe.

Let~$F$ be a functor from $G\cross X$-$C^*$-algebras to $H\cross
Y$-$C^*$-algebras.  The relevant functors~$F$ satisfy $F(A\otimes B) \approx
F(A)\otimes B$ for any nuclear $C^*$\nbd{}algebra~$B$ equipped with the
trivial representation of~$G$ (recall that $\approx$ denotes isomorphism of
$G\cross X$-$C^*$-algebras).  This implies immediately that~$F$ is stable and
homotopy invariant and commutes with suspensions.  Suppose, in addition,
that~$F$ maps extensions with a completely positive, contractive,
$G$\nbd{}equivariant linear section again to such extensions.  This is the
case in the above examples.  By the universal property, $F$ induces a functor
$\KK^{G\cross X}\to\KK^{H\cross Y}$.  Our mild exactness hypothesis guarantees
that~$F$ maps mapping cone triangles again to mapping cone triangles.  This
suffices to conclude that we have got a triangulated functor.  This argument
provides a very quick existence proof for the functors above and also shows
that they are triangulated.  It is also easy to check that they commute with
countable direct sums on $\KK^{G\cross X}$ (recall that direct sums in
$\KK^{G\cross X}$ are just $C^*$-direct sums).

Green's Imprimitivity Theorem and its reduced version assert that
\begin{equation}  \label{eq:Green}
  G\cross \Ind_H^G(A) \sim_M H\cross A,
  \qquad
  G\rcross \Ind_H^G(A) \sim_M H\rcross A.
\end{equation}

The functors $f_*$ and~$f^*$ are compatible with~$\otimes$ (without~$X$) in
the evident sense: $f_*(A\otimes B)\cong f_*(A)\otimes B$, $f^*(A\otimes
B)\cong f^*(A)\otimes B$.  We have a natural $G\cross Y$\nbd{}equivariant
isomorphism
\begin{equation}  \label{eq:adjointness_map}
  f_*(A) \otimes_Y B \approx f_*\bigl(A\otimes_X f^*(B)\bigr)
\end{equation}
for $f\colon X\to Y$ and $A\in\KK^{G\cross X}$, $B\in\KK^{G\cross Y}$
because~$\otimes_X$ is associative and $A\otimes_X C_0(X)\approx A$.
Equation~\eqref{eq:adjointness_map} asserts for the constant map $p_X\colon
X\to\pt$ that
\begin{equation}  \label{eq:tensorX_tensor}
  A\otimes_X p_X^*(B) \approx A\otimes B.
\end{equation}
The isomorphisms in \eqref{eq:adjointness_map} and~\eqref{eq:tensorX_tensor}
are natural, even in the formal sense.  For~\eqref{eq:tensorX_tensor},
naturality means that the isomorphisms intertwine
$$
x\otimes_X p_X^*(y)\in \KK^G(A\otimes_X p_X^*B, A'\otimes_X p_X^*B')
\quad \text{and} \quad
x\otimes y\in\KK^G(A\otimes B,A'\otimes B')
$$
if $x\in\KK^G(A,A')$, $y\in\KK^G(B,B')$.  By the universal property of
$\KK$, it suffices to verify this in the (easy) special case where $x$ and~$y$
are ordinary $*$\nbd{}homomorphisms; the general case then follows because two
functors agree on $\KK$ once they agree for ordinary $*$\nbd{}homomorphisms.
All the isomorphisms that follow are also natural in this sense, for the same
reason.

There are obvious compatibility conditions
\begin{equation}  \label{eq:res_tensorX}
  \Res_G^H(A\otimes_{(X)} B) \approx \Res_G^H(A)\otimes_{(X)}\Res_G^H(B),
\end{equation}
where we write $\otimes_{(X)}$ for either $\otimes_X$ or $\otimes$, and
\begin{equation}  \label{eq:resind_trivial_compatibility}
  \Ind_H^G\circ f_*\approx f_*\circ \Ind_H^G,
  \quad
  \Res_G^H\circ f_*\approx f_*\circ \Res_G^H,
  \quad
  \Res_G^H\circ f^*\approx f^*\circ \Res_G^H,
\end{equation}
because in each case one of the functors is a forgetful functor.  The
relation $\Ind_H^G\circ f^*\approx f^*\circ \Ind_H^G$ also holds.  The easiest
way to prove this isomorphism and many others is to replace $\Ind_H^G$ by a
forgetful functor as follows.

The groupoid $G\cross (G/H\times X)$ is Morita equivalent to $H\cross X$.
Therefore, the categories of $H\cross X$-$C^*$-algebras and of $G\cross
(G/H\times X)$-$C^*$-algebras are equivalent.  We may view induction as a
functor between these two categories.  This is an equivalence of categories.
Its inverse simply restricts a $G\cross (G/H\times X)$-$C^*$-algebra to
$\{H\}\times X\subseteq G/H\times X$.  By the universal property of $\KK$,
these functors induce an equivalence of categories $\KK^{G\cross (G/H\times
  X)} \cong \KK^{H\cross X}$ (see also~\cite{LeGall:KK_groupoid}).

Let $\pi_X\colon G/H\times X\to X$ be the projection.  When we reinterpret
$\Ind_H^G$ and $\Res_G^H$ as functors between $\KK^{G\cross X}$ and
$\KK^{G\cross (G/H\times X)}$, we get
\begin{equation}  \label{eq:resind_as_induced}
  \Ind_H^G \approx \pi_{X,*},
  \qquad
  \Res_G^H \approx \pi_X^*.
\end{equation}
This is useful for understanding the formal properties of these functors.
Using the properties of $f_*$ and $f^*$ shown above we get
\begin{align}
  \label{eq:ind_compatibility}
  \Ind_H^G\circ f^* &\approx f^*\circ \Ind_H^G,
  \\
  \label{eq:ind_tensorX}
  (\Ind_H^G A)\otimes_{(X)} B &\approx \Ind_H^G(A\otimes_{(X)}\Res_G^H B)
  \\
  \label{eq:ind_after_res}
  \Ind_H^G\circ \Res_G^H(A) &\approx C_0(G/H) \otimes A.
\end{align}

Our next goal is to prove the adjointness relation
\begin{equation}  \label{eq:map_induced_adjoint}
  \KK^{G\cross X}(f^*(A),B) \cong \KK^{G\cross Y}(A,f_*(B))
\end{equation}
for a \emph{proper} continuous $G$\nbd{}map $f\colon X\to Y$,
$A\in\KK^{G\cross Y}$, $B\in\KK^{G\cross X}$.  Experts on $\KK$\nbd{}theory
can verify~\eqref{eq:map_induced_adjoint} easily by showing that both sides
are defined by equivalent classes of cycles.  Category theorists may prefer
the following argument, which requires no knowledge of $\KK$ except the
existence of $f_*$ and~$f^*$ as functors on $\KK$.  Let $f\colon X\to Y$ be a
continuous $G$\nbd{}map.  Let~$B$ be a $G\cross X$-$C^*$-algebra.  There is a
natural homomorphism
$$
\pi_B\colon f^*f_*(B)\cong C_0(X)\otimes_Y B\to C_0(X)\otimes_X B \cong B.
$$
Let $A$ be a $G\cross Y$-$C^*$-algebra.  We have a natural map~$\iota_A$
from~$A$ to the multiplier algebra of $f^*(A)=f_*f^*(A)$, which sends
$a\mapsto 1\otimes a\in C_b(X)\otimes_Y A$ or, equivalently, $h\cdot a\mapsto
h\otimes_Y a$ for $h\in C_0(Y)$, $a\in A$.  The second description shows that
we have a map $\iota_A\colon A\to f_*f^*(A)$ if~$f$ is proper.  The composite
maps
\begin{gather*}
  f^*(A) \overset{f^*(\iota_A)}\longrightarrow f^*\bigl(f_*f^*(A)\bigr)
  = f^*f_*\bigl(f^*A) \overset{\pi_{f^*(A)}}\longrightarrow f^*A,
  \\
  f_*(B) \overset{\iota_{f_*B}}\longrightarrow f_*f^*(f_*B)
  = f_*(f^*f_*B) \overset{f_*(\pi_B)}\longrightarrow f_*(B)
\end{gather*}
are the identity.  Thus the maps $\pi_B$ and~$\iota_A$ form (co)units of
adjunction between the functors $f^*$ and~$f_*$ (see \cite{MacLane:Categories}
for this notion).  This holds regardless of whether we use homomorphism or
$\KK$\nbd{}morphisms.  Thus we get the desired adjointness
relation~\eqref{eq:map_induced_adjoint} and a corresponding statement about
equivariant $*$\nbd{}homomorphisms.

Combining~\eqref{eq:map_induced_adjoint} with~\eqref{eq:resind_as_induced}, we
get a Frobenius reciprocity isomorphism
\begin{equation}  \label{eq:res_cocompact}
  \KK^{G\cross X}(A,\Ind_H^G B) \cong \KK^{H\cross X}(\Res_G^H A,B)
\end{equation}
if $H\subseteq G$ is a \emph{cocompact} closed subgroup and $A\in\KK^{G\cross
  X}$, $B\in\KK^{H\cross X}$.  Dually, there is a natural isomorphism
\begin{equation}  \label{eq:ind_open}
  \KK^{G\cross X}(\Ind_U^G A,B) \cong \KK^{U\cross X}(A,\Res_G^U B)
\end{equation}
for an \emph{open} subgroup $U\subseteq G$.  This can also be proven by
writing down explicitly the units of adjunction.  We can decompose $\Res_G^U
\Ind_U^G(A)$ as a direct sum of $U\cross X$-$C^*$-algebras over the discrete
space of double cosets $G//U$.  The summand for the identity coset can be
identified with~$A$, so that we get a natural map $\iota_A\colon A\to \Res_G^U
\Ind_U^G(A)$.  We can represent $C_0(G/U)$ on the Hilbert space $\ell^2(G/U)$
by multiplication operators.  Since~$U$ is open in~$G$, this maps $C_0(G/U)$
into the $C^*$-algebra of compact operators $\Comp(\ell^2(G/U))$.  Hence we
get a natural morphism
$$
\Ind_U^G \Res_G^U (B)
\approx C_0(G/U)\otimes B
\to \Comp(\ell^2(G/U)) \otimes B
\sim_M B
$$
for $B\in\KK^{G\cross X}$.  This defines an element $\pi_B\in\KK^{G\cross
  X}(\Ind_U^G\Res_G^U(B),B)$ because $\KK$ is stable.  One verifies easily
that the morphisms $\pi_B$ and $\iota_A$ are units of adjunction, so that we
get~\eqref{eq:ind_open}.

After these purely formal manipulations of functors, we now come to a much
deeper assertion, which is due to Gennadi Kasparov.

\begin{proposition}[\cite{Kasparov:Novikov}*{Theorem 5.8}]
  \label{pro:ind_res}
  Let~$G$ be almost connected and let $H\subseteq G$ be a maximal compact
  subgroup.  If one of $X$, $A$ and~$B$ is a proper $G$\nbd{}space or a proper
  $G$-$C^*$\brd{}algebra, then
  \begin{equation}  \label{eq:res_mc}
    \Res_G^H\colon \KK^{G\cross X}(A,B)
    \to \KK^{H\cross X}(\Res_G^H A,\Res_G^H B)
  \end{equation}
  is an isomorphism.
\end{proposition}

\begin{lemma}  \label{lem:ind_res}
  Let $H\subseteq G$ be a large compact subgroup and let $U\defeq U_H$.  There
  is a natural isomorphism
  \begin{equation}  \label{eq:ind_res}
    \KK^{G\cross X}(\Ind_H^G A,B)
    \cong \KK^{H\cross X}(\Res_U^H \Ind_H^U A,\Res_G^H B).
  \end{equation}
  Define $J_H^G(A)\defeq \Ind_H^G (C_0((U/H)^7)\otimes A)$.  Then there
  is a natural isomorphism
  $$
  \KK^{G\cross X}(J_H^G A,B)\cong \KK^{H\cross X}(A,\Res_G^H B)
  $$
  if~$A$ belongs to the essential range of the functor $\Res_U^H$.
  Furthermore, the functors
  $$
  \Res_U^H\colon \KK^{U\cross X}\to\KK^{H\cross X},
  \qquad
  \Res_U^H\Ind_H^U\colon \KK^{H\cross X}\to\KK^{H\cross X},
  $$
  have the same essential range.
\end{lemma}

The \emph{essential range} of a functor $F\colon \Cat\to\Cat'$ is defined as
the class of all objects of~$\Cat'$ that are isomorphic to an object of the
form $F(X)$ with~$X$ an object of~$\Cat$.

\begin{proof}
  Induction in stages and~\eqref{eq:ind_open} yield
  $$
  \KK^{G\cross X}(\Ind_H^G A,B)
  \cong \KK^{U\cross X}(\Ind_H^U A,\Res_G^U B).
  $$
  Since $\Ind_H^U A$ is proper, Proposition~\ref{pro:ind_res}
  yields~\eqref{eq:ind_res}.  Let us abbreviate $R\defeq \Res_U^H$ and
  $I\defeq \Ind_H^U$.  If~$A$ belongs to the essential range of~$R$, then
  $RI(A) \cong C_0(U/H)\otimes A$ by~\eqref{eq:ind_after_res}.  Recall that
  $U/H$ is $H$\nbd{}equivariantly diffeomorphic to a real vector space with a
  linear action of~$H$.  Hence Bott periodicity provides a
  $\KK^H$\nbd{}equivalence between $(RI)^8(A) \cong C_0((U/H)^8) \otimes A$
  and~$A$.  This yields the second isomorphism and shows that the essential
  range of~$R$ is contained in the essential range of $RI$.  The converse
  inclusion is trivial.
\end{proof}

Hence $\Ind_H^G$ and $\Res_G^H$ for a large compact subgroup $H\subseteq G$
become adjoint functors if we replace $\KK^{H\cross X}$ by the essential range
of $\Res_U^H$.  There is no analogue of this for arbitrary compact subgroups.

\section{A decomposition of the Kasparov category}
\label{sec:decomposition}

\begin{definition}  \label{def:CC}
  We call $A\in\KK^{G\cross X}$ \emph{weakly contractible} if
  $\Res_G^H(A)\cong 0$ for all compact subgroups $H\subseteq G$.  Let
  $\CC\subseteq\KK^{G\cross X}$ be the full subcategory of weakly contractible
  objects.
  
  A morphism~$f$ in $\KK^{G\cross X}(A,B)$ is called a \emph{weak
    equivalence} if $\Res_G^H(f)$ is invertible in $\KK^{H\cross
    X}$ for all compact subgroups $H\subseteq G$.  We say that~$f$
  \emph{vanishes for compact subgroups} if $\Res_G^H(f)=0$ for all
  compact subgroups $H\subseteq G$.
  
  We call a $G\cross X$-$C^*$\brd{}algebra \emph{compactly induced} if it is
  isomorphic in $\KK^{G\cross X}$ to $\Ind_H^G(A)$ for some compact subgroup
  $H\subseteq G$ and some $A\in\KK^{H\cross X}$.  We let
  $\CI\subseteq\KK^{G\cross X}$ be the full subcategory of compactly induced
  objects.
\end{definition}

In all these definitions, it suffices to consider large compact subgroups
because any compact subgroup is contained in a large one by
Lemma~\ref{lem:MCAC_enough}.  Our next goal is to prove that $(\gen{\CI},\CC)$
is a complementary pair of localising subcategories of $\KK^{G\cross X}$, so
that we can apply Proposition~\ref{pro:complementary}.

\begin{lemma}  \label{lem:CC_CI_localising}
  The subcategories $\CC$ and $\gen{\CI}$ of $\KK^{G\cross X}$ are localising.
  
  The subcategories $\CC$, $\CI$ and $\gen{\CI}$ are closed under tensor
  products with arbitrary objects of $\KK^G$ and $\KK^{G\cross X}$.
\end{lemma}

\begin{proof}
  Since the functor $\Res_G^H$ is triangulated and commutes with direct sums,
  its kernel is localising.  Hence $\CC$ is localising as an intersection of
  localising subcategories.  The subcategory $\gen{\CI}$ is localising by
  construction.  The subcategories $\CC$ and $\CI$ are closed under tensor
  products because of the compatibility of restriction and induction with
  tensor products discussed in Section~\ref{sec:preliminaries}.  Since the
  functor $\blank\otimes_{(X)} B$ is triangulated and commutes with direct
  sums, it leaves $\gen{\CI}$ invariant as well.
\end{proof}

\begin{lemma}  \label{lem:weak_contractible_equivalence}
  A morphism in $\KK^{G\cross X}$ is a weak equivalence if and only if its
  mapping cone is weakly contractible.
\end{lemma}

\begin{proof}
  Since the functor $\Res_G^H$ is triangulated, it maps an exact triangle
  $\Sigma B\to C\to A\overset{f}\to B$ again to an exact triangle.
  Lemma~\ref{lem:triangle_lemma} implies that $\Res_G^H f$ is invertible if
  and only if $\Res_G^H C\cong0$.
\end{proof}

\begin{proposition}  \label{pro:CC_CI_orthogonal}
  An object~$N$ of $\KK^{G\cross X}$ is weakly contractible if and only if
  $\KK^{G\cross X}(P,N)\cong 0$ for all $P\in\CI$.
  
  A morphism $f\in\KK^{G\cross X}(A,B)$ is a weak equivalence if and only if
  it induces an isomorphism $f_*\colon \KK^{G\cross X}(P,A) \to \KK^{G\cross
    X} (P,B)$ for all $P\in\CI$.
  
  A morphism $f\in\KK^{G\cross X}(A,B)$ vanishes for compact subgroups if and
  only if it induces the zero map $f_*\colon \KK^{G\cross
    X}(P,A)\to\KK^{G\cross X}(P,B)$ for all $P\in\CI$.
  
  In the first two assertions, we can replace $\CI$ by $\gen{\CI}$.
\end{proposition}

\begin{proof}
  Let $H\in\MCAC$ be maximal compact and let $U\defeq U_H$.  Let $P=\Ind_H^G
  A$ for some $A\in\KK^{H\cross X}$.  Any object of $\CI$ is of this form for
  some $H$, $A$ by Lemma~\ref{lem:MCAC_enough}.  We use~\eqref{eq:ind_res} to
  rewrite $\KK^{G\cross X}(P,N)\cong\KK^{H\cross X}(A',\Res_G^H N)$ with
  $A'\defeq\Res_U^H\Ind_H^U A$.  If $\Res_G^H N\cong0$, then the right hand
  side vanishes, so that $\KK^{G\cross X}(P,N)=0$.  Conversely, if
  $\KK^{G\cross X}(P,N)=0$ for all $P\in\CI$, then $\KK^{H\cross X}(\Res_G^H
  N,\Res_G^H N)=0$ and hence $\Res_G^H N=0$.  We have used that $\Res_G^H N$
  belongs to the essential range of $\Res_U^H\Ind_H^U$ by
  Lemma~\ref{lem:ind_res}.  This proves the first assertion.  We can enlarge
  $\CI$ to $\gen{\CI}$ because the class of objects~$P$ with $\KK^{G\cross
    X}(P,N)=0$ for all $N\in\CC$ is localising.  The remaining assertions are
  proven similarly.
\end{proof}

\begin{definition}  \label{def:approximation}
  A \emph{$\CI$\nbd{}simplicial approximation} of $A\in\KK^{G\cross X}$ is a
  weak equivalence $\tilde{A}\to A$ with $\tilde{A}\in\gen{\CI}$.  A
  $\CI$\brd{}simplicial approximation of $C_0(X)$ is also called a \emph{Dirac
    morphism} for $G\cross X$.
\end{definition}

\begin{proposition}  \label{pro:Dirac}
  A Dirac morphism exists for any $G\cross X$.
\end{proposition}

The existence of a Dirac morphism is the main technical result needed for our
approach to the Baum-Connes conjecture.  Logically, we should now prove the
existence of the Dirac morphism (we postpone this until
Section~\ref{sec:representability}) and then compute the localisation
$\KK^{G\cross X}/\CC$ (we do this in Section~\ref{sec:derived_proper}) before
we dare to localise the functor $\K_*(G\rcross\blank)$.  Instead, we head for
the Baum-Connes assembly map as quickly as possible.

The following theorem uses the notation of
Proposition~\ref{pro:complementary}.

\begin{theorem}  \label{the:KKG_decomposition}
  The localising subcategories $\gen{\CI},\CC$ of $\KK^{G\cross X}$ are
  complementary.  Let $\Dirac\in\KK^{G\cross X}(\ADir,C_0(X))$ be a Dirac
  morphism for $G\cross X$ and form the exact triangle
  \begin{equation}  \label{eq:Dirac_triangle}
     \Sigma\AN \to \ADir \overset{\Dirac}\to C_0(X) \to \AN.
  \end{equation}
  Then $P(A)\cong \ADir\otimes_X A$ and $N(A)\cong \AN \otimes_X A$ naturally
  for all $A\in\KK^{G\cross X}$, and the natural transformations $\Sigma
  N(A)\to P(A)\to A\to N(A)$ are induced by the maps
  in~\eqref{eq:Dirac_triangle}.  We have $A\in\gen{\CI}$ if and only if
  $\KK^{G\cross X}(A,B)=0$ for all $B\in\CC$ if and only if $\ADir\otimes_X
  A\cong A$; and $B\in\CC$ if and only if $\KK^{G\cross X}(A,B)=0$ for all
  $A\in\gen{\CI}$ if and only if $\ADir\otimes_X B\cong0$.  In particular,
  $\ADir\otimes_X\ADir\cong\ADir$.
\end{theorem}

\begin{proof}
  Since $\Dirac$ is a weak equivalence, $\AN$ is weakly contractible by
  Lemma~\ref{lem:weak_contractible_equivalence}.  The tensor product
  of~\eqref{eq:Dirac_triangle} with~$A$ is another exact triangle
  because $\blank\otimes_X A$ is a triangulated functor.  Since $\CC$ and
  $\gen{\CI}$ are closed under tensor products by
  Lemma~\ref{lem:CC_CI_localising}, we get an exact triangle as in the
  definition of a complementary pair of subcategories in
  Section~\ref{sec:localisation}.  This yields the assertions in the first
  paragraph.  Those in the second paragraph follow from
  Proposition~\ref{pro:complementary}.
\end{proof}

\begin{definition}  \label{def:Dirac_triangle}
  An exact triangle as in~\eqref{eq:Dirac_triangle} is called a
  \emph{Dirac triangle}.
\end{definition}

Using Proposition~\ref{pro:complementary}.\ref{pro::PN_adjoint}, we can now
compute localisations and obstruction functors from a Dirac triangle.  The
morphisms in $\KK^G/\CC$ are given by
$$
\KK^{G\cross X}/\CC (A,B) \cong \KK^{G\cross X} (\ADir\otimes_X A,B).
$$
The localisation and the obstruction functor of a functor
$F\colon\KK^{G\cross X}\to\Cat$ are
$$
\Left F(A)\cong F(\ADir\otimes_X A),
\qquad
\Obs F(A)\cong F(\AN\otimes_X A),
$$
and the natural transformations $\Left F(A)\to F(A)\to \Obs F(A)$ are
induced by the maps $\ADir\to\pt\to\AN$ in the Dirac triangle.

We are particularly interested in the functor
$$
\KK^{G\cross X}\to\KK,\qquad  A\mapsto (G\cross X)\rcross A.
$$
We denote its localisation and obstruction functor by $A\mapsto (G\cross
X)\Lrcross A$ and $A\mapsto (G\rcross X)\Obscross A$, respectively.

\section{The Baum-Connes assembly map}
\label{sec:assembly}

We now relate our analysis of $\KK^{G\cross X}$ to the Baum-Connes assembly
map.  Since we consider transformation groups $G\cross X$, we first have to do
some work related to the space~$X$.  Jérôme Chabert, Siegfried Echterhoff and
Hervé Oyono-Oyono show in~\cite{Chabert-Echterhoff-Oyono:Shapiro} that there
is a commutative diagram
$$
\xymatrix{
  {\Ktop_*(G\cross X,A)} \ar[r]^{\mu_A} \ar[d]^{\cong} &
  {\K_*((G\cross X)\rcross A)} \ar[d]^{\cong} \\
  {\Ktop_*(G,A)} \ar[r]^{\mu_A} &
  {\K_*(G\rcross A).}
}
$$
That is, the Baum-Connes assembly map just ignores the space~$X$.  We need
a similar result in our setup.

\begin{lemma}  \label{lem:forget_X}
  The functor $p_X^*\colon \KK^G\to\KK^{G\cross X}$ maps $\CC,\CI,\gen{\CI}
  \subseteq\KK^G$ to the corresponding subcategories in $\KK^{G\cross X}$.  If
  $f\in\KK^G(A,B)$ is a weak equivalence or vanishes for compact subgroups, so
  does $p_X^*(f)$.  If $\Dirac\in\KK^G(\ADir,\pt)$ is a Dirac morphism
  for~$G$, then $p_X^*(\Dirac)\in\KK^{G\cross X}(C_0(X,\ADir),C_0(X))$ is a
  Dirac morphism for $G\cross X$.  There are natural isomorphisms
  $$
  (G\cross X)\Lrcross A\cong G\Lrcross A,
  \qquad
  (G\cross X)\Obscross A\cong G\Obscross A.
  $$
\end{lemma}

\begin{proof}
  Recall from Section~\ref{sec:preliminaries} that the functor~$p_X^*$ is
  compatible with restriction and induction.  This implies
  $p_X^*(\CC)\subseteq\CC$ and $p_X^*(\CI)\subseteq\CI$.  The same holds for
  $\gen{\CI}$ because~$p_X^*$ is triangulated and commutes with direct sums.
  This implies the assertions about weak equivalences and Dirac morphisms.
  Now~\eqref{eq:tensorX_tensor} yields
  $$
  (G\cross X)\Lrcross A
  \cong (G\cross X) \rcross (p_X^*(\ADir) \otimes_X A)
  \approx G\rcross (\ADir\otimes A)
  \cong G\Lrcross A.
  $$
  For the same reason, $(G\cross X)\Obscross A\cong G\Obscross A$.
\end{proof}

For our purposes, we do not need the details of the standard definition of the
Baum-Connes conjecture.  We only need to know the following two facts: the
Baum-Connes conjecture holds for compactly induced coefficient algebras (see
\cite{Chabert-Echterhoff:Permanence}), and weak equivalences induce
isomorphisms on $\Ktop_*(G,\blank)$ (see
\cite{Chabert-Echterhoff-Oyono:Going_down}).  This second assertion also
follows immediately from the definition of $\Ktop_*$ and
Corollary~\ref{cor:CC_proper} below.  Thus the only substantial result about
the standard definition of the Baum-Connes conjecture that we have to import
is that it holds for compactly induced coefficient algebras.

\begin{theorem}  \label{the:Baum-Connes_assembly}
  Let $\tilde{A}\to A$ be a $\CI$\brd{}simplicial approximation of
  $A\in\KK^{G\cross X}$.  Then the indicated maps in the commutative diagram
  $$
  \xymatrix{
    {\Ktop_*(G\cross X,\tilde{A})} \ar[r]^{\cong}
    \ar[d]^{\mu_{\tilde{A}}}_{\cong} &
    {\Ktop_*(G\cross X,A)} \ar[d]^{\mu_A} \\
    {\K_*((G\cross X)\rcross \tilde{A})} \ar[r] &
    {\K_*((G\cross X)\rcross A)}
  }
  $$
  are isomorphisms.  Hence the Baum-Connes assembly map is naturally
  isomorphic to the canonical map $\K_*((G\cross X)\Lrcross A)\to
  \K_*((G\cross X)\rcross A)$.  It is an isomorphism if and only if
  $\K_*((G\cross X)\Obscross A)\cong 0$.
\end{theorem}

\begin{proof}
  Lemma~\ref{lem:forget_X} shows that we may assume without loss of generality
  that $X=\pt$.  The left vertical map is the assembly map for the coefficient
  algebra~$\tilde{A}$.  Since $\Ktop_*(G,\blank)$ is a homological functor
  that commutes with direct sums, the class of coefficient algebras~$B$ for
  which~$\mu_B$ is an isomorphism is a localising subcategory of $\KK^G$.
  Therefore, the Baum-Connes conjecture holds for all coefficient algebras in
  $\gen{\CI}$ because it holds for compactly induced coefficient algebras
  by~\cite{Chabert-Echterhoff:Permanence}.  As a result, the left vertical map
  in our diagram is an isomorphism.  It is shown
  in~\cite{Chabert-Echterhoff-Oyono:Going_down} that weak equivalences induce
  isomorphisms on $\Ktop_*(G,\blank)$.  Hence the top horizontal map is an
  isomorphism as well.
\end{proof}

Therefore, it is legitimate to view the map $\Left F(A)\to F(A)$ for a
covariant functor~$F$ defined on $\KK^G$ as an assembly map for $F(A)$.

As we explained in Section~\ref{sec:localisation}, there are two pictures of
$\Left F(A)$: either as the best possible approximation to $F(A)$ that
vanishes on $\CC$ or as the best possible approximation to $F(A)$ that only
uses the values $F(B)$ for $B\in\gen{\CI}$.  In particular, the map $\Left
F(A)\to F(A)$ is an isomorphism for all $A\in\KK^{G\cross X}$ if and only if
$F|_{\CC}=0$.  As we remarked in the introduction, this yields an equivalent,
elementary formulation of the Baum-Connes conjecture when applied to
$\K_*(G\rcross A)$.

Alain Connes has asked recently whether it is possible to improve upon the
Baum-Connes conjecture, finding better approximations to $\K_*(G\rcross A)$.
Like the approaches in \cites{Balmer-Matthey:Foundations,
  Balmer-Matthey:Cofibrant, Balmer-Matthey:Model_BC, Davis-Lueck:Assembly},
our construction of the assembly map is sufficiently flexible to say something
about this, though our answer may not be very satisfactory.  The Baum-Connes
conjecture asserts that the objects of $\gen{\CI}$ are general enough to
predict everything that happens in the $\K$\nbd{}theory of reduced crossed
products.  If it fails, this means that there are some phenomena in
$\K_*(G\rcross A)$ that do not yet occur for $A\in\gen{\CI}$.  To get a better
conjecture, we have to add some of the coefficient algebras for which
Baum-Connes fails to $\CI$.  Then the general machinery of localisation yields
again a best possible approximation to $\K_*(G\rcross A)$ based on what
happens for coefficients in $\gen{\CI'}$.  The new conjecture expresses
$\K_*(G\rcross A)$ for arbitrary~$A$ in terms of $\K_*(G\rcross A)$ for
$A\in\CI'$.  However, such a reduction of the problem is only as good as our
understanding of what happens for $A\in\CI'$.  At the moment, it does not seem
that we have a sufficient understanding of the failure of the Baum-Connes
conjecture to make any progress in this direction.

\section{The Brown Representability Theorem and the Dirac morphism}
\label{sec:representability}

Recall that a morphism $\Dirac\in\KK^G(\ADir,\pt)$ is a weak equivalence if
and only if the induced map $\KK^G(A,\ADir)\to\KK^G(A,\pt)$ is an isomorphism
for all $A\in\gen{\CI}$.  Since~$\ADir$ is supposed to lie in the same
subcategory $\gen{\CI}$, the Dirac morphism exists if and only if the functor
$A\mapsto \KK^G(A,\pt)$ on the category $\gen{\CI}$ is representable.  In the
classical case of simplicial approximation of arbitrary topological spaces by
simplicial complexes, one can either write down explicitly such a representing
object or appeal to the Brown Representability Theorem.  We shall prove the
existence of the Dirac morphism using the second method.

There are several representability theorems for triangulated categories that
use different hypotheses.  It seems that none of them applies directly to the
category $\gen{\CI}$ that we need.  To circumvent this, we choose a smaller
set of generators $\CIZ\subseteq\CI$ which is small enough so that a general
representability theorem is available in $\gen{\CIZ}$ and large enough so that
the representing object in $\gen{\CIZ}$ actually represents the functor on the
whole of $\gen{\CI}$.  A byproduct of this proof technique is that we get
$\ADir\in\gen{\CIZ}$.  This is used in Section~\ref{sec:strong_BC}.

Since $\KK^G$ only has countable direct sums, we have to do some cardinality
bookkeeping.  Let~$\Tri$ be a triangulated category and let~$\aleph$ be an
infinite regular cardinal number.  We will only need the countable cardinal
number~$\aleph_0$.  We suppose that~$\Tri$ has direct sums of
cardinality~$\aleph$.

Recall that $\Ab$ denotes the category of Abelian groups.  A contravariant
functor $F\colon\Tri\to\Ab$ is called \emph{representable} if it is isomorphic
to $X\mapsto \Tri(X,Y)$ for some $Y\in\Tri$.  Representable functors are
cohomological and compatible with direct sums of cardinality~$\aleph$.  We now
formulate conditions on~$\Tri$ that ensure that these necessary conditions
plus an extra cardinality hypothesis are also sufficient.

An object $X\in\Tri$ is called \emph{$\aleph$\nbd{}compact} if $\Tri(X,Y)$ has
cardinality at most~$\aleph$ for all $Y\in\Tri$ and the covariant functor
$\Tri(X,\blank)$ is compatible with direct sums of cardinality~$\aleph$.  The
reader who consults~\cite{Neeman:Triangulated} on direct sums and
representability should beware that our notation differs slightly.  The axiom
(TR5$\aleph$) and the notion of $\aleph$\nbd{}compactness
in~\cite{Neeman:Triangulated} deal with direct sums of cardinality
\emph{strictly less} than~$\aleph$.

\begin{theorem}  \label{the:Brown_rep}
  Let~$\aleph$ be an infinite regular cardinal number and let~$\Tri$ be a
  triangulated category with direct sums of cardinality~$\aleph$.  Let~$\GEN$
  be a set of $\aleph$\nbd{}compact objects of~$\Tri$ with
  $\abs{\GEN}\le\aleph$.  Suppose that $\Tri(X,Y)=0$ for all $X\in \GEN$
  already implies $Y=0$.  Let $F\colon \Tri\to\Ab$ be an additive,
  contravariant functor.

  Then~$F$ is representable if and only if it satisfies the following
  conditions:
  \begin{enumerate}[(i)]
  \item $F$ is cohomological;
  \item $F$ is compatible with $\aleph$\nbd{}direct sums;
  \item $F(C)$ has cardinality at most~$\aleph$ for all $C\in \GEN$.
  \end{enumerate}
  Moreover, the hypothesis that $\Tri(X,Y)=0$ for all $X\in \GEN$ implies
  $Y=0$ can be replaced by the hypothesis that $\Tri=\gen{\GEN}^\aleph$.
\end{theorem}

\begin{proof}
  The conditions (i)--(iii) are clearly necessary.  The interesting assertion
  is that they are also sufficient.  If we leave out the cardinality
  restriction~$\aleph$, this is proven by Neeman in
  \cite{Neeman:Grothendieck_duality}, and it also follows from
  \cite{Neeman:Triangulated}*{Theorem 8.3.3}.  Since we do not have direct
  sums of arbitrary cardinality in~$\Tri$, we have to check that the proof
  does not require direct sums of cardinality strictly greater than~$\aleph$.
  This is indeed the case, as the determined reader may check for himself.  It
  turns out that the largest sums we need have cardinality at most
  $\N\times\aleph\times\aleph$.  This is dominated by~$\aleph$
  because~$\aleph$ is a regular cardinal.
\end{proof}

\subsection{Construction of Dirac morphisms}
\label{sec:construct_Dirac}

Now we introduce a subcategory of $\gen{\CI}$ with a hand-selected set of
generators $\CIZ$.  Let $H\in\MCAC$ and let $U\defeq U_H$ as in
Lemma~\ref{lem:MCAC_enough}.  We define $J_H^G$ as in
Lemma~\ref{lem:ind_res} and let
\begin{equation}  \label{eq:RC_concrete}
  \RC_H\defeq
  J_H^G(\pt) =
  \Ind_H^G C_0\bigl((U/H)^7\bigr) \cong
  C_0(G\times_H (U/H)^7).
\end{equation}
This compactly induced $G$\nbd{}$C^*$\brd{}algebra satisfies
\begin{equation}  \label{eq:RC_represents}
  \KK^G(\RC_H,B) \cong
  \KK^H(\pt,B) \cong
  \K(H\cross B).
\end{equation}
by Lemma~\ref{lem:ind_res}.  Equation~\eqref{eq:RC_represents} says that
$\RC_H$ (co)represents the covariant functor $\K(H\cross\blank)$ and thus
determines~$\RC_H$ uniquely up to $\KK^G$\brd{}equivalence.
Equation~\eqref{eq:RC_concrete} merely is a convenient choice of representing
object.  For an arbitrary compact subgroup, the functor $\K(H\cross\blank)$
may fail to be representable.  This is why we work with large compact
subgroups.

If $H,H'\subseteq U$ are two maximal compact subgroups, they are conjugate
in~$U$ by Lemma~\ref{lem:MCAC_enough}, so that the $G$-$C^*$\brd{}algebras
$\RC_H$ and~$\RC_{H'}$ are isomorphic.  Hence it suffices to choose one
maximal compact subgroup in each almost connected open subgroup.  We let
$\CIZ$ be the set of~$\RC_H$ for the chosen subgroups~$H$.

\begin{lemma}  \label{lem:CIZ_compact}
  The set $\CIZ$ is (at most) countable and consists of
  $\aleph_0$\nbd{}compact objects of $\KK^G$, where~$\aleph_0$ denotes the
  countable cardinal.
\end{lemma}

\begin{proof}
  It is well-known that $\K_*(A)$ is countable if~$A$ is a separable
  $C^*$-algebra and that $\K$\nbd{}theory for $C^*$\nbd{}algebras commutes
  with direct sums (\cite{Blackadar:Book}).  Hence~\eqref{eq:RC_represents}
  implies that $\KK^G(\RC_H,B)$ is countable for each $B\in\KK^G$ and that
  $\KK^G(\RC_H,\blank)$ commutes with countable direct sums.  That is, $\CIZ$
  consists of $\aleph_0$\nbd{}compact objects.  Since the objects of $\CIZ$
  are in bijection with the almost connected open subgroups of~$G$, it remains
  to prove that there are at most countably many such subgroups.  Let
  $G_0\subseteq G$ be the connected component of the identity element.  The
  open almost connected subgroups of~$G$ are in bijection with the compact
  open subgroups of $G/G_0$.  Since the latter group is second countable as a
  topological space, even the set of all compact open sub\emph{sets} of
  $G/G_0$ is countable.
\end{proof}

\begin{corollary}  \label{cor:CIZ_rep}
  For any $B\in\KK^G$, there is $\tilde{B}\in\gen{\CIZ}$ and
  $f\in\KK^G(\tilde{B},B)$ such that $f_*\colon
  \KK^G(A,\tilde{B})\to\KK^G(A,B)$ is an isomorphism for all
  $A\in\gen{\CIZ}$.
\end{corollary}

\begin{proof}
  Lemma~\ref{lem:CIZ_compact} implies that the category $\gen{\CIZ}$ with
  generating set $\CIZ$ satisfies the conditions of
  Theorem~\ref{the:Brown_rep}.  The functor $F(A)\defeq \KK^G(A,B)$ fulfils
  the necessary and sufficient conditions for representability because it is
  already represented on the larger category $\KK^G$.  Hence it is
  representable on $\gen{\CIZ}$.
\end{proof}

We can now prove the existence of the Dirac morphism.

\begin{proof}[Proof of Proposition~\ref{pro:Dirac}]
  We may assume $X=\pt$ by Lemma~\ref{lem:forget_X}.
  Corollary~\ref{cor:CIZ_rep} for $B=\pt$ yields $\Dirac\in\KK^G(\ADir,\pt)$
  with $\ADir\in\gen{\CIZ}\subseteq\gen{\CI}$.  We claim that~$\Dirac$ is a
  weak equivalence.  Any compact subgroup $H\subseteq G$ is subconjugate to a
  subgroup $L\subseteq G$ with $\RC_L\in\CIZ$ by Lemma~\ref{lem:MCAC_enough}.
  Hence it suffices to prove that $\Res_G^L(\Dirac)$ is invertible for
  those~$L$.  By construction of~$\Dirac$, it induces an isomorphism on
  $\KK^G_*(\RC_L,\blank)$.  Equation~\eqref{eq:RC_represents} yields
  that~$\Dirac$ induces an isomorphism $\KK^L(\pt,\ADir)\to\KK^L(\pt,\pt)$.
  To check that~$\Dirac$ is an isomorphism in $\KK^L$, it suffices to check
  that~$\Dirac$ induces an isomorphism $\KK^L(\ADir,\ADir)\to\KK^L(\ADir,\pt)$
  as well.  Since $\ADir\in\gen{\CIZ}$, this follows if we have isomorphisms
  $\KK^L_*(A,\ADir)\cong\KK^L_*(A,\pt)$ for all $A\in\CIZ$.  Thus we have to
  fix another large compact subgroup~$H$ and show that~$\Dirac$ induces an
  isomorphism $\KK^L_*(\RC_H,\ADir)\to\KK^L_*(\RC_H,\pt)$.
  
  Let $V\defeq U/H$, then we have $\RC_H=C_0(G\times_H V^7) \cong
  C_0(G\times_U V^8)$.  We only need the action of~$L$ on this space.  The
  space $G\times_U V^8$ decomposes into a disjoint union of the spaces
  $LgU\times_U V^8$ over the double cosets $g\in L\backslash G/U$.  We have a
  natural isomorphism $LgU\times_U V^8\cong L\times_{L\cap gUg^{-1}} V^8$,
  where we use the conjugation automorphism $gUg^{-1}\cong U$ to let
  $gUg^{-1}$ act on~$V$.  Since the action of $L\cap gUg^{-1}$ on~$V$ is
  diffeomorphic to a linear action, equivariant Bott periodicity yields that
  $C_0(LgU\times_U V^8)$ is $\KK^L$\brd{}equivalent to $\Ind_{L\cap
  gUg^{-1}}^L(\pt)$.  Thus
  $$
  \Res_G^L \RC_H
  \cong \bigoplus_{g\in L\backslash G/U} \Ind_{L\cap gUg^{-1}}^L(\pt)
  $$
  The subgroups $L\cap gUg^{-1}$ are open in~$L$ and again large by
  Lemma~\ref{lem:MCAC_enough}.  It follows from~\eqref{eq:ind_open} that
  \begin{multline*}
    \KK^L_*(\RC_H,B)
    \cong \bigoplus_{g\in L\backslash G/U}
    \KK^L_*(\Ind_{L\cap gUg^{-1}}^L(\pt),B)
    \\ \cong \bigoplus_{g\in L\backslash G/U}
    \KK^{L\cap gUg^{-1}}_*(\pt,B)
    \cong
    \bigoplus_{g\in L\backslash G/U} \KK^G_*(\RC_{L\cap gUg^{-1}},B).
  \end{multline*}
  By construction, $\Dirac$ induces an isomorphism on the right hand side and
  hence on $\KK^L_*(\RC_H,\blank)$.
\end{proof}

\begin{remark}
  Incidentally, the proof above shows that
  $$
  J_L^G\Res_G^L J_H^G(\pt) =
  J_L^G\Res_G^L(\RC_H)
  \cong \bigoplus_{g\in L\backslash G/U} \RC_{L\cap gUg^{-1}}.
  $$
\end{remark}

\subsection{A localisation related to the Universal Coefficient
  Theorem}
\label{sec:UCT}

\begin{definition}  \label{def:K_contractible}
  A $C^*$\nbd{}algebra~$A$ is called \emph{$\K$\nbd{}contractible} if
  $\K_*(A)=0$.  A morphism $f\in\KK(A,B)$ is called a
  \emph{$\K$\nbd{}equivalence} if $f_*\colon \K_*(A)\to\K_*(B)$ vanishes.
\end{definition}

We write $\Null\subseteq\KK$ for the full subcategory of
$\K$\nbd{}contractible objects.  This subcategory is localising because
$\K$\nbd{}theory is a homological functor compatible with direct sums.  A
morphism is a $\K$\nbd{}equivalence if and only if its mapping cone is
$\K$\nbd{}contractible.

\begin{theorem}  \label{the:KK_decomposition}
  The localising subcategories $\gen{\pt}$ and~$\Null$ in $\KK$ are
  complementary.
\end{theorem}

\begin{proof}
  We have $B\in\Null$ if and only if $\KK_*(\pt,B)\cong\K_*(B)=0$ for all
  $*\in\Z$, if and only if $\KK(A,B)=0$ for all $A\in\gen{\pt}$.  Similarly,
  $f\in\KK(\tilde{B},B)$ is a $\K$\nbd{}equivalence if and only if $f_*\colon
  \KK(A,\tilde{B})\to\KK(A,B)$ is an isomorphism for all $A\in\gen{\pt}$.  We
  have to construct a $\K$\nbd{}equivalence $\tilde{B}\to B$ with
  $\tilde{B}\in\gen{\pt}$ for any $B\in\KK$.  This is equivalent to the
  representability of the functor $A\mapsto \KK(A,B)$ on $\gen{\pt}$.  The
  object~$\pt$ of $\KK$ is clearly compact and generates $\gen{\pt}$ by
  definition.  Hence we can apply the Brown Representability
  Theorem~\ref{the:Brown_rep} to get the assertion.
\end{proof}

As we observed in Section~\ref{sec:functors_subcategories},
$\gen{\pt}\subseteq\KK$ is just the bootstrap category.  The simplicial
approximations in this context are usually called \emph{geometric resolutions}
(see~\cite{Blackadar:Book}).  Let $\UCT\defeq\KK/\Null$ be the localisation of
$\KK$ at the $\K$\nbd{}contractible objects.  This is a triangulated category
with countable direct sums and equipped with a triangulated functor
$\KK\to\UCT$ commuting with direct sums.  It has the same objects as $\KK$.
Morphisms are computed using geometric resolutions:
$$
\UCT(A,B)\cong \KK(\tilde{A},B)\cong\KK(\tilde{A},\tilde{B}).
$$
The group $\UCT(A,B)$ can always be computed using the Universal
Coefficient Theorem (see \cite{Blackadar:Book}) because the latter applies to
$\tilde{A}\in\gen{\pt}$.  Moreover, $A$ satisfies the Universal Coefficient
Theorem if and only if $\KK(A,B)\cong\UCT(A,B)$ for all~$B$ if and only if
$A\in\gen{\pt}$.  Thus we have translated the Universal Coefficient Theorem as
an isomorphism statement.  This is often convenient.

The functor $A\otimes\blank$ preserves $\K$\nbd{}equivalences for
$A\in\gen{\pt}$ because this holds for the generator~$\pt$.  Hence the natural
maps from $\tilde{A}\otimes\tilde{B}$ to $A\otimes\tilde{B}$ and
$\tilde{A}\otimes B$ are both $\K$\nbd{}equivalences, so that the various ways
of localising $A\otimes B$ give the same result $A\Lotimes B$ in the category
$\UCT$.  Since $A\Lotimes B$ only involves $C^*$\nbd{}algebras from the
bootstrap category, $\K_*(A\Lotimes B)$ can always be computed by the Künneth
Formula (see \cite{Blackadar:Book}).  We remark also that
$\gen{\pt}\otimes\gen{\pt}\subseteq\gen{\pt}$ because $\pt\otimes\pt=\pt$.

Thus localisation of $\KK$ at~$\Null$ yields a natural map
$\KK_*(A,B)\to\UCT_*(A,B)$, which is an isomorphisms for all~$B$ if and only
if~$A$ satisfies the Universal Coefficient Theorem, and a natural map
$\K_*(A\Lotimes B)\to \K_*(A\otimes B)$, which is an isomorphism for all~$B$
if and only if~$A$ satisfies the Künneth Formula.  We want to emphasise that,
on a formal level, these maps are analogous to the Baum-Connes assembly map
$\K_*(G\Lrcross A)\to \K_*(G\rcross A)$.

\section{The derived category and proper actions}
\label{sec:derived_proper}

We want to describe the localisation of $\KK^{G\cross X}$ at $\CC$ in more
classical terms.  Let $A$ and~$B$ be $G\cross X$-$C^*$\brd{}algebras.  Let~$Y$
be a locally compact $G$\nbd{}space.  Generalising~\eqref{eq:RKKG} slightly,
we let
$$
\RKK^{G\cross X}(Y;A,B)\defeq
\KK^{G\cross (X\times Y)}(C_0(Y,A), C_0(Y,B)).
$$
We let $\RKK^{G\cross X}(Y)$ be the category with the same objects as
$\KK^{G\cross X}$ and with morphisms as above.  That is, $\RKK^{G\cross X}(Y)$
identifies with the range of the functor
\begin{equation}  \label{eq:pY_induced}
  p_Y^*\colon \KK^{G\cross X}\to \KK^{G\cross (X\times Y)}
\end{equation}
induced by the projection map $X\times Y\to X$ and thus with a subcategory of
$\KK^{G\cross (X\times Y)}$.  (There is no reason to expect this subcategory
to be triangulated.)

Let $f\colon Y\to Y'$ be a continuous $G$\nbd{}equivariant map.  Since
$p_{Y'}\circ f=p_Y$, the functor $f^*\colon \KK^{G\cross (X\times
  Y')}\to\KK^{G\cross (X\times Y)}$ yields natural maps
\begin{equation}  \label{eq:RKKG_functor}
  f^*\colon \RKK^{G\cross X}(Y';A,B)\to\RKK^{G\cross X}(Y;A,B).
\end{equation}
For the constant map $p_Y\colon Y\to\pt$ this reproduces the functor $p_Y^*$
in~\eqref{eq:pY_induced}.  The maps in~\eqref{eq:RKKG_functor} turn $Y\mapsto
\RKK^{G\cross X} (Y;A,B)$ into a contravariant functor.

If~$S$ is a compact $G$\nbd{}space, then~\eqref{eq:map_induced_adjoint} yields
a natural isomorphism
$$
\RKK^{G\cross X} (Y\times S;A,B) \cong \KK^{G\cross X}(Y;A,C(S,B)).
$$
For $S=[0,1]$, we see that homotopy invariance of $\KK^G(A,B)$ in the
variable~$B$ implies homotopy invariance in the variable~$Y$; that is,
$f_1^\ast=f_2^\ast$ if $f_1, f_2\colon Y\to Y'$ are $G$\nbd{}equivariantly
homotopic.  Let $\EG$ be a second countable, locally compact universal proper
$G$\nbd{}space.  Then $X\times\EG$ is a universal proper $G\cross
X$\brd{}space.  The category $\RKK^{G\cross X}(\EG)$ and the functor
$p_\EG^*\colon \KK^{G\cross X}\to\RKK^{G\cross X}(\EG)$ do not depend on the
choice of $\EG$ because $\EG$ is unique up to homotopy.

\begin{theorem}  \label{the:derived_concrete}
  The functor $p_\EG^*\colon \KK^{G\cross X}\to\RKK^{G\cross X}(\EG)$ descends
  to an isomorphism of categories $\KK^{G\cross X}/\CC\cong\RKK^{G\cross
    X}(\EG)$.
  
  More explicitly, let $\pi\colon \tilde{A}\to A$ be a $\CI$\brd{}simplicial
  approximation of $A\in\KK^{G\cross X}$.  Then the indicated maps in the
  following commutative diagram are isomorphisms:
  $$
  \xymatrix{
    {\KK^{G\cross X}(\tilde{A},B)}
    \ar[d]^{\cong}_{p_\EG^\ast} &
    {\KK^{G\cross X}(A,B)} \ar[d]^{p_\EG^\ast}  \ar[l]_{\pi^*} \\
    {\RKK^{G\cross X}(\EG;\tilde{A},B)} &
    {\RKK^{G\cross X}(\EG;A,B).} \ar[l]^{\pi^*}_{\cong}
  }
  $$
\end{theorem}

\begin{proof}
  The first assertion follows from the second one and
  Proposition~\ref{pro:complementary}.\ref{pro::PN_adjoint}.  Hence we only
  have to prove that the two indicated maps are isomorphisms.  Consider
  $p_\EG^*$ first.  Fix~$B$.  Both $\KK^{G\cross X}(\blank,B)$ and
  $\RKK^{G\cross X}(\EG;\blank,B)$ are cohomological functors compatible with
  direct sums.  Thus the class of objects~$\tilde{A}$ for which the natural
  transformation $p_\EG^\ast$ between them is an isomorphism is localising.
  Hence we have an isomorphism for all $\tilde{A}\in\gen{\CI}$ once we have an
  isomorphism for $\tilde{A}\in\CI$.  This is what we shall prove.  Thus we
  let $\tilde{A}\defeq \Ind_H^G A'$ for some large compact subgroup~$H$ and
  some $A'\in\KK^{H\cross X}$.  Let $U\defeq U_H$ and let~$Y$ be a
  $G$\nbd{}space as above.  We use Lemma~\ref{lem:ind_res} and the
  compatibility of $\Ind_H^G$ with $p_Y^*$ to rewrite
  \begin{multline*}
    \RKK^{G\cross X}(Y;\Ind_H^G A',B)
    = \KK^{G\cross (X\times Y)} (\Ind_H^G p_Y^*A', p_Y^*B)
    \\ \cong \KK^{H\cross (X\times Y)} (p_Y^* \Ind_H^U A', p_Y^* B)
    = \RKK^{H\cross X}(Y; \Res_U^H \Ind_H^U A',\Res_G^H B);
  \end{multline*}
  we have dropped restriction functors from the notation except in the final
  result.  These isomorphisms are natural and especially compatible with the
  functoriality in~$Y$.  Since~$H$ is compact and $\EG$ is
  $H$\nbd{}equivariantly contractible, homotopy invariance implies that
  $p_\EG^\ast$ is an isomorphism
  $$
  \RKK^{H\cross X}(\pt; \Res_U^H \Ind_H^U A',\Res_G^H B)
  \congto \RKK^{H\cross X}(\EG; \Res_U^H \Ind_H^U A',\Res_G^H B).
  $$
  Hence $p_\EG^\ast\colon \RKK^{G\cross X}(\pt;\Ind_H^G A',B) \to
  \RKK^{G\cross X}(\EG;\Ind_H^G A',B)$ is an isomorphism as well.
  
  We claim that any weak equivalence $\pi\colon \tilde{A}\to A$ induces an
  isomorphism
  $$
  \RKK^{G\cross X}(\EG;\tilde{A},B)\cong\RKK^{G\cross X}(\EG;A,B).
  $$
  The proof of this claim will finish the proof of the theorem.  We remark
  that the usual definition of $\Ktop_*(G,A)$ is functorial for elements of
  $\RKK^G(\EG;A,A')$ in a rather obvious way.  Therefore, $f\in\KK^G(A,A')$
  induces an isomorphism on $\Ktop_*(G,A)$ once $p_\EG^*(f)$ is an isomorphism
  in $\RKK^G(\EG;A,A')$.  Thus the claim above implies that weak equivalences
  induce isomorphisms on $\Ktop_*(G,A)$.  We have already used this result
  of~\cite{Chabert-Echterhoff-Oyono:Going_down} in the proof of
  Theorem~\ref{the:Baum-Connes_assembly} above.
  
  The claim is equivalent to $\RKK^{G\cross X}(\EG;A,B)=0$ for $A\in\CC$ by
  Lemma~\ref{lem:weak_contractible_equivalence} because $\RKK^{G\cross
    X}(\EG;\blank,B)$ is a cohomological functor.  This condition for all~$B$
  is equivalent to $p_\EG^*(A)=0$ for $A\in\CC$, that is, $\CC\subseteq \ker
  p_\EG^*$.  Let~$\mathcal{S}$ be the class of proper $G$\nbd{}spaces~$Y$ for
  which $\CC\subseteq \ker p_Y^*$.  We shall use the following trivial
  observation.  If $Y\to Y'$ is a $G$\nbd{}equivariant map and
  $Y'\in\mathcal{S}$, then $Y\in\mathcal{S}$ as well because $p_Y^*$ factors
  through $p_{Y'}^*$.  Therefore, $\EG\in\mathcal{S}$ if and only if
  $Y\in\mathcal{S}$ for all proper $G$\nbd{}spaces~$Y$.  This is what we are
  going to prove.

  Let $H\subseteq G$ be a compact subgroup and let~$Y'$ be a locally compact
  $H$\nbd{}space.  Then we can form a $G$\nbd{}space $Y=G\times_H Y'$.  We
  call such $G$\nbd{}spaces \emph{compactly induced}.  The groupoid $G\cross
  (X\times Y)$ is Morita equivalent to $H\cross (X\times Y')$.  This yields an
  isomorphism
  $$
  \RKK^{G\cross X}(G\times_H Y';A,B)
  \cong \RKK^{H\cross X}(Y';\Res_G^H A,\Res_G^H B),
  $$
  (see \cite{Kasparov:Novikov}*{Theorem 3.6}) and hence factors~$p_Y^*$
  through $\Res_G^H$.  Thus $Y\in\mathcal{S}$, that is, $\mathcal{S}$ contains
  all compactly induced $G$\nbd{}spaces.
  
  Any locally compact proper $G$\nbd{}space is locally compactly induced, that
  is, can be covered by open $G$\nbd{}invariant subsets that are isomorphic to
  compactly induced spaces.  This result of Herbert
  Abels~\cite{Abels:Universal} is rediscovered
  in~\cite{Chabert-Echterhoff-Meyer:Deux}.  It ought to imply our claim by a
  Mayer-Vietoris argument.  However, the proof is somewhat delicate because it
  is unclear whether $\RKK^{G\cross X}(Y;A,B)$ as a functor of~$Y$ is
  excisive.
  
  We let~$\Cat_n$ be the class of proper $G$\nbd{}spaces that can be covered
  by at most~$n$ compactly induced, $G$\nbd{}invariant open subsets.
  Thus~$\Cat_1$ consists of the compactly induced $G$\nbd{}spaces.  We prove
  $\Cat_n\subseteq \mathcal{S}$ by induction on~$n$.  We already know
  $\Cat_1\subseteq\mathcal{S}$.  If $Y\in \Cat_n$, then $Y=Y_0\cup Y_1$ with
  open subsets $Y_0,Y_1$ such that $Y_0\in\Cat_1$ and $Y_1\in\Cat_{n-1}$.
  Hence $Y_0,Y_1\in\mathcal{S}$ by induction hypothesis.  Let $Y_\cap\defeq
  Y_0\cap Y_1$.  Then $Y_\cap\in\mathcal{S}$ as well because~$Y_\cap$ maps
  to~$Y_0$.  The idea of the following proof is to replace~$Y$ by a homotopy
  push-out~$Z$ of the diagram $Y_0\leftarrow Y_\cap \to Y_1$.  It is easy to
  see that $Z\in\mathcal{S}$.  Since there is a $G$\nbd{}equivariant map $Y\to
  Z$ this implies $Y\in\mathcal{S}$.
  
  In detail, let $\phi_0,\phi_1$ be a $G$\nbd{}invariant partition of unity
  subordinate to the covering $Y_0,Y_1$.  This can be constructed by working
  in $G\backslash Y$.  Let
  $$
  Z \defeq \bigl(Y_0\sqcup Y_1 \sqcup ([0,1]\times Y_\cap)\bigr)/\sim
  $$
  where we identify $Y_\cap\subseteq Y_j$ with $\{j\}\times Y_\cap
  \subseteq [0,1]\times Y_\cap$ for $j=0,1$.  We define a map $\phi_*\colon
  Y\to Z$ by $\phi_*(y)\defeq (\phi_1(y),y)\in [0,1]\times Y_\cap$ for $y\in
  Y_\cap$, $\phi_*(y)\defeq y\in Y_0$ for $y\in Y_0\setminus Y_\cap$ and
  $\phi_*(y)\defeq y\in Y_1$ for $y\in Y_1\setminus Y_\cap$.  Notice that this
  is a continuous, $G$\nbd{}equivariant map.  Thus $Z\in\mathcal{S}$ implies
  $Y\in\mathcal{S}$.
  
  A cycle for $\RKK^{G\cross X}(Z;A,A)$ is a triple $(f_0,f_1,f_\cap)$,
  where~$f_j$ is a cycle for $\RKK^{G\cross X} (Y_j;A,A)$ for $j=0,1$
  and~$f_\cap$ is a homotopy in $\RKK^{G\cross X}(Y_\cap;A,A)$ between
  $f_0|_{Y_\cap}$ and $f_1|_{Y_\cap}$.  Fix any such cycle.  The cycles $f_0$
  and~$f_1$ are homotopic to~$0$ because $Y_0,Y_1\in\mathcal{S}$ and
  $A\in\CC$.  This yields a homotopy between $(f_0,f_1,f_\cap)$ and
  $(0,0,f_\cap')$, where $f_\cap'$ is some cycle for $\RKK^{G\cross
    X}(Y_\cap\times [0,1];A,A)$ whose restrictions to $0$ and~$1$ vanish.
  Thus~$f_\cap'$ is equivalent to a cycle for
  $$
  \RKK^{G\cross X\times [0,1]}(Y_\cap;C([0,1],A),\Sigma A)
  \cong \RKK^{G\cross X}(Y_\cap;A,\Sigma A).
  $$
  Apply~\eqref{eq:map_induced_adjoint} to the coordinate projection
  $X\times Y_\cap \times [0,1]\to X\times Y_\cap$ to get this isomorphism.  We
  have $\RKK^{G\cross X}(Y_\cap;A,\Sigma A)=0$ because $Y_\cap\in\mathcal{S}$.
  Thus~$f_\cap'$ is homotopic to~$0$ and $\RKK^{G\cross X}(Z;A,A)=0$ for all
  $A\in\CC$, that is, $Z\in\mathcal{S}$.
  
  So far we have proven that $\Cat_n\subseteq\mathcal{S}$ for all $n\in\N$.
  Now let~$Y$ be an arbitrary proper $G$\nbd{}space.  Since~$Y$ is locally
  compactly induced, there is a locally finite covering by compactly induced
  $G$\nbd{}invariant open subsets $(U_j)_{j\in\N}$.  We let
  $Y_j=\bigcup_{k=0}^j U_j$ and write $Y=\bigcup Y_n$.  Thus $Y_n\in
  \Cat_n\subseteq \mathcal{S}$ for all $n\in\N$.  We use the following variant
  of the mapping telescope (compare Section~\ref{sec:homotopy_limits}):
  $$
  Z\defeq
  \{(y,t)\in Y\times\R_+\mid \text{$y\in Y_m$ whenever $t<m+1$}\}
  = \bigcup_{m\in\N} Y_m\times [m,m+1].
  $$
  This is a closed $G$\nbd{}invariant subset of $Y\times\R_+$.  There
  exists a partition of unity by $G$\nbd{}invariant functions subordinate to
  $(U_j)$ because $G\backslash Y$ is paracompact.  We use this to construct a
  $G$\nbd{}invariant function $\phi\colon Y\to\R_+$ with $\phi(y)\ge m$ for
  all $y\in Y_m\setminus Y_{m-1}$.  We get an embedding $Y\to Z$, $y\mapsto
  (y,\phi(y))$.  Thus $Y\in\mathcal{S}$ follows if $Z\in\mathcal{S}$.  The
  proof of $Z\in\mathcal{S}$ is analogous to the argument in the preceding
  paragraph.  Therefore, we are rather brief.
  
  A cycle for $\RKK^{G\cross X}(Z;A,A)$ is equivalent to sequences of cycles
  $(f_m)_{m\in\N}$ for $\RKK^{G\cross X}(Y_m;A,A)$ and homotopies
  $(H_m)_{m\in\N}$ between $f_m$ and $f_{m+1}|_{Y_m}$.  The assumption that
  $Y_m\in\mathcal{S}$ for all~$m$ allows us to find a homotopy between~$f_m$
  and~$0$ for all~$m$.  Thus the cycle described by the data
  $(f_m,H_m)_{m\in\N}$ is homotopic to a cycle $(0,H_m')_{m\in\N}$.
  Each~$H_m'$ is equivalent to a cycle for $\RKK^{G\cross X}(Y_m;A,\Sigma A)
  \cong 0$ and thus homotopic to~$0$.  Hence $\RKK^{G\cross X}(Z;A,A)=0$, that
  is, $Z\in\mathcal{S}$.  This finishes the proof.
\end{proof}

\begin{corollary}  \label{cor:CC_proper}
  An object of $\KK^{G\cross X}$ is weakly contractible if and only if its
  image in $\KK^{G\cross (X\times\EG)}$ is isomorphic to~$0$, if and only if
  its image in $\KK^{G\cross (X\times Y)}$ is isomorphic to~$0$ for all proper
  $G$\nbd{}spaces~$Y$.  A morphism in $\KK^{G\cross X}$ is a weak equivalence
  if and only if its image in $\KK^{G\cross (X\times\EG)}$ is invertible, if
  and only if its image in $\KK^{G\cross (X\times Y)}$ is invertible for all
  proper $G$\nbd{}spaces~$Y$.
\end{corollary}

\begin{proof}
  By the universal property of $\EG$ the map $Y\to\pt$ for any proper
  $G$\nbd{}space factors through $\EG$.  Hence assertions about $\RKK^G(\EG)$
  as in the statement of the corollary imply the corresponding assertions
  about $\RKK^G(Y)$ for all proper $G$\nbd{}spaces~$Y$.  An object is weakly
  contractible if and only if its image in the localisation vanishes and a
  morphism is a weak equivalence if and only if its image in the localisation
  is invertible.  Thus everything follows from
  Theorem~\ref{the:derived_concrete}.
\end{proof}

Recall that a $G\cross X$-$C^*$\brd{}algebra is called \emph{proper} if it is
a $G\cross (X\times Y)$\brd{}algebra for some proper $G$\nbd{}space~$Y$.

\begin{corollary}  \label{cor:proper_CI}
  All proper $G\cross X$-$C^*$\brd{}algebras belong to $\gen{\CI}$.
\end{corollary}

\begin{proof}
  Let~$A$ be a proper $G\cross X$-$C^*$\brd{}algebra.  Then~$A$ is a $G\cross
  (X\times\EG)$-$C^*$\brd{}algebra.  Let $\Dirac\colon \ADir\to C_0(X)$ be a
  Dirac morphism for $G\cross X$.  Since~$\Dirac$ is a weak equivalence,
  $p^*_\EG(\Dirac)$ is an invertible morphism in $\KK^{G\cross (X\times\EG)}$
  by Corollary~\ref{cor:CC_proper}.  Hence
  $$
  p^*_\EG(\Dirac)\otimes_{X\times\EG} \ID_A
  \in \KK^{G\cross (X\times\EG)}(p^*_\EG(\ADir)
  \otimes_{X\times\EG} A, C_0(X\times \EG)\otimes_{X\times\EG} A)
  $$
  is invertible.  If we forget the $\EG$\nbd{}structure, we still have an
  invertible element in $\KK^{G\cross X}$.  Equation~\eqref{eq:tensorX_tensor}
  implies $C_0(X\times \EG)\otimes_{X\times\EG} A\cong A$ and $p^*_\EG(\ADir)
  \otimes_{X\times\EG} A\cong \ADir\otimes_X A \in \gen{\CI}$.  These
  isomorphisms identify $p^*_\EG(\Dirac)\otimes_{X\times\EG} \ID_A$ with
  $\Dirac_*\in\KK^G(\ADir\otimes_X A,A)$.  Thus~$\Dirac$ is invertible and
  $A\in\gen{\CI}$.
\end{proof}

We do not know whether, conversely, any object in $\gen{\CI}$ is isomorphic in
$\KK^{G\cross X}$ to a proper $G$-$C^*$\brd{}algebra.  Since $A\in\gen{\CI}$
implies $A\cong \ADir\otimes_X A$, this holds if and only if the
source~$\ADir$ of the Dirac morphism for $G\cross X$ has this property.  Thus
the question is whether we can find a Dirac morphism whose source is proper.
This can be done for many groups.  For instance, if~$G$ is almost connected
with maximal compact subgroup~$K$, then the cotangent bundle $T^*(G/K)$ always
has a $K$\nbd{}equivariant spin structure, so that its Dirac operator is
defined.  It is indeed a Dirac morphism for~$G$ by results of Gennadi
Kasparov~\cite{Kasparov:Novikov}.  This is where our terminology comes from.
Generalising this construction to non-Hausdorff manifolds, one can construct
concrete Dirac morphisms of a similar sort for totally disconnected groups
with finite dimensional $\EG$ (see \cites{Kasparov-Skandalis:Buildings,
  Emerson-Meyer:Descent}).

\section{Dual Dirac morphisms}
\label{sec:dual_Dirac}

Let $\Sigma\AN\to \ADir \overset{\Dirac}\to C_0(X) \to\AN$ be a Dirac
triangle.

\begin{definition}  \label{def:dual_Dirac}
  We call $\eta\in\KK^{G\cross X}(C_0(X),\ADir)$ a \emph{dual Dirac morphism}
  for $G\cross X$ if $\eta\circ\Dirac=\ID_\ADir$.  The composition
  $\gamma\defeq\Dirac\eta\in\KK^{G\cross X}(C_0(X),C_0(X))$ is called a
  \emph{$\gamma$\nbd{}element} for $G\cross X$.
\end{definition}

Kasparov's Dirac dual Dirac method is the main tool for proving injectivity
and bijectivity of the Baum-Connes assembly map.  The following theorem shows
that a dual Dirac morphism in the above sense exists whenever the Dirac dual
Dirac method in the usual sense applies.  Our reformulation has the advantage
that the Dirac morphism is fixed, so that we only have to find one piece of
structure.  This is quite useful for analysing the existence of a dual Dirac
morphism (see~\cite{Emerson-Meyer:Descent}).

\begin{theorem}  \label{the:dual_Dirac_classical}
  Let~$A$ be a $\Ztwo$\nbd{}graded $G\cross X$-$C^*$\brd{}algebra,
  $\alpha\in\KK^{G\cross X}(A,C_0(X))$ and $\beta\in\KK^{G\cross
    X}(C_0(X),A)$.  If $\gamma\defeq \alpha\beta\in\KK^{G\cross
    X}(C_0(X),C_0(X))$ satisfies $p^*_\EG(\gamma)=1$ and~$A$ is proper, then
  there is a dual Dirac morphism for $G\cross X$.  Moreover, $\gamma$ is equal
  to the $\gamma$\nbd{}element.
  
  If, in addition, $A$ is trivially graded and $\beta\alpha=1$, then $\alpha$
  and~$\beta$ themselves are Dirac and dual Dirac morphisms for~$G$.
\end{theorem}

\begin{proof}
  Let $\Dirac\in\KK^{G\cross X}(\ADir,C_0(X))$ be a Dirac morphism.  Even
  if~$A$ is graded, the same argument as in the proof of
  Corollary~\ref{cor:proper_CI} shows that $\Dirac_*\in\KK^G(\ADir\otimes_X
  A,A)$ is invertible---provided~$A$ is proper.  We claim that the composite
  morphism
  $$
  \eta\colon C_0(X)
  \overset{\beta}\longrightarrow A
  \overset{\Dirac^{-1}_*}\longrightarrow \ADir\otimes_X A
  \overset{\alpha_*}\longrightarrow \ADir
  $$
  is a dual Dirac morphism, that is, $\eta\circ\Dirac=1_\ADir$.  We have
  $\Dirac\circ\eta=\beta\circ\alpha=\gamma$ because exterior products are
  graded commutative.  Since~$\Dirac$ is a weak equivalence, $p_\EG^*(\Dirac)$
  is invertible.  Since $1=p_\EG^*(\gamma)=p_\EG^*(\eta\Dirac)$, we get
  $p_\EG^*(\eta) = p_\EG^*(\Dirac)^{-1}$.  Therefore,
  $p_\EG^*(\eta\Dirac)=1=p_\EG^*(1)$.  The map
  $$
  p_\EG^*\colon \KK^{G\cross X}(\ADir,\ADir)
  \to \RKK^{G\cross X}(\EG;\ADir,\ADir)
  $$
  is bijective by Theorem~\ref{the:derived_concrete} because
  $\ADir\in\gen{\CI}$.  Hence $\eta\Dirac=1$.
  
  If~$A$ is trivially graded, then $A\in\gen{\CI}$ by
  Corollary~\ref{cor:proper_CI}.  The morphisms $\alpha$ and~$\beta$ are weak
  equivalences because $\beta\alpha=1$ and $\alpha\beta=\gamma$ are.  This
  implies that~$\alpha$ is a Dirac morphism and that~$\beta$ is a dual Dirac
  morphism.
\end{proof}

\begin{theorem}  \label{the:dual_Dirac_splitting}
  The following assertions are equivalent:
  \begin{enumerate}[\ref{the:dual_Dirac_splitting}.1.]
  \item there is a dual Dirac morphism ($\eta\in\KK^{G\cross X}(C_0(X),\ADir)$
    with $\eta\Dirac=\ID_{\ADir}$);

  \item $\KK^{G\cross X}_*(\AN,\ADir)=0$ (for all $*\in\Z$);
    
  \item the natural map $p_\EG^*\colon \KK^{G\cross X}_*
    (C_0(X),\ADir)\to\RKK^{G\cross X}_*(\EG;C_0(X),\ADir)$ is an
    isomorphism (for all $*\in\Z$);

  \item $\KK^{G\cross X}_*(A,B)=0$ for all $A\in\CC$, $B\in\gen{\CI}$;
    
  \item the natural map $\KK^{G\cross X}_*(A,B)\to\RKK^{G\cross
      X}_*(\EG;A,B)$ is an isomorphism for all $A\in\KK^{G\cross X}$,
    $B\in\gen{\CI}$;
    
  \item there is an equivalence of triangulated categories $\KK^{G\cross
      X}\cong\gen{\CI}\times\CC$.

  \end{enumerate}

  Suppose these equivalent conditions to be satisfied and let
  $$
  \gamma\defeq\Dirac\circ\eta\in\KK^{G\cross X}(C_0(X),C_0(X)).
  $$
  Then $\gamma_A\defeq \gamma\otimes_X A\in\KK^{G\cross X}(A,A)$ is an
  idempotent for all $A\in\KK^{G\cross X}$.  We have $\gamma_A=0$ if and only
  if $A\in\CC$ and $\gamma_A=\ID$ if and only if $A\in\gen{\CI}$.
\end{theorem}

\begin{proof}
  We often use the isomorphism $\RKK^{G\cross X}_*(\EG;A,B)\cong \KK^{G\cross
    X}_*(\ADir\otimes_X A,B)$ proven in Theorem~\ref{the:derived_concrete}.  A
  long exact sequence argument shows that \ref{the:dual_Dirac_splitting}.2 and
  \ref{the:dual_Dirac_splitting}.3 are equivalent.  Conditions
  \ref{the:dual_Dirac_splitting}.4 and \ref{the:dual_Dirac_splitting}.5 are
  two ways of expressing that objects of $\gen{\CI}$ are $\CC$\brd{}injective
  and hence equivalent.  The implications
  \ref{the:dual_Dirac_splitting}.6$\Longrightarrow
  $\ref{the:dual_Dirac_splitting}.4$\Longrightarrow
  $\ref{the:dual_Dirac_splitting}.2 and
  \ref{the:dual_Dirac_splitting}.3$\Longrightarrow
  $\ref{the:dual_Dirac_splitting}.1 are trivial.  It remains to prove that
  \ref{the:dual_Dirac_splitting}.1 implies \ref{the:dual_Dirac_splitting}.6.
  Along the way we show the additional assertions about~$\gamma$ (and part of
  the following corollary).
  
  Since $\eta\Dirac=\ID_\ADir$, the map $\Sigma\AN\to\ADir$ in the Dirac
  triangle vanishes.  Hence Lemma~\ref{lem:triangle_lemma} yields an
  isomorphism $C_0(X)\cong \ADir\oplus\AN$ such that the maps $\ADir\to
  C_0(X)\to\AN$ become the obvious ones.  Any two choices for~$\eta$ differ by
  a morphism $\AN\to\ADir$.  Therefore, \ref{the:dual_Dirac_splitting}.2
  implies that~$\eta$ is unique.  We cannot use this so far because we still
  have to prove that \ref{the:dual_Dirac_splitting}.2 follows from
  \ref{the:dual_Dirac_splitting}.1.  We may, however, choose~$\eta$ such that
  $\gamma=\Dirac\eta$ is the orthogonal projection onto~$\ADir$ that vanishes
  on~$\AN$.  We get a direct sum decomposition (in $\KK$)
  $$
  A \cong C_0(X) \otimes_X A \cong \ADir\otimes_X A \oplus \AN\otimes_X A
  $$
  such that $\Dirac\otimes_X\ID_A$ is the inclusion of the first summand
  and~$\gamma_A$ is the orthogonal projection onto $\ADir\otimes_X A$.
  Theorem~\ref{the:KKG_decomposition} yields $\gamma_A=1$ if and only if
  $A\in\gen{\CI}$, and $\gamma_A=0$ if and only if $A\in\CC$.  Since
  $$
  \gamma_B\circ f = \gamma\otimes_X f = f\circ\gamma_A
  $$
  for all $f\in\KK^{G\cross X}(A,B)$, there are no non-zero morphisms
  between $\CC$ and $\gen{\CI}$.  The above decomposition of~$A$ respects
  suspensions and exact triangles because the tensor product functors
  $\ADir\otimes_X\blank$ and $\AN\otimes_X\blank$ are triangulated.  Hence we
  get an equivalence of triangulated categories
  $\gen{\CI}\times\CC\cong\KK^G$.
\end{proof}

\begin{corollary}  \label{cor:dual_Dirac}
  Fix a Dirac morphism $\Dirac\in\KK^{G\cross X}(\ADir,C_0(X))$.  Then the
  dual Dirac morphism and the $\gamma$\nbd{}element are unique if they exist.

  A morphism $\eta\in\KK^{G\cross X}(C_0(X),\ADir)$ is a dual Dirac morphism
  if and only if $p^*_\EG(\eta)$ is inverse to $p^*_\EG(\Dirac)$ if and only
  if $p^*_\EG(\Dirac\eta)=1$.
\end{corollary}

\begin{proof}
  We have already shown the uniqueness of~$\eta$ and hence of~$\gamma$ during
  the proof of Theorem~\ref{the:dual_Dirac_splitting}.  The map $p^*_\EG\colon
  \KK^{G\cross X}_*(\ADir,\ADir)\to \RKK^{G\cross X}_*(\EG;\ADir,\ADir)$ is an
  isomorphism by Theorem~\ref{the:derived_concrete}.  Hence $\eta\Dirac=\ID$
  if and only if $p^*_\EG(\eta\Dirac)=\ID$.  Since~$\Dirac$ is a weak
  equivalence, $p^*_\EG(\Dirac)$ is an isomorphism.  Hence there is no
  difference between left, right and two-sided inverses for $p_\EG^*(\Dirac)$.
\end{proof}

Suppose now that a dual Dirac morphism exists.  It induces a canonical section
for the map $\ADir\otimes_X A\to A$.  Hence the natural transformation $\Left
F(A)\to F(A)$ for a covariant functor~$F$ is naturally split injective.
Similarly, the natural transformation $F(A)\to F(\ADir\otimes_X A)$ is
naturally split surjective for a contravariant functor~$F$.

It is clear from Theorem~\ref{the:dual_Dirac_splitting} that $\gamma=1$ if and
only if $\CC=0$.  In this case, $\Left F(A)\cong F(A)$ for any functor on
$\KK^{G\cross X}$, that is, any functor~$F$ satisfies the analogue of the
Baum-Connes conjecture.  Nigel Higson and Gennadi Kasparov show
in~\cite{Higson-Kasparov:Amenable} that all groups with the Haagerup property
and in particular all amenable groups have a dual Dirac element and satisfy
$\gamma=1$.  Jean-Louis Tu generalises their argument to groupoids that
satisfy an analogue of the Haagerup property in~\cite{Tu:Amenable}.  In
particular, this applies to the special groupoids $G\cross X$.  We get:

\begin{theorem}  \label{the:Haagerup}
  Suppose that the groupoid $G\cross X$ is amenable or, more generally, acts
  continuously and isometrically on a continuous field of affine Euclidean
  spaces over~$X$.  Then weakly contractible objects of $\KK^{G\cross X}$ are
  already isomorphic to~$0$ and weak equivalences are isomorphisms.  The
  assembly map is an isomorphism for any functor defined on $\KK^{G\cross X}$.
\end{theorem}

\subsection{Approximate dual Dirac morphisms}
\label{sec:approx_dual_Dirac}

In some cases of interest, for instance, for groups acting on bolic spaces,
one cannot construct an actual dual Dirac morphism but only approximations to
one.

\begin{definition}  \label{def:approximate_dual_Dirac}
  Suppose that for each $G$\nbd{}compact proper $G$\nbd{}space~$Y$ there is
  $\eta_Y\in\KK^G(\pt,\ADir)$ such that $p_Y^*(\Dirac\circ\eta_Y)=1\in
  \RKK^G(Y;\pt,\pt)$.  Then we call the family $(\eta_Y)$ an \emph{approximate
    dual Dirac morphism} for~$G$.  We also let
  $\gamma_Y\defeq\Dirac\circ\eta_Y$.
\end{definition}

\begin{lemma}  \label{lem:get_approximate_dual_Dirac}
  Suppose that for each $G$\nbd{}compact proper $G$\nbd{}space~$Y$ there are a
  possibly $\Ztwo$\nbd{}graded, proper $G$\nbd{}$C^*$\brd{}algebra~$A_Y$ and
  $\alpha_Y\in\KK^G(A_Y,\pt)$, $\beta_Y\in\KK^G(\pt,A_Y)$ such that
  $p_Y^*(\alpha_Y\circ\beta_Y)=1\in \RKK^G(Y;\pt,\pt)$.  Then~$G$ has an
  approximate dual Dirac morphism with $\gamma_Y=\alpha_Y\beta_Y$.
\end{lemma}

\begin{proof}
  Proceed as in the proof of Theorem~\ref{the:dual_Dirac_classical}.
\end{proof}

The situation of Lemma~\ref{lem:get_approximate_dual_Dirac} occurs
in~\cite{Kasparov-Skandalis:Bolic}.  It follows that a discrete group~$G$ has
an approximate dual Dirac morphism if it acts properly and by isometries on a
weakly bolic, weakly geodesic metric space.  Clearly, $G$ has an approximate
dual Dirac morphism once it has a dual Dirac morphism.  The converse holds
if~$G$ does not have too many compact subgroups:

\begin{proposition}  \label{pro:approximate_dual_Dirac}
  Suppose that there exist finitely many compact subgroups of~$G$ such that
  any compact subgroup is subconjugate to one of them.  If~$G$ has an
  approximate dual Dirac morphism, then it already has a dual Dirac morphism.
\end{proposition}

\begin{proof}
  Let $\Dirac\in\KK^G(\ADir,\pt)$ be a Dirac morphism for~$G$.  Let~$S$ be a
  finite set of compact subgroups such that any other compact subgroup is
  subconjugate to one of them.  Let~$Y$ be the disjoint union of the spaces
  $G/H$ for $H\in S$.  By hypothesis, there is $\eta_Y\in\KK^G(\pt,\ADir)$
  such that $\gamma_Y\defeq\Dirac\circ\eta_Y$ satisfies $p_Y^*(\gamma_Y)=1$.
  This means that $p_{G/H}^*(\gamma_Y)=1$ for all $H\in S$.
  By~\eqref{eq:resind_as_induced}, this is equivalent to
  $\Res_G^H(\gamma_Y)=1$ for all $H\in S$.  By hypothesis, any compact
  subgroup of~$G$ is subconjugate to one in~$S$.  Thus~$\gamma_Y$ is a weak
  equivalence.  Since~$\Dirac$ is a weak equivalence as well, it follows
  that~$\eta_Y$ is a weak equivalence.  Hence the composition
  $\eta_Y\circ\Dirac\in\KK^G(\ADir,\ADir)$ is a weak equivalence.  Since
  $\ADir\in\gen{\CI}$, it is projective with respect to weak equivalences by
  Proposition~\ref{pro:CC_CI_orthogonal}.  Hence $\eta_Y\circ\Dirac$ is
  invertible; $\eta\defeq (\eta_Y\circ\Dirac)^{-1}\eta_Y\in\KK^G(\pt,\ADir)$
  is the desired dual Dirac morphism for~$G$.
\end{proof}

It is unclear whether the condition on compact subgroups in
Proposition~\ref{pro:approximate_dual_Dirac} can be removed.  Our next goal is
a weakening of Theorem~\ref{the:dual_Dirac_splitting}.5, which still holds
if~$G$ has an approximate dual Dirac morphism and which is used
in~\cite{Emerson-Meyer:Descent}.

\begin{lemma}  \label{lem:rotation_trick}
  Let $\Dirac\in\KK^G(\ADir,\pt)$ be a Dirac morphism and let
  $\alpha\in\KK^G(\pt,\ADir)$.  Define $\beta\defeq
  \Dirac\circ\alpha\in\KK^G(\pt,\pt)$ and $\beta_A\defeq \beta\otimes \ID_A
  \in\KK^G(A,A)$ for all $A\in\KK^G$.  Then $\alpha\circ\Dirac=\beta_\ADir$.
  For $A,B\in\KK^G$, the composites
  \begin{align*}
    \alpha^*\Dirac^*\colon
    \KK^G(A,B) \to \KK^G(\ADir\otimes A,B) \to \KK^G(A,B),
    \\
    \Dirac^*\alpha^*\colon
    \KK^G(\ADir\otimes A,B) \to \KK^G(A,B) \to \KK^G(\ADir\otimes A,B),
  \end{align*}
  are both given by $f\mapsto \beta_B \circ f$.
\end{lemma}

\begin{proof}
  Since~$\Dirac$ is a weak equivalence, the map
  $$
  \Dirac_*\colon \KK^G(\ADir,\ADir)\to\KK^G(\ADir,\pt)
  $$
  is an isomorphism.  It maps both $\alpha\circ\Dirac$ and $\beta_\ADir$ to
  $\beta\circ\Dirac=\Dirac\otimes\beta$.  Hence
  $\alpha\circ\Dirac=\beta_\ADir$.  The second assertion now follows from
  $\beta_B\circ f=\beta\otimes f=f\circ\beta_{A'}$ for all $A',B\in\KK^G$,
  $f\in\KK^G(A',B)$ (applied to $A'=A$ and $A'=\ADir\otimes A$).
\end{proof}

\begin{lemma}  \label{lem:gamma_one}
  Let $\beta\in\KK^G(\pt,\pt)$, let~$Y$ be a locally compact
  $G$\nbd{}space and let~$A$ be a $G\cross Y$-$C^*$\brd{}algebra.  If
  $p_Y^*(\beta)=1$, then $\beta_A=1$.
\end{lemma}

\begin{proof}
  We have a natural isomorphism $B\otimes A\cong p_Y^*(B)\otimes_Y A$ for all
  $A$ and~$B$.  Hence $\beta_A \defeq \beta\otimes\ID_A = p_Y^*(\beta)
  \otimes_Y \ID_A = 1 \otimes_Y \ID_A = 1$.
\end{proof}

For a finite set of compact subgroups~$S$, let $\CI(S)\subseteq\CI$ be the
class of $G$\nbd{}$C^*$\brd{}algebras that are $\KK^G$\brd{}equivalent to
$\Ind_H^G(A)$ for some $H\in S$ and some $A\in\KK^H$.  Let $\gen{\CI(S)}$ be
the localising subcategory generated by $\CI(S)$.  These subcategories form a
directed set of localising subcategories.  Let $\Prfin$ be their union, that
is, $A\in\Prfin$ if and only if $A\in\CI(S)$ for some \emph{finite} set of
compact subgroups~$S$.  This is a thick, triangulated subcategory of $\KK^G$,
but it need not be localising: it is only closed under countable direct sums
if all summands lie in the same category $\CI(S)$ for some~$S$.  The
hypothesis of Proposition~\ref{pro:approximate_dual_Dirac} ensures that
$\Prfin=\CI(S)$ for some~$S$, so that $\Prfin$ is localising as well.  Thus
$\Prfin=\gen{\CI}$ in this case.  In general,
$\CI\subseteq\Prfin\subseteq\gen{\CI}$ and both containments may be strict.

\begin{proposition}  \label{pro:approximate_dual_Dirac_II}
  If~$G$ has an approximate dual Dirac morphism, then the map
  \begin{equation}   \label{eq:approximate_dual_Dirac_II}
    p_\EG^*\colon \KK^G_*(A,B)\to\RKK^G_*(\EG;A,B)
  \end{equation}
  is an isomorphism for all $B\in\Prfin$, $A\in\KK^G$.
\end{proposition}

\begin{proof}
  Fix $B\in\Prfin$ and let~$S$ be a finite set of compact subgroups of~$G$
  such that $B\in\gen{\CI(S)}$.  Let~$Y$ be the disjoint union of the spaces
  $G/H$ for $H\in S$.  This is a $G$\nbd{}compact proper $G$\nbd{}space.
  Since~$G$ has an approximate dual Dirac morphism, there is
  $\eta_Y\in\KK^G(\pt,\ADir)$ such that $\gamma\defeq\Dirac\eta_Y$ satisfies
  $p_Y^*(\gamma)=1$.  This yields $\gamma_{B'}=1$ in $\KK^G(B',B')$ for
  $B'\in\CI(S)$ by Lemma~\ref{lem:gamma_one}.  In particular, $\gamma_{B'}$ is
  invertible if $B'\in\CI(S)$.  The class of $B'\in\KK^G$ for
  which~$\gamma_{B'}$ is invertible is localising by the Five Lemma and the
  functoriality of direct sums.  Hence~$\gamma_{B'}$ is invertible for all
  $B'\in\gen{\CI(S)}$ and, especially, for our chosen~$B$.  By
  Lemma~\ref{lem:rotation_trick}, $\Dirac^*\colon
  \KK^G(A,B)\to\KK^G(\ADir\otimes A,B)$ is invertible because both
  $\Dirac^*\eta_Y^*$ and $\eta_Y^*\Dirac^*$ are equal to invertible maps of
  the form $(\gamma_B)_*$.  Theorem~\ref{the:derived_concrete} allows us to
  replace~$\Dirac^*$ by the map $p_\EG^*$
  in~\eqref{eq:approximate_dual_Dirac_II}.
\end{proof}

\section{The strong Baum-Connes conjecture}
\label{sec:strong_BC}

\begin{definition}  \label{def:strong_BC}
  We say that~$G$ satisfies the \emph{strong Baum-Connes conjecture with
    coefficients $A\in\KK^G$} if the assembly map $G\Lrcross A\to G\rcross A$
  is a $\KK$\brd{}equivalence.
\end{definition}

The strong Baum-Connes conjecture implies that the assembly map is an
isomorphism for any functor defined on $\KK$.  In particular, this covers
$\K$\nbd{}theory, $\K$\nbd{}homology and local cyclic homology and cohomology
of the reduced crossed product.

Suppose that~$G$ has a dual Dirac morphism and resulting
$\gamma$\nbd{}element~$\gamma$.  Applying descent, we get $G\rcross
\gamma_A\in\KK(G\rcross A,G\rcross A)$.  The strong Baum-Connes conjecture
amounts to the assertion that $G\rcross \gamma_A=1$.  This is known to be
false for quite some time if $A=\pt$ and~$G$ is a discrete subgroup of finite
covolume in $\mathrm{Sp}(n,1)$ (\cite{Skandalis:Knuclear}).

For general~$G$, the Baum-Connes conjecture with coefficients~$A$ holds if and
only if ${G\Obscross A}$ is $\K$\nbd{}contractible as in
Definition~\ref{def:K_contractible}, whereas the \emph{strong} Baum-Connes
conjecture with coefficients~$A$ holds if and only if $G\Obscross A\cong0$ in
$\KK$.  These two assertions are equivalent if $G\Obscross A$ belongs to the
bootstrap category $\gen{\pt}$.  A sufficient condition for $G\Obscross A
\in\gen{\pt}$ is that both $G\Lrcross A$ and $G\rcross A$ belong to the
bootstrap category.

Now we use the notion of smooth compact subgroup introduced in
Section~\ref{sec:compact_subgroups}.  If~$G$ is discrete, any finite subgroup
of~$G$ is smooth.  Let $\CI_1\subseteq\CI$ be the set of all
$G$\nbd{}$C^*$\brd{}algebras of the form $C_0(G/H)$ for smooth, compact
subgroups $H\subseteq G$.  This is a variant of the subcategory
$\CIZ\subseteq\CI$ introduced in Section~\ref{sec:construct_Dirac}.

The following lemma is motivated by work of Jérôme Chabert and Siegfried
Echterhoff (see, for instance,
\cite{Chabert-Echterhoff-Oyono:Going_down}*{Lemma 4.20}).

\begin{proposition}  \label{pro:RC_decompose}
  The localising category $\gen{\CI_1}$ generated by $\CI_1$ contains
  $\gen{\CIZ}$ and hence also contains the source of the Dirac
  morphism.
\end{proposition}

\begin{proof}
  Our existence proof for Dirac morphisms shows that $\ADir\in\gen{\CIZ}$.  If
  the generators $\RC_H$ defined in~\eqref{eq:RC_concrete} belong to
  $\gen{\CI_1}$, then $\gen{\CIZ}\subseteq\gen{\CI_1}$, and we are done.  Let
  $H\subseteq G$ be a large compact subgroup, $U\defeq U_H$ and $V\defeq U/H$.
  Recall that $\RC_H=\Ind_H^G C_0(V^7)$.  Since~$U$ is almost connected, there
  is a compact normal subgroup $N\subseteq U$ such that $U/N$ is a Lie group.
  By maximality, $H=NH\supseteq N$.  The quotient $H/N$ is a compact Lie
  group.  It acts linearly on the $\R$\nbd{}vector space~$V^7$.  One can show
  that~$V^7$ is an $H/N$\nbd{}CW\brd{}complex; this is a special case
  of~\cite{Illman:Equivariant_Triangulations}.  Hence $C_0(V^7)$ belongs to
  the localising subcategory of $\KK^{H/N}$ generated by $C_0(H/K)$ with
  $N\subseteq K\subseteq H$.  Hence~$\RC_H$ belongs to the localising
  subcategory of $\KK^G$ generated by $\Ind_H^G C_0(H/K) \cong C_0(G/K)$ for
  such~$K$.  Since~$N$ is a strongly smooth compact subgroup of~$G$ contained
  in each~$K$, the assertion follows.
\end{proof}

\begin{theorem}  \label{the:Lrcross_gen}
  For any $A\in\KK^G$, the $C^*$\nbd{}algebra $G\Lrcross A$ belongs to the
  localising subcategory of $\KK$ generated by $H\cross A$ for smooth compact
  subgroups $H\subseteq G$.  In particular, $G\Lrcross A\in\gen{\pt}$ once
  $H\cross A\in\gen{\pt}$ for all smooth compact subgroups~$H$.
  
  If $H\cross A \cong 0$ in $\KK$ for all smooth compact subgroups~$H$, then
  $G\Lrcross A\cong0$ as well.  If $f\in\KK^G(A,B)$ induces
  $\KK$\brd{}equivalences $H\cross A\cong H\cross B$ for all smooth compact
  subgroups~$H$, then it induces a $\KK$\brd{}equivalence $G\Lrcross A\cong
  G\Lrcross B$.
  
  If $H\cross A$ is $\K$\nbd{}contractible for all smooth compact
  subgroups~$H$, so is $G\Lrcross A$.  If $f\in\KK^G(A,B)$ induces a
  $\K$\nbd{}equivalence $H\cross A\to H\cross B$ for all smooth compact
  subgroups~$H$, then it induces a $\K$\nbd{}equivalence $G\Lrcross A\to
  G\Lrcross B$.
\end{theorem}

\begin{proof}
  Proposition~\ref{pro:RC_decompose} implies that $G\Lrcross A\cong G\rcross
  (\ADir\otimes A)$ belongs to the localising subcategory of $\KK$ generated
  by $G\rcross C_0(G/H,A)$ for smooth compact subgroups $H\subseteq G$.
  Equation~\eqref{eq:Green} yields $G\rcross C_0(G/H,A)\sim_M H\cross A$.
  This implies the criteria for $G\Lrcross A$ to be in $\gen{\pt}$, to be
  $\KK$\brd{}contractible and to be $\K$\nbd{}contractible because all these
  conditions define localising subcategories of $\KK$.  The assertions about
  morphisms follow if we replace~$f$ by its mapping cone.
\end{proof}

The following corollary is originally due to Jean-Louis Tu
(\cite{Tu:Amenable}).  It applies to amenable groups by
Theorem~\ref{the:Haagerup}.

\begin{corollary}  \label{cor:crossed_bootstrap}
  Let~$G$ be a locally compact group, let~$X$ be a $G$\nbd{}space, and let
  $A\in\KK^{G\cross X}$.  Suppose that $G\cross X$ has a dual Dirac morphism
  with $\gamma=1$ or, more generally, $G\rcross\gamma_A=1\in\KK(G\rcross
  A,G\rcross A)$.  If $H\cross A\in\gen{\pt}$ for all smooth compact
  subgroups~$H$, then $G\rcross A\in\gen{\pt}$.
\end{corollary}

\begin{proof}
  If $\gamma=1\in\KK^{G\cross X}(C_0(X),C_0(X))$, then $(G\cross X)\rcross
  \gamma_A=1$.  This implies $(G\cross X)\Lrcross A\cong (G\cross X)\rcross A$
  in $\KK$.  Now use Lemma~\ref{lem:forget_X} to get rid of the space~$X$ and
  apply Theorem~\ref{the:Lrcross_gen}.
\end{proof}

Theorem~\ref{the:Lrcross_gen} describes other interesting localising
subcategories of $\KK^G$ on which $\Ktop_*(G,\blank)$ vanishes.  Hence it
gives a variant of the rigidity formulation of the Baum-Connes conjecture.
Namely, $G$ satisfies the Baum-Connes conjecture with arbitrary coefficients
if and only if $\K_*(G\rcross A)\cong0$ whenever $A\in\KK^G$ satisfies
$\K_*(H\rcross A)\cong0$ for all smooth compact subgroups $H\subseteq G$.

\begin{proposition}  \label{pro:Lrcross_genpt}
  If the $G$-$C^*$\brd{}algebra~$A$ is commutative (or just type~$I$), then
  $G\Lrcross A\in\gen{\pt}$.  In particular, $G\Lrcross\pt\in\gen{\pt}$.
  Suppose $G\Lrcross A\in\gen{\pt}$ (for instance, $A=\pt$).  Then the strong
  Baum-Connes conjecture with coefficients~$A$ holds if and only if $G\rcross
  A\in\gen{\pt}$ and the usual Baum-Connes conjecture with coefficients~$A$
  holds.
\end{proposition}

\begin{proof}
  If~$A$ is a type~I $C^*$-algebra and~$H$ is compact, then $H\cross A$ is a
  type~I $C^*$\nbd{}algebra as well (this follows easily from
  \cite{Takesaki:Covariant}*{Theorem 6.1}).  Therefore, it belongs to
  $\gen{\pt}$ (see \cite{Blackadar:Book}*{22.3.5}).  Thus $G\Lrcross
  A\in\gen{\pt}$ by Theorem~\ref{the:Lrcross_gen}.  The strong Baum-Connes
  conjecture is stronger than the Baum-Connes conjecture and implies that
  $G\rcross A\in\gen{\pt}$ once $G\Lrcross A\in\gen{\pt}$.  The converse also
  holds because a $\K$\nbd{}equivalence between objects of $\gen{\pt}$ is
  already a $\KK$\brd{}equivalence.
\end{proof}

Therefore, if we already know that $\Cred(G)\in\gen{\pt}$, then the strong and
the usual Baum-Connes conjecture with trivial coefficients are equivalent.
Jérôme Chabert, Siegfried Echterhoff and Hervé Oyono-Oyono show
in~\cite{Chabert-Echterhoff-Oyono:Going_down} that $\Cred(G)\in\gen{\pt}$
if~$G$ is almost connected or a linear algebraic group over the $p$\nbd{}adic
numbers or over the adele ring of a number field.  The Baum-Connes conjecture
with trivial coefficients for these groups is also known, see
\cites{Chabert-Echterhoff-Nest:Connes_Kasparov,
  Chabert-Echterhoff-Oyono:Going_down}.  Hence these groups satisfy the strong
Baum-Connes conjecture with trivial coefficients.

\section{Permanence properties of the (strong) Baum-Connes conjecture}
\label{sec:permanence}

Let $\Tri$ and~$\Tri'$ be triangulated categories, let $F\colon \Tri\to\Tri'$
be a triangulated functor, and let $(\Null,\Proj)$ and~$(\Null',\Proj')$ be
complementary pairs of localising subcategories in $\Tri$ and~$\Tri'$,
respectively.  Suppose $F(\Proj)\subseteq\Proj'$.  Then $\Left (F'\circ F)
\cong \Left F'\circ\Left F$ up to isomorphism for any covariant functor~$F'$
defined on~$\Tri'$.  This trivial observation has lots of applications.  When
applied to restriction and induction functors, partial crossed product
functors and the complexification functor, we get permanence properties of the
(strong) Baum-Connes conjecture.  We remark that Lemma~\ref{lem:forget_X} is
another such result that logically belongs into this section.

\subsection{Restriction and induction}
\label{sec:restrict_induct_permanence}

\begin{proposition}  \label{pro:res_ind_exact}
  Let $H\subseteq G$ be a closed subgroup.  The functors
  $$
  \Res_G^H\colon \KK^{G\cross X}\to\KK^{H\cross X}
  \quad\text{and}\quad
  \Ind_H^G\colon \KK^{H\cross X}\to\KK^{G\cross X}
  $$
  preserve weak contractibility and weak equivalences and map $\gen{\CI}$
  to $\gen{\CI}$.  Therefore, $\Res_G^H$ maps a Dirac triangle for $G\cross X$
  to a Dirac triangle for $H\cross X$ and $\Ind_H^G$ maps a Dirac triangle for
  $H\cross X$ to a Dirac triangle for $G\cross X$.
\end{proposition}

\begin{proof}
  Restriction and induction in stages yield $\Res_G^H(\CC)\subseteq\CC$ and
  $\Ind_H^G(\CI)\subseteq\CI$ and hence
  $\Ind_H^G(\gen{\CI})\subseteq\gen{\CI}$.  To prove
  $\Res_G^H(\gen{\CI})\subseteq\gen{\CI}$, it suffices to show
  $\Res_G^H(\CI)\subseteq\gen{\CI}$ because $\Res_G^H$ is triangulated and
  commutes with direct sums.  It is clear that $\Res_G^H$ maps compactly
  induced $G$-$C^*$\brd{}algebras to proper $H$\nbd{}$C^*$\brd{}algebras.
  This implies the assertion by Corollary~\ref{cor:proper_CI}.  As a
  consequence, $\Res_G^H$ maps a Dirac triangle for $G\cross X$ to one for
  $H\cross X$.
  
  Next we prove that $\Ind_H^G(\CC)\subseteq\CC$.  Let $\Dirac\in\KK^{G\cross
    X}(\ADir,C_0(X))$ be a Dirac morphism for $G\cross X$.  We have just seen
  that $\Res_G^H\Dirac$ is a Dirac morphism for $H\cross X$.  Let
  $A\in\KK^{H\cross X}$.  Equation~\eqref{eq:ind_tensorX} yields
  $$
  \ADir\otimes_X \Ind_H^G A \approx \Ind_H^G(\Res_G^H \ADir\otimes_X A).
  $$
  By Theorem~\ref{the:KKG_decomposition}, $\Ind_H^G A\in\CC$ is equivalent
  to $\ADir\otimes_X \Ind_H^G A\cong0$ and $A\in\CC$ is equivalent to
  $\Res_G^H \ADir\otimes_X A\cong 0$.  Thus $\Ind_H^G(\CC)\subseteq\CC$.  As a
  consequence, $\Ind_H^G$ maps a Dirac triangle for $H\cross X$ to one for
  $G\cross X$.
\end{proof}

It follows immediately from Proposition~\ref{pro:res_ind_exact} that
\begin{alignat*}{2}
  \Left (F\circ\Ind_H^G) &\cong (\Left F)\circ\Ind_H^G,
  &\qquad
  \Obs (F\circ\Ind_H^G) &\cong (\Obs F)\circ\Ind_H^G,
  \\
  \Left (F\circ\Res_G^H) &\cong (\Left F)\circ\Res_G^H,
  &\qquad
  \Obs (F\circ\Res_G^H) &\cong (\Obs F)\circ\Res_G^H.
\end{alignat*}
Since $G\rcross \Ind_H^G A\sim_M H\rcross A$ by~\eqref{eq:Green}, this yields
natural $\KK$\brd{}equivalences
\begin{equation}  \label{eq:induce_rcross}
  G\Lrcross \Ind_H^G A\cong H\Lrcross A,
  \qquad
  G\Obscross\Ind_H^G A\cong H\Obscross A.
\end{equation}
Hence the (strong) Baum-Connes conjectures for $G\rcross\Ind_H^G A$ and
$H\rcross A$ are equivalent.  As a result, the (strong) Baum-Connes conjecture
with coefficients and the (strong) Baum-Connes conjecture with commutative
coefficients are both hereditary for subgroups.  For the usual Baum-Connes
conjecture, this is due to Jérôme Chabert and Siegfried Echterhoff
(\cite{Chabert-Echterhoff:Permanence}).

\subsection{Full and reduced crossed products and functoriality}
\label{sec:functoriality_Ktop}

Let $\phi\colon G_1\to G_2$ be a continuous group homomorphism.  It induces a
functor $\phi^*\colon \KK^{G_2}\to\KK^{G_1}$.  Of course, $\phi^*(\pt)=\pt$.
If~$\phi$ is open, then the universal property of full crossed products yields
a natural transformation
\begin{equation}  \label{eq:full_crossed_functor}
  \phi_\ast\colon G_1\cross \phi^\ast(A)\to G_2\cross A
\end{equation}
for $A\in\KK^{G_2}$ (if~$\phi$ is not open, we only get a map to the
multiplier algebra of $G_2\cross A$).  There is no analogue
of~\eqref{eq:full_crossed_functor} for reduced crossed products.  For
instance, the homomorphism from~$G$ to the trivial group induces a
homomorphism $\Cred(G)\to\Cred(\{1\})$ if and only if~$G$ is amenable.
Nevertheless, $\Ktop(G)$ has the same functoriality as full crossed products.
We can reprove this easily in our setup.

\begin{theorem}  \label{the:full_versus_reduced}
  The natural map $G\cross A\to G\rcross A$ is a $\KK$\brd{}equivalence for
  $A\in\gen{\CI}$.  Hence $G\Lcross A\cong G\Lrcross A$ (in $\KK$) for any
  $A\in\KK^G$.
\end{theorem}

\begin{proof}
  Since full and reduced crossed products agree for compact groups,
  \eqref{eq:Green} yields that the map $G\cross A\to G\rcross A$ is an
  isomorphism in $\KK$ for $A\in\CI$.  Since both crossed products are
  triangulated functors that commute with direct sums, this extends from $\CI$
  to $\gen{\CI}$.  This implies the second statement because the localisations
  only see $\gen{\CI}$.
\end{proof}

\begin{corollary}  \label{cor:Ktop_functorial}
  There exists a natural map $\phi_\ast\colon G_1\Lrcross \phi^\ast(A)\to
  G_2\Lrcross A$ for any open, continuous group homomorphism $\phi\colon
  G_1\to G_2$ and any $A\in\KK^{G_2}$.
\end{corollary}

\begin{proof}
  Let $\tilde{A}\to A$ be a $\CI$\brd{}simplicial approximation in
  $\KK^{G_2}$, so that $G_2\cross \tilde{A}\cong G_2\Lcross A$.  Since~$\phi$
  maps compact subgroups in~$G_1$ to compact subgroups in~$G_2$, the functor
  $\phi^\ast\colon \KK^{G_2}\to\KK^{G_1}$ preserves weak equivalences.  Hence
  $\phi^\ast(\tilde{A})\to\phi^\ast(A)$ is a weak equivalence in $\KK^{G_1}$.
  As such it induces an isomorphism on $\Left F$ for any functor~$F$.
  Theorem~\ref{the:full_versus_reduced} and~\eqref{eq:full_crossed_functor}
  yield canonical maps
  \begin{multline*}
    G_1\Lrcross \phi^\ast(A)
    \cong G_1\Lcross \phi^\ast(A)
    \cong G_1\Lcross \phi^\ast(\tilde{A})
    \\ \to G_1\cross \phi^\ast(\tilde{A})
    \to G_2\cross \tilde{A}
    \cong G_2\Lcross A
    \cong G_2\Lrcross A.
    \qedhere
  \end{multline*}
\end{proof}

\subsection{Unions of open subgroups}
\label{sec:direct_unions}

Let $G=\bigcup G_n$ be a union of a sequence of open subgroups.  For instance,
adelic groups are of this form.  Then $G\rcross A\cong
\varinjlim G_n\rcross A$ for any $A\in\KK^G$.  Since restriction to
$G_n\subseteq G$ is a completely positive map $G\rcross A\to G_n\rcross A$,
the inductive system $(G_n\rcross A)_{n\in\N}$ is admissible.  Hence we can
replace the direct limit by the homotopy direct limit (see
Section~\ref{sec:homotopy_limits}).

Let $\Sigma\AN\to\ADir\overset{\Dirac}\to\pt\to\AN$ be a Dirac triangle
for~$G$.  By Proposition~\ref{pro:res_ind_exact}, the functor $\Res_G^{G_n}$
maps this to a Dirac triangle in $\KK^{G_n}$.  Hence
$$
G_n \rcross (\ADir \otimes A) \cong G_n\Lrcross A,
\qquad
G_n \rcross (\AN \otimes A) \cong G_n\Obscross A.
$$
Taking limits, we obtain
\begin{equation}  \label{eq:union}
  G \Lrcross A \cong \hoinjlim G_n \Lrcross A,
  \qquad
  G \Obscross A \cong \hoinjlim G_n \Obscross A.
\end{equation}
We have omitted restriction functors from our notation to avoid clutter.  The
following result is due to Paul Baum, Stephen Millington and Roger Plymen
(\cite{Baum-Millington-Plymen}) for the usual Baum-Connes conjecture.

\begin{theorem}  \label{the:union}
  If the groups~$G_n$ satisfy the (strong) Baum-Connes conjecture with
  coefficients~$A$ for all $n\in\N$, then so does~$G$.
\end{theorem}

\begin{proof}
  Recall that~$G$ satisfies the Baum-Conjecture (or the strong Baum-Connes
  conjecture) with coefficients~$A$ if and only if $G\Obscross A$ is
  $\K$\nbd{}contractible (or $\KK$\brd{}contractible).  Since the category of
  $\K$\nbd{}contractible $C^*$\nbd{}algebras is localising, it is closed under
  homotopy direct limits.  Hence the assertions follow from~\eqref{eq:union}.
\end{proof}

\subsection{Direct products of groups}
\label{sec:product_groups}

Let $G_1$ and~$G_2$ be locally compact groups and let $G\defeq G_1\times G_2$.
Let $\Dirac_j\in\KK^{G_j}(\ADir_j,\pt)$ be Dirac morphisms for the factors.
Then $\Dirac_1\otimes\Dirac_2\in \KK^{G_1\times G_2}
(\ADir_1\otimes\ADir_2,\pt)$ is a Dirac morphism for $G_1\times G_2$ because
$$
\CI(G_1)\otimes\CI(G_2)\subseteq\CI(G_1\times G_2)
\qquad\text{and}\qquad
\CC(G_1)\otimes\CC(G_2)\subseteq\CC(G_1\times G_2).
$$
Let $A_j\in\KK^{G_j}$ for $j=1,2$ and put $A\defeq A_1\otimes A_2\in\KK^G$.
We have a natural isomorphism
$$
G\rcross A \approx (G_1\rcross A_1)\otimes (G_2\rcross A_2)
$$
(because we use minimal $C^*$-tensor products) and hence
\begin{multline*}
  G\Lrcross A
  \cong G\rcross (A_1\otimes\ADir_1)\otimes(A_2\otimes\ADir_2)
  \\ \approx (G_1\rcross A_1\otimes\ADir_1)\otimes
  (G_2\rcross A_2\otimes\ADir_2)
  \cong (G_1\Lrcross A_1)\otimes (G_2\Lrcross A_2).
\end{multline*}
Furthermore, the assembly map $G\Lrcross A\to G\rcross A$ is the exterior
tensor product of the assembly maps $G_j\Lrcross A_j\to G_j\rcross A_j$ for
the factors.

There are, of course, similar isomorphisms for $G\Obscross A$.  Therefore, if
the strong Baum-Connes conjecture holds for $G_1\rcross A_1$ and $G_2\rcross
A_2$, then also for $G\rcross A$.  The corresponding assertion about the usual
Baum-Connes conjecture needs further hypotheses
(see~\cite{Chabert-Echterhoff-Oyono:Going_down}) because we cannot always
compute the $\K$\nbd{}theory of a tensor product by the Künneth Formula.  We
can formulate this as
$$
(G_1\Obscross A_1)\otimes (G_2\Obscross A_2)
\cong (G_1\Obscross A_1)\Lotimes (G_2\Obscross A_2),
$$
using the localised tensor product~$\Lotimes$ introduced in
Section~\ref{sec:UCT}.

Combining the results on finite direct products and unions of groups, we get
assertions about restricted direct products as
in~\cite{Chabert-Echterhoff-Oyono:Going_down}.

\subsection{Group extensions}
\label{sec:group_extensions}

Next we consider a group extension $N\into G\prto G/N$.  If~$A$ is a
$G$-$C^*$\brd{}algebra, then $N\rcross A$ carries a canonical twisted action
of $G/N$.  In~\cite{Chabert-Echterhoff:Permanence}, Jérôme Chabert and
Siegfried Echterhoff use this to construct a \emph{partial crossed product
  functor}
$$
N\rcross\blank\colon \KK^G\to\KK^{G/N}.
$$
This functor is triangulated and commutes with direct sums.  We have a natural
isomorphism $G/N\rcross (N\rcross A)\cong G\rcross A$ in $\KK$.  The following
result is due to Jérôme Chabert, Siegfried Echterhoff and Hervé
Oyono-Oyono~\cite{Chabert-Echterhoff-Oyono:Going_down} for the usual
Baum-Connes conjecture.

\begin{theorem}  \label{the:extension}
  The functor $N\rcross\blank\colon \KK^G\to\KK^{G/N}$ maps $\CI$ to $\CI$ and
  hence $\gen{\CI}$ to $\gen{\CI}$.  Therefore, there is a natural isomorphism
  $$
  G/N\Lrcross (N\Lrcross A)\cong G\Lrcross A,
  $$
  which is compatible with the isomorphism $G/N\rcross (N\rcross A)\cong
  G\rcross A$.
  
  Suppose that the (strong) Baum-Connes conjecture holds for $HN\subseteq G$
  with coefficients~$A$ for any smooth compact subgroup $H\subseteq G/N$.
  Then the (strong) Baum-Connes conjecture holds for~$G$ with coefficients~$A$
  if and only if it holds for $G/N$ with coefficients $N\rcross A$.

  Suppose that $G/N$ and $HN$ for compact subgroups $H\subseteq G/N$ have a
  dual Dirac morphism and satisfy $\gamma=1$.  Then the same holds for~$G$.
\end{theorem}

\begin{proof}
  Let~$A$ be compactly induced from, say, the compact subgroup $H\subseteq G$.
  By~\eqref{eq:resind_as_induced}, this means that~$A$ is a $G\cross
  G/H$-$C^*$\brd{}algebra.  We still have a canonical homomorphism from
  $C_0(G/HN)$ to the central multiplier algebra of $N\rcross A$.  This means
  that $N\rcross A$ is compactly induced as a $G/N$\brd{}algebra.  Therefore,
  $N\rcross\blank$ preserves $\CI$ and hence $\gen{\CI}$.  This implies
  $G/N\Lrcross (N\Lrcross A)\cong G\Lrcross A$.
  
  Proposition~\ref{pro:res_ind_exact} implies that a Dirac morphism for $HN$
  is one for~$N$ as well.  Hence the hypothesis of the second paragraph is
  equivalent to the condition that the assembly map $N\Lrcross A\to N\rcross
  A$ in $\KK^{G/N}$ induces a $\K$\nbd{}equivalence (or a
  $\KK$\brd{}equivalence) $H\cross (N\Lrcross A)\to H\cross (N\rcross A)$ for
  all smooth compact subgroups $H\subseteq G/N$.  By
  Theorem~\ref{the:Lrcross_gen}, the map $G/N\Lrcross (N\Lrcross A)\to
  G/N\Lrcross (N\rcross A)$ is a $\K$\nbd{}equivalence (or a
  $\KK$\brd{}equivalence) as well.  Together with $G\Lrcross A\cong
  G/N\Lrcross (N\Lrcross A)$ this yields the assertions in the second
  paragraph.

  Now assume that $G/N$ and the subgroups $HN\subseteq G$ for $H\subseteq G/N$
  compact have dual Dirac morphisms and satisfy $\gamma=1$.  We show that~$G$
  has the same properties.  Recall that this is equivalent to
  $\gen{\CI(G)}=\KK^G$.  The group homomorphism $\pi\colon G\to G/N$ induces a
  triangulated functor commuting with direct sums $\pi^*\colon
  \KK^{G/N}\to\KK^G$.  Of course, $\pi^*(\pt)=\pt$.  Since
  $\gen{\CI(G/N)}=\KK^{G/N}$, the essential range of~$\pi^*$ is generated by
  objects of the form $\pi^*(\Ind_H^{G/N} A)$, where $H\subseteq G/N$ is
  compact.  We have $\pi^*(\Ind_H^{G/N} A) \cong \Ind_{HN}^G \pi_H^*(A)$,
  where $\pi_H\colon HN\to H$ is the restriction of~$\pi$.  Hence
  $\pt\in\KK^G$ belongs to the localising subcategory generated by the ranges
  of the functors $\Ind_{HN}^G$ for compact subgroups $H\subseteq G/N$.

  By hypothesis, $\KK^{HN}=\gen{\CI(HN)}$.  Since induction is a triangulated
  functor that commutes with direct sums, the range of $\Ind_{HN}^G$ is
  contained in the localising subcategory of $\KK^G$ generated by objects of
  the form $\Ind_{HN}^G \Ind_L^{HN}(D) \cong \Ind_L^G(D)$ for compact
  subgroups $L\subseteq HN$ and $D\in\KK^L$.  As a result, $\pt\in\KK^G$
  belongs to $\gen{\CI(G)}$.  This implies that the Dirac morphism is
  invertible, that is, $G$ has a dual Dirac morphism and $\gamma=1$.
\end{proof}

\subsection{Real versus complex assembly maps}
\label{sec:real_complex}

Now we reprove a result of Paul Baum and Max Karoubi (\cite{Baum-Karoubi}) and
Thomas Schick (\cite{Schick:Real}).  In order to compare the real and complex
assembly maps, we have to distinguish between the real and complex Kasparov
theories in our notation.  We denote them by $\KK^{G\cross X}_\R$ and
$\KK^{G\cross X}_\C$, respectively.  We write $A\mapsto A_\C$ for the
complexification functor $\KK^{G\cross X}_\R\to\KK^{G\cross X}_\C$.  This
functor is obviously triangulated and commutes with direct sums and tensor
products, that is, $(A\otimes_{(X)} B)_\C\cong A_\C \otimes_{(X)} B_\C$.

\begin{proposition}  \label{pro:real_complex}
  The complexification functor $\KK^{G\cross X}_\R\to\KK^{G\cross X}_\C$
  preserves weak contractibility and weak equivalences, and it maps
  $\gen{\CI}$ to $\gen{\CI}$.  Hence it maps a Dirac triangle in $\KK^{G\cross
    X}_\R$ to one in $\KK^{G\cross X}_\C$.
\end{proposition}

\begin{proof}
  Since complexification commutes with restriction and induction, it maps
  $\CC_\R$ to~$\CC_\C$ and~$\CI_\R$ to~$\CI_\C$.  Being triangulated and
  compatible with direct sums, it also maps $\gen{\CI_\R}$ to $\gen{\CI_\C}$.
  This implies the assertion about Dirac triangles.
\end{proof}

There is a long exact sequence that relates real and complex $\K$\nbd{}theory.
This exact sequence is generalised in~\cite{Schick:Real} to a similar long
exact sequences
\begin{multline}  \label{eq:compare_real_complex}
  \dotso
  \overset{\delta}\to \KK^{G\cross X}_{\R,q-1}(A,B)
  \overset{\chi  }\to \KK^{G\cross X}_{\R,q}(A,B)
  \overset{c     }\to \KK^{G\cross X}_{\C,q}(A_\C,B_\C)
  \\
  \overset{\delta}\to \KK^{G\cross X}_{\R,q-2}(A,B)
  \overset{\chi  }\to \KK^{G\cross X}_{\R,q-1}(A,B)
  \overset{c     }\to \KK^{G\cross X}_{\C,q-1}(A_\C,B_\C)
  \overset{\delta}\to \dotso,
\end{multline}
for any $A,B\in\KK^{G\cross X}_\R$.  The map~$c$ is the complexification
functor, $\chi$ is the product with the generator of $\KK_1(\R,\R)\cong\Ztwo$
and~$\delta$ is the composition of the inverse of the Bott periodicity
isomorphism with ``forgetting the complex structure''.  In~\cite{Schick:Real},
\eqref{eq:compare_real_complex} is only written down for $\KK^G$.  The same
proof works for $\KK^{G\cross X}$, even for equivariant Kasparov theory for
groupoids.  It is easy to see that~\eqref{eq:compare_real_complex} is natural
with respect to morphisms in $\KK^{G\cross X}_\R$ (see~\cite{Schick:Real}).
Hence the maps
$$
\KK^{G\cross X}_{\C,q}(A_\C,B_\C)\to \KK^{G\cross X}_{\C,q}(A'_\C,B'_\C)
$$
induced by elements of $\KK^G_\R(A',A)$ and $\KK^G_\R(B,B')$ are
isomorphisms for all $q\in\Z$ once the corresponding maps
$$
\KK^{G\cross X}_{\R,q}(A,B)\to \KK^{G\cross X}_{\R,q}(A',B')
$$
are isomorphisms for all $q\in\Z$.  Remarkably, the converse also holds by
\cite{Schick:Real}*{Lemma 3.1}.  A special case is Karoubi's result that
$\K_*(A)\cong0$ if and only if $\K_*(A_\C)\cong0$ (\cite{Karoubi:Descent}).
Moreover, $A\cong0$ in $\KK_\R^{G\cross X}$ if and only if $A_\C\cong0$ in
$\KK_\C^{G\cross X}$ (because $A\cong0$ if and only if~$0$ induces the
identity map on $\KK^{G\cross X}_*(A,A)$).

\begin{theorem}  \label{the:real_complex_BC}
  Let $A\in\KK^G_\R$.  The (strong) Baum-Connes conjecture for~$G$ holds with
  coefficients~$A$ if and only if it holds with coefficients~$A_\C$.
\end{theorem}

\begin{proof}
  The (strong) Baum-Connes conjecture with coefficients~$A$ is equivalent to
  the statement that $\K_*(G\Obscross A)\cong0$ (or $G\Obscross A\cong0$ in
  $\KK$).  Proposition~\ref{pro:real_complex} implies $G\Obscross A_\C \cong
  (G\Obscross A)_\C$.  Hence the assertion follows from the results
  of~\cite{Schick:Real} discussed above.
\end{proof}

\begin{theorem}  \label{the:dual_Dirac_real}
  Let~$G$ be a locally compact group and~$X$ a locally compact $G$\nbd{}space.
  If there is a dual Dirac morphism in $\KK^{G\cross X}_\C$, then there is one
  in $\KK^{G\cross X}_\R$, and vice versa.  In this case, we have
  $\gamma_\C=1$ if and only if $\gamma_\R=1$.
\end{theorem}

\begin{proof}
  By Theorem~\ref{the:dual_Dirac_splitting}, a dual Dirac morphism exists if
  and only if~$\Dirac$ induces an isomorphism $\KK^{G\cross
    X}_*(C_0(X),\ADir)\cong \KK^{G\cross X}_*(\ADir,\ADir)$.  This holds both
  in the real and complex case.  By Proposition~\ref{pro:real_complex}, the
  Dirac morphism in $\KK^{G\cross X}_\C$ is the complexification of the Dirac
  morphism in $\KK^{G\cross X}_\R$.  Hence the existence of a dual Dirac
  morphism in $\KK^{G\cross X}_\C$ and $\KK^{G\cross X}_\R$ are equivalent by
  the results of~\cite{Schick:Real} discussed above.  Since the
  complexification of a dual Dirac morphism in $\KK^{G\cross X}_\R$ is one in
  $\KK^{G\cross X}_\C$, $\gamma_\C$ is the complexification of~$\gamma_\R$.
  We have $\gamma=1$ if and only if~$\gamma$ is invertible if and only if
  multiplication by~$\gamma$ on $\KK^{G\cross X}_*(C_0(X),C_0(X))$ is an
  isomorphism.  Again it follows from~\cite{Schick:Real} that $\gamma_\R=1$ if
  and only if $\gamma_\C=1$.
\end{proof}

\begin{appendix}

\section{The equivariant Kasparov category is triangulated}
\label{sec:KK_triangulated}

We have defined a translation automorphism and a class of exact triangles on
$\widetilde{\KK}{}^{G\cross X}$ in Section~\ref{sec:preliminaries}.  Here we
prove that these data verify the axioms of a triangulated category
(see~\cite{Neeman:Triangulated}).  More precisely, we prove the equivalent
assertion that the opposite category of $\widetilde{\KK}{}^{G\cross X}$ is
triangulated.

By definition, the class of exact triangles is closed under isomorphism and
the translation functor is an automorphism.  The zeroth axiom (TR~0) requires
triangles of the form $\Sigma X\to 0\to X\overset{\ID_X}\to X$ to be exact.
This follows from the contractibility of $\cone(\ID_X)\cong
C_0(\mathopen]0,1])\otimes X$.

Axiom (TR~1) asks that for any morphism $f\colon A\to B$ there should be an
exact triangle $\Sigma B\to C\to A\overset{f}\to B$.  If~$f$ is an equivariant
$*$\nbd{}homomorphism, we may take the mapping cone triangle of~$f$.  In
general, we claim that any morphism in $\widetilde{\KK}{}^G$ is isomorphic to
an equivariant $*$\nbd{}homomorphism.  We can first replace~$f$ by a morphism
in $\KK^G$ because $\KK^G$ and $\widetilde{\KK}{}^G$ are equivalent
categories.  By~\cite{Meyer:KKG} we can represent~$f$ by an equivariant
$*$\nbd{}homomorphism $f_*\colon q_s A\to q_s B$, where
$$
q_s A\defeq  \Comp(L^2(G\times\N)) \otimes q(\Comp(L^2G)\otimes A).
$$
If $X=\pt$, then the $C^*$\nbd{}algebra $qA$ is the usual one from the
Cuntz picture for Kasparov theory.  Otherwise, we have to modify its
definition so as to get a $G\cross X$-$C^*$-algebra.  Namely, let $A\ast_X A$
be the free product of~$A$ with itself in the category of $G\cross
X$-$C^*$-algebras.  That is, it comes equipped with two natural maps
$\iota_1,\iota_2\colon A\to A\ast_X A$ with the universal property that pairs
of $G\cross X$\brd{}equivariant $*$\nbd{}homomorphisms $A\to B$ correspond
bijectively to $G\cross X$\brd{}equivariant $*$\nbd{}homomorphisms $A\ast_X
A\to B$.  We can construct $A\ast_X A$ as the quotient of $A\ast A$ by the
ideal generated by the relations $\iota_1(fa_1)\iota_2(a_2)
\sim\iota_1(a_1)\iota_2(fa_2)$ for all $a_1,a_2\in A$, $f\in C_0(X)$.  The
pair $(\ID_A,\ID_A)$ induces a natural homomorphism $A\ast_X A\to A$.  Let
$q_XA$ be its kernel.  With this modified definition of $qA$, the assertions
of~\cite{Meyer:KKG} remain true for $\KK^{G\cross X}$.  In particular, there
is a natural $\KK^{G\cross X}$\brd{}equivalence $q_s A\cong A$.  Therefore,
any morphism in $\KK$ is isomorphic to an equivariant $*$\nbd{}homomorphism.
Thus axiom (TR~1) holds.

Axiom (TR~2) asks that a triangle  $\Sigma B\to C\to A\to B$ be exact if and
only if $\Sigma A\to\Sigma B\to C\to A$ (with certain signs) is exact.  It
suffices to prove one direction because suspensions and desuspensions
evidently preserve exact triangles.  Thus axiom (TR~2) is equivalent to the
statement that the rotated mapping cone triangle
$$
\Sigma A
\overset{-\Sigma f}\longrightarrow \Sigma B
\overset{\iota}\longrightarrow \cone(f)
\overset{\epsilon}\longrightarrow A
$$
is exact for any equivariant $*$\nbd{}homomorphism $f\colon A\to B$.  We
claim that this triangle is the extension triangle for the natural extension
$$
0
\longrightarrow \Sigma B
\overset{\iota}\longrightarrow \cone(f)
\overset{\epsilon}\longrightarrow A
\longrightarrow 0
$$
and hence exact.  Build the diagram~\eqref{eq:extension_cone} for this
extension.  The resulting map $\Sigma B\to\cone(\epsilon)$ is a homotopy
equivalence in a natural and hence equivariant fashion.  Thus the above
extension is admissible and gives rise to an exact triangle.  One easily
identifies the map $\Sigma A\to\Sigma B$ in the extension triangle with
$-\Sigma f$.  This finishes the proof of axiom (TR~2).

Suppose that we are given the solid arrows in the diagram
\begin{equation}  \label{eq:TR_three}
\begin{gathered}
  \xymatrix{
    {\Sigma B} \ar[d]^{\Sigma \beta} \ar[r] &
    {C} \ar[r] \ar@{.>}[d]^{\gamma} &
    {A} \ar[r] \ar[d]^{\alpha} &
    {B} \ar[d]^{\beta} \\
    {\Sigma B'} \ar[r] &
    {C'} \ar[r] &
    {A'} \ar[r] &
    {B'}
  }
\end{gathered}
\end{equation}
and that the rows in this diagram are exact triangles.  Axiom (TR~3)
asks that we can find~$\gamma$ making the diagram commute.  We may first
assume that the rows are mapping cone triangles for certain maps $f\colon A\to
B$ and $f'\colon A'\to B'$ because any exact triangle is isomorphic to
one of this form.

We represent $\alpha$ and~$\beta$ by Kasparov cycles, which we again denote by
$\alpha$ and~$\beta$.  Since~\eqref{eq:TR_three} commutes, the Kasparov cycles
$f'_*(\alpha)$ and $f^*(\beta)$ are homotopic.  Choose a homotopy~$H$ between
them.  Now we glue together $\beta$, $H$ and~$\alpha$ to obtain a cycle for
$\KK^{G\cross X}(\cone(f),\cone(f'))$ with the required properties.  Since
$(\ev_1)_*(H)=f'_*(\alpha)$, the pair $(H,\alpha)$ defines a Kasparov cycle
for~$A$ and $\cyl(f')$.  The constant family~$\beta$ defines a cycle~$C\beta$
for $\KK^{G\cross X} (C_0(\mathopen]0,1],B), C_0(\mathopen]0,1],B'))$.
Reparametrisation gives a canonical isomorphism
$$
\cone(f') \cong
\{(x,y)\in C_0(\mathopen]0,1],B') \oplus \cyl(f') \mid
x(1)=\tilde{f}'(y)\}.
$$
Since $\tilde{f}'_*(H,\alpha)=(\ev_0)_*(H)=f^*(\beta)$, we can glue
together $(H,\alpha)$ and $C\beta$ to get a cycle for $\KK^{G\cross
  X}(\cone(f),\cone(f'))$.  It is straightforward to see that it has the
required properties.  This finishes the verification of axiom (TR~3).

It remains to verify Jean-Louis Verdier's octahedral axiom, which is crucial
to localise triangulated categories.  Amnon Neeman formulates it rather
differently in~\cite{Neeman:Triangulated}.  We shall use Verdier's original
octahedral axiom (see \cite{Verdier:Thesis} or
\cite{Neeman:Triangulated}*{Proposition 1.4.6}) because it can be applied more
directly and because its meaning is more transparent in the applications we
have met so far.

\begin{proposition}  \label{pro:TRfour}
  For any pair of morphisms $f\in\KK^G(B,D)$, $g\in\KK^G(A,B)$ there is a
  commuting diagram as in Figure~\ref{fig:TRfour}
  \begin{figure}
    $$
    \xymatrix{
      \Sigma^2 D \ar[d] \ar[r] & \Sigma C_f \ar[d] \ar[r] &
      \Sigma B \ar[d] \ar[r]^{\Sigma f} & \Sigma D \ar[d] \\
      0 \ar[r] \ar[d] & C_g \ar@{=}[r] \ar[d] &
      C_g \ar[r] \ar[d] & 0 \ar[d] \\
      \Sigma D \ar[r] \ar[d] & C_{fg} \ar[r] \ar[d] &
      A \ar[r]^{fg} \ar[d]^{g} & D \ar@{=}[d] \\
      \Sigma D \ar[r] & C_f \ar[r] &
      B \ar[r]^{f} & D
    }
    $$
    \caption{The octahedral axiom}
    \label{fig:TRfour}
  \end{figure}
  whose rows and columns are exact triangles.  Moreover, the two maps $\Sigma
  B\to\Sigma D\to C_{fg}$ and $\Sigma B\to C_g\to C_{fg}$ in this diagram
  coincide.
\end{proposition}

\begin{proof}
  Replacing all $C^*$\nbd{}algebras by appropriate universal algebras, we can
  achieve that $f$ and~$g$ are equivariant $*$\nbd{}homomorphisms.  We assume
  this in the following.  We shall use the mapping cones and mapping cylinders
  defined in Section~\ref{sec:triangulated_categories}.  We define a natural
  $G$\nbd{}$C^*$\brd{}algebra
  $$
  \cyl(f,g)
  \defeq \bigl\{(a,b,d)\in A\oplus C([0,1],B)\oplus C([0,1],D) \mid
  g(a)=b(1),\ f\bigl(b(0)\bigr)=d(1)\bigr\}
  $$
  and natural equivariant $*$\nbd{}homomorphisms
  \begin{alignat*}{2}
    p_A &\colon \cyl(f,g)\to A,
    &\qquad (a,b,d) &\mapsto a,
    \\
    j_A &\colon A\to\cyl(f,g),
    &\qquad a &\mapsto \bigl(a,\const g(a), \const fg(a)\bigr),
    \\
    \tilde{g} &\colon \cyl(f,g) \to \cyl(f),
    &\qquad  (a,b,d) &\mapsto (b(0),d).
  \end{alignat*}
  We have $p_Aj_A=\ID_A$, and $j_Ap_A$ is homotopic to the identity map in a
  natural and hence equivariant way.  Thus $\cyl(f,g)$ is homotopy equivalent
  to~$A$.  Moreover, $\tilde{g}j_A=j_B g$, where $j_B\colon B\to\cyl(f)$ is
  the natural map, which is a homotopy equivalence.  That is, the
  map~$\tilde{g}$ is isomorphic to $g\colon A\to B$.  Recall also that the map
  $\tilde{f}\colon \cyl(f)\to D$ is isomorphic to $f\colon B\to D$.
  
  The maps $\tilde{g}\colon \cyl(f,g)\to\cyl(f)$, $\tilde{f}\colon \cyl(f)\to
  D$ and $\tilde{f}\circ\tilde{g}\colon \cyl(f,g)\to D$ are all surjective.
  The kernel of~$\tilde{f}$ is $\cone(f)$, the kernel of~$\tilde{g}$ is
  naturally isomorphic to $\cone(g)$.  We let $\cone(f,g)$ be the kernel of
  $\tilde{f}\tilde{g}$.  Thus we obtain a commutative diagram of
  $G$\nbd{}$C^*$\brd{}algebras whose rows and columns are extensions:
  \begin{equation}  \label{fig:octahedral}
    \begin{gathered}
    \xymatrix@C=1.3cm{
      \cone(g) \ar@{=}[r] \ar[d] &
      \cone(g) \ar[r] \ar[d] &
      0 \ar[d] \\
      \cone(f,g)\ \ar@{>->}[r] \ar[d] &
      \cyl(f,g) \ar@{->>}[r]^{\tilde{f}\tilde{g}} \ar[d]^{\tilde{g}} &
      D \ar@{=}[d] \\
      \cone(f)\ \ar@{>->}[r] &
      \cyl(f) \ar@{->>}[r]^{\tilde{f}} &
      D.
    }
    \end{gathered}
  \end{equation}
  We claim that all rows and columns in this diagram are admissible
  extensions.  (Even more, the maps $K\to\cone(p)$
  in~\eqref{eq:extension_cone} for these extensions are all homotopy
  equivalences.)  We have already observed this for the third row in
  Section~\ref{sec:extension_triangles}, and the argument for the second row
  is similar.  The assertion is trivial for the first row and the third
  column.  The remaining two columns can be treated in a similar fashion.  A
  conceptual reason for this is that they are pull backs of the standard
  extension $\cone(g)\into\cyl(g)\prto B$ along the natural projections
  $\cone(f)\to B$ and $\cyl(f)\to B$, respectively.  The projection
  $\cyl(g)\to B$ is a cofibration in the notation
  of~\cite{Schochet:Axiomatic}; this property implies that $K\to\cone(p)$ is a
  homotopy equivalence and is hereditary for pull backs
  (see~\cite{Schochet:Axiomatic}).
  
  We can now write down extension triangles for the rows and columns
  in~\eqref{fig:octahedral} and replace $A$ and~$B$ by the homotopy equivalent
  algebras $\cyl(f,g)$ and $\cyl(f)$, respectively.  This yields a diagram as
  in Figure~\ref{fig:TRfour}.
  
  The composite map $\Sigma B\to \cone(g)\to \cone(f,g)$ is just the
  restriction of the canonical map $\cone(g)\to\cone(f,g)$ to $\Sigma B$.
  There is a natural homotopy from this map to the composition $\Sigma
  B\to\Sigma D\to\cone(f,g)$ via translations involving~$f$.  This finishes
  the proof of Proposition~\ref{pro:TRfour}.
\end{proof}

We have now verified that $\widetilde{\KK}^{G\cross X}$ is a triangulated
category.

\end{appendix}

\end{document}